\documentclass[11pt]{amsart}
\usepackage{amsfonts} 
\usepackage{amscd}
\usepackage{amssymb}
\usepackage{float}
\usepackage{amsmath}
\usepackage{graphicx}
\usepackage{amsxtra}
\usepackage{amsthm}

\newcommand{\ring}[1]{\mathbb{#1}}
\newcommand{\C}{\ring{C}} \newcommand{\Q}{\ring{Q}}
\newcommand{\Z}{\ring{Z}} 
\newcommand{\A}{\ring{A}} \newcommand{\G}{\ring{G}} 
\newcommand{\F}{\ring{F}} \newcommand{\hi}{\ring{H}}
\newcommand{\bu}{\bullet}

\newcommand{\ssp}{\smallskip}

\newcommand{\ra}{\rightarrow}
\newcommand{\ub}{\underline}
\newcommand{\hra}{\hookrightarrow}
\newcommand{\zim}{\xrightarrow{\sim}}

\newcommand{\gm}{\ring{G}_m}
\newcommand{\be}{\begin{equation}}
\newcommand{\ee}{\end{equation}}
\newcommand{\nd}{\noindent}
\newcommand{\End}{\operatorname{End}}
\newcommand{\Ker}{\operatorname{Ker}}
\newcommand{\im}{\operatorname{Im}}
\newcommand{\Coker}{\operatorname{Coker}}
\def\Hom{\textup{Hom}}
\def\1{{\mu\mkern-6mu\mu}}
\title[motives]{One-motives and a conjecture of Deligne}
\author{Niranjan Ramachandran}
\address{Department of Mathematics, University of Maryland, College Park
  MD 20742, USA}
 \address{Max-Planck-Institut f\"ur Mathematik\\
    Vivatsgasse 7\\ D-53111 Bonn, Germany}
\email{atma@math.umd.edu}
\date{\today}

\begin{document}
\theoremstyle{plain}
\newtheorem{thm}{Theorem}[section]
\newtheorem{lem}[thm]{Lemma}
\newtheorem{cor}[thm]{Corollary}
\newtheorem{prop}[thm]{Proposition}
\newtheorem{conj}[thm]{Conjecture}
\newtheorem{quest}[thm]{Question}
 
\theoremstyle{definition}
\newtheorem{rem}[thm]{Remark}
\newtheorem{defn}[thm]{Definition}
\newtheorem{ex}[thm]{Example}
\newtheorem{obser}[thm]{Observation}
 
\begin{abstract}
We introduce new motivic invariants of arbitrary varieties over a perfect
field. These cohomological invariants 
take values in the category of one-motives (considered up to isogeny
in positive characteristic). The
algebraic definition of 
these invariants  
presented here proves a conjecture of
Deligne. Other applications include some cases of conjectures of Serre, Katz,
and Jannsen on the independence of $\ell$ of
parts of the \'etale cohomology of arbitrary varieties  over number
fields and finite fields. 
\end{abstract}
\thanks{Funded in part by grants from Hewlett-Packard, Graduate Research Board
(U. Maryland) and MPIM (Bonn)} 
\maketitle
\setcounter{tocdepth}{1}
\tableofcontents

\hfill $\frak{Wenn~die~K\ddot{o}nige~bau'n,}$

\hfill $\frak{haben~die~K\ddot{a}rrner~zu~tun.}$

\hfill $\frak{F.~Schiller}$
\medskip


\nd {\bf Introduction.} P. Deligne \cite[10.4.1]{h} has attached one-motives to complex
 algebraic varieties 
 using the theory of mixed Hodge structures.  He has conjectured  that
these one-motives admit a \emph{purely algebraic} definition.
The aim of  this article is to prove his  conjecture (Theorem
 \ref{peddha}).

Recall the well known result of Riemann \cite[4.4.3]{h2}, presented
here  in modern
guise: the
``Hodge realization" $T_{\Z}$  ---  this is $A \mapsto H_1(A, \Z)$  ---   
defines an equivalence from the 
category of complex abelian varieties to the category of 
torsion-free polarizable Hodge structures of type $\{(0,-1),
	(-1,0)\}$. In particular, any such Hodge structure arises as
the $H_1$ of an essentially unique complex 
abelian variety. 

Deligne \cite[\S 10.1]{h} has introduced the algebraic notion of a
one-motive  over a field $k$, generalizing that of an abelian
variety   ---  \S \ref{rev} contains the precise 
definitions; he has also generalized Riemann's result by showing that the
``Hodge realization'' $T_{\Z}$ 
defines an equivalence  from the 
category of one-motives over $\C$ to the
category of torsion-free mixed Hodge structures $H$ of type 
$$(*) \qquad {} \qquad {} \qquad {} \qquad \{(-1, -1), (-1, 0), (0, -1), (0,0)\}
\qquad {}$$ 
with $Gr^W_{-1}H$ polarizable. Thus, 
any such mixed Hodge structure $H$ arises from an essentially unique
one-motive 
$I(H)$ over $\C$. The functor $I$ is a quasi-inverse to $T_{\Z}$.

For any complex variety $V$ and any integer $n \ge 0$, consider the largest
mixed Hodge substructure  $t^n(V)$ of type $(*)$ of $H^n(V,
\Z(1))/{\rm torsion}$; there exists a well-defined one-motive $I^n(V)$
over $\C$ whose
Hodge realization is $t^n(V)$; so $I^n(V):= I(t^n(V))$. Deligne
\cite[10.4.1]{h} has conjectured that $I^n(V)$
admits a purely algebraic definition. His proof (ibid. 10.3.  --- 
Interpr\'etation alg\'ebrique du $H^1$ mixte: cas des courbes) of his 
conjecture  for arbitrary curves suggests a precise formulation of the
conjecture. Namely, we 
have the following (this formulation  is due to the
referee): 

\begin{conj}\label{dc1} {\rm (Deligne)}
For an arbitrary variety $V$ over an arbitrary field $k$  and integer
$n$, define a one-motive $L^n(V/k)$ and homomorphisms\footnote{Here
$\bar{k}$ is an algebraic closure of $k$.} 
\begin{eqnarray*} T_{\ell}(L^n(V/k)) & \to & H^n(V \times \bar{k}, \Z_{\ell}(1))/{\rm
torsion}, \\ T_{DR}(L^n(V/k)) & \to & H^n_{DR}(V/k)
\end{eqnarray*}
from the $\ell$-adic and de Rham realizations of $L^n(V/k)$. The
definitions of $L^n(V/k)$ and the homomorphisms should be algebraic,
canonical, and 
functorial in $V$ and $k$.

 Furthermore, 
$L^n(V/{\C})$ should be canonically isomorphic to $I^n(V)$.
\end{conj}

(Clearly, $V$ can be replaced by a simplicial scheme.)

The prototype is A. Weil's construction \cite{aweil} of the Jacobian;
his construction proves the conjecture for smooth projective curves and
$n=1$. The conjecture  is true for smooth projective varieties
(\ref{dcon}): it amounts to an algebraic construction of
the Picard variety and the N\'eron-Severi group.

The case $n=1$ of (\ref{dc1}) is known (up to $p$-isogeny in
characteristic $p >0$) for arbitrary varieties over
perfect fields \cite{bs3, h, ra, jp2}; the case $n=2$ is known 
for complex proper surfaces \cite{ca,
ca2}. No general results were  known for higher cohomology (i.e., for
$n >2$).  

A natural approach to Conjecture \ref{dc1} is to use
 proper hypercoverings \cite[6.2]{h} by smooth simplicial schemes; namely, to mimic
 Deligne's approach \cite{h} 
 to the construction of the mixed Hodge structure on $H^*(V,\Z)$ of a
 complex algebraic variety $V$. This approach, which  we follow here, 
 gives a two-step strategy to prove (\ref{dc1}):\ssp

\nd {\bf Step 1.} Construct one-motives $L^n$ ($n \ge 0$) for smooth
 simplicial schemes 
 arising from simplicial pairs (\ref{simpxy}) and show that they have
 the properties given in (\ref{dc1}).\ssp

\nd {\bf Step 2.} Prove cohomological descent for these one-motives; more
precisely, show that the one-motives $L^n$, given by (i), of a proper
hypercovering of a variety $V$ are ``independent'' of the proper
hypercovering; and, thus, $L^n$ depend only on $V$.\ssp

Sections \ref{constr}, \ref{Tmt}, \ref{oakland} are devoted to the
first step, but only for fields of characteristic zero; the case of positive
characteristic is relegated to Section \ref{pos+}. Our construction of
the requisite one-motives $L^n$, inspired by \cite{ca}, relies on the theory
of the Picard scheme 
\cite[Chapter 8]{blr}; the techniques are
those of \cite{ra} but here applied to truncated
simplicial schemes. The realizations of $L^n$ are treated in Sections
\ref{Tmt} (Hodge, de Rham), \ref{oakland} (\'etale); here a crucial
use is made of the validity of the Hodge conjecture for divisors
(\ref{h11}).     

Section \ref{varieties} is devoted to the second step.
It turns out that, because an important 
spectral sequence \cite[8.1.19.1]{h} degenerates only with rational
coefficients, the method of proper hypercoverings only provides 
a theory of isogeny one-motives $L^*(-)\otimes\Q$. More precisely,
given two proper hypercoverings $U_{\bu}$ and $'U_{\bu}$ 
of $V$, we can only show that the
associated one-motives $L^n$ and $'L^n$ are isogenous; the isogeny
one-motive $L^n\otimes\Q$ depends only on $V$.  Thus, a new
ingredient is necessary to complete the second step, i.e., 
to endow these isogeny one-motives with integral structures. This is
done, as in \cite{milram}, via the integral structure on \'etale
cohomology. Thus, we provide a complete proof (\ref{peddha}) of Conjecture
\ref{dc1} for an arbitrary field of characteristic zero.

We now turn to the case of Conjecture
\ref{dc1} for a  field $k$ of characteristic $p >0$; let us begin by
indicating why the conjecture must be weakened slightly.

First, in \cite[Appendix]{jamo}, A. Grothendieck notes that,
 for a curve $C$ 
over $k$, the construction of Deligne \cite[10.3]{h}
provides a one-motive $H^1_m(C) = L^1(C/k)$ defined over the perfection
$k^{perf}$ of $k$; thus, \cite[10.3]{h} proves the case $n=1$ of (\ref{dc1})  
only for curves over a perfect field. Second, he (loc. cit) expresses
 doubts  about the
 existence of a $\Z$-linear (i.e., integral) category of mixed motives
 over an imperfect field $k$; he anticipates only a
 $\Z[1/p]$-linear category,  i.e., a category of mixed motives up to
 $p$-isogeny.  
Third, the existence of proper hypercoverings (by 
smooth simplicial schemes) for an arbitrary variety $V$ over $k$ 
is known only when $k$ is perfect \cite{dj}. If one expects that the
 method of proper 
hypercoverings provides, as in characteristic zero, one-motives up to
isogeny associated with $V$, then the methods of
\cite{milram} allow a refinement to one-motives defined up to
$p$-isogeny: controlling $p$-isogeny requires an integral $p$-adic 
cohomology theory for arbitrary varieties over $k$. These
 considerations\footnote{V. Voevodsky \cite{vv} works over perfect
 fields and neglects $p$-torsion in characteristic $p >0$.} 
lead us to a weak version  of (\ref{dc1}) by only 
requiring  one-motives $L^*(V/k)\otimes\Z[1/p]$
(up to $p$-isogeny) over $k^{perf}$.

While the first step  can be carried out in positive characteristic in
the same way as in characteristic zero, the second step cannot
be unless, as it seems, one assumes the Tate conjecture
(\ref{tatec}) for divisors. Note that our proof of (\ref{dc1}) 
in characteristic zero depends on the validity of the Hodge
conjecture (\ref{h11}) for divisors; thus, the appearance of the Tate
conjecture is rather natural. 

For a perfect field $k$, (\ref{dc1})  ---  up to $p$-isogeny  ---  holds (\ref{groth})
under the assumption of the Tate conjecture
(\ref{tatec}) for  
surfaces. The analogous result  
is also valid for an imperfect field $k$ under the additional assumption of
``resolution of singularities'' over $k$.

Using a suggestion of M. Marcolli, we provide an
\emph{unconditional} construction  (\ref{mem}) of 
$J^n(-)\otimes\Z[1/p]$ ($n \ge 0$) and $L^n(-)\otimes\Z[1/p]$ ($ 2\ge
n \ge 0$) of one-motives (up to $p$-isogeny) for arbitrary
varieties over a perfect field $k$. The $J^n(V)\otimes\Z[1/p]$ are
good substitutes for the (as yet conditional) $L^n(V)\otimes\Z[1/p]$; for
instance, $W_{-1}J^n(V) = W_{-1}L^n(V)$. 

In particular, we generalize Carlson's results \cite{ca} on $L^2$ of
a complex projective surface to any variety over a perfect field (and
up to $p$-isogeny in characteristic $p>0$). 
 
In Section \ref{applix}, we use these new invariants $L^n$ and
$J^n$ to provide affirmative answers to special cases of 
questions \cite{ja, ka, se2} in the motivic 
folklore. These concern
``independence of $\ell$'' of 
$\ell$-adic \'etale cohomology of \emph{arbitrary} 
varieties over number fields and finite fields. 

Since the circulation of this manuscript (circa 1998), other authors
\cite{brs} have independently obtained some of the results presented
here; \cite{banff} is a leisurely introduction to our results.\medskip

\begin{small} 
\nd {\bf Notation.}
We work over $S:= $~Spec $k$; here $k$ is a field of
characteristic $p$ (except in \S \ref{rev}, \ref{pos+}, $k$ is perfect unless
indicated otherwise);  in sections
\S \ref{constr}, \ref{Tmt}, \ref{oakland}, and \ref{varieties} we assume $p$
to be zero. In \S \ref{pos+}, we assume $p>0$. 

We fix an algebraic closure $\bar{k}$ of $k$; $\bar{S}:=$
Spec~$\bar{k}$ and $\G:= {\rm{Gal}}(\bar{k}/k)$.  
For any scheme $X$ over $S$, we set 
$\bar{X}:= X \times_S \bar{S}$. All schemes will be supposed to be 
separated and locally noetherian. A variety is a geometrically
 integral scheme of
finite type over $S$.

\nd For any set $B$, $\Z(B)$ is the free abelian group
generated by the elements of $B$.

\nd $\pi_0(X):= $ the set (with a $\G$-action) of connected
components of $\bar{X}$.

\nd $D_X:= $ the \'etale group scheme corresponding to $\Z(\pi_0(\bar{X}))$.

\nd $T_X:= Hom (D_X,\gm)$, the algebraic torus associated with  $D_X$.

\nd $w(X):= $ the set of irreducible components of $\bar{X}$. 

\nd $W_X:= $ the \'etale group scheme corresponding to the $\G$-module
$\Z({w(X)})$. 
 
\nd For any group scheme $\mathcal G$, $\pi_0(\mathcal G)$ is a group
with an action of $\G$; we shall also use $\pi_0(\mathcal G)$ for the 
corresponding \'etale group scheme.

For
a variety $V$ over $\C$, $V(\C)$ (resp. $V^{an}$) denotes the
associated topological space with the classical topology 
(resp. analytic variety). Given $X$ over $S$ and an imbedding $\iota: k
\hookrightarrow \C$, we denote by $X_{\iota}$ the scheme over $\C$
obtained by base change.

\nd $MHS:=$ the abelian category of ($\Z$)-mixed Hodge structures.

\nd $\mathcal M_F:=$ the additive category of one-motives over a field
$F$. 

\nd $S_{fppf}$ (resp. $S_{fpqc}$) is the big site over $S$ with the $fppf$ 
(resp. $fpqc$) topology \cite[pp.~200-201]{blr}. 

\nd $\hat{\Z}$ is the profinite
completion ${\varprojlim}_r \Z/{r\Z}$ of $\Z$.

\nd  $\mathbb A =
\hat{\Z}\otimes \Q$ is the ring of finite ad\`eles of $\Q$.

\nd $\Z^p: = {\varprojlim}_r \Z/{r\Z}$ with $r$ coprime to $p$ and
$\A^p:= \Z^p\otimes\Q$.

We refer to section x.y by \S x.y and to specific results by 
Theorem x.y or Remark x.y or simply (x.y). We use $\square$ to denote the end
of a remark or a proof. 
\end{small}

\section{Preliminaries}\label{rev}
Let us begin by reviewing well-known results, which will be of use
in the paper.\ssp 

\nd{\bf One-motives.} \cite[\S 10]{h}\ssp

 A one-motive $M:= [B
\xrightarrow{u} G]$ over $S$ (or over $k$) is a two-term complex  
consisting of a semi-abelian variety $G$ over $k$ (i.e., $G$ is an
extension of an abelian variety $A$ by a torus $T$), a finitely
generated torsion-free abelian 
group $B$ with a structure of a discrete $\G$-module, and a homomorphism $u:B
 \rightarrow G(\bar{k})$ of $\G$-modules.  In particular, if $k$ is
 algebraically closed, then $u$ is a homomorphism of abelian
 groups. It is
convenient to regard $B$ as an \'etale group scheme  (locally
constant)  on $S$.  A morphism of one-motives is a morphism of
complexes. From  the category $M_k$ of one-motives over $k$, there are 
``realization'' functors  (Hodge) $T_{\Z}$  ---  for
each $\iota: k \hra \C$  ---  to
(the category of) torsion-free $\Z$-mixed Hodge structures of type
$(*)$, (\'etale) $T_{\ell}$  ---  for each $\ell \neq p$  ---  to 
$\Z_{\ell}$-modules with an action of $\G$, and (de Rham) $T_{DR}$  --- 
if $p =0$  ---  to $k$-vector spaces. 

A morphism $\phi:M_1 \to M_2$ is called an isogeny if $\phi_B: B_1
\to B_2$ is injective with finite cokernel and $\phi_G: G_1 \to G_2$
is surjective with finite kernel. 
The (additive) category $\mathcal M_k$ of one-motives over
$k$ enjoys Cartier duality \cite[10.2.11]{h}. 
The dual of an isogeny is also an
isogeny. A $p$-isogeny is an isogeny $\phi$ such that the orders of
$\Coker(\phi_B)$ and $\Ker(\phi_G)$ are powers of $p$. 

The $\Q$-linear
abelian category $\mathcal M_k\otimes \Q$ 
of isogeny one-motives over $k$ is obtained from $\mathcal M_k$ by
inverting isogenies;  $\mathcal M_k\otimes \Q$ inherits
realization functors (Hodge) $T_{\Z}$  to $\Q$-mixed Hodge structures,
(\'etale) 
$T_{\ell}$ to $\Q_{\ell}$-vector spaces,  (de
Rham) $T_{DR}$ to 
$k$-vector spaces, weight filtration $W$, and Cartier duality from
$\mathcal M_k$. Every one-motive $M$ defines an isogeny one-motive $M
\otimes\Q$. The weight filtration $W$ on $[B \xrightarrow{u}
G]\otimes\Q$ is $W_{-3} =0$, $W_{-2} = [0 \to T]\otimes\Q$, $W_{-1} = [0 \to
G]\otimes\Q$, and $W_0 = [B \xrightarrow{u}
G]\otimes\Q$.  

Finally, if $p =0$, the functors $T_{\ell}$ can be combined to a
functor $M \mapsto TM = 
M\otimes\hat{\Z}:= {\prod_{\ell}}~T_{\ell}M$; here $TM$ is a
$\hat{\Z}$-module with a $\G$-action. The ad\`elic realization functor $M
\mapsto TM\otimes\Q$ from $\mathcal M_k$ to $\A$-modules with a
$\G$-action factorizes via $\mathcal M_k\otimes\Q$; this gives the
ad\`elic realization functor of an isogeny one-motive: $M\otimes\Q
\mapsto M\otimes\A$.

If $p > 0$, then the $\Z[1/p]$-linear category $\mathcal M_k\otimes
\Z[1/p]$ of one-motives up to $p$-isogeny over $k$ is obtained from   $\mathcal M_k$ by
inverting $p$-isogenies. The realization functor $M \mapsto T^pM  =
{\prod_{\ell \neq p}}~T_{\ell}M$ from $\mathcal M_k$ to the category of $\Z^p$-modules
with a $\G$-action extends to the category $\mathcal M_k\otimes 
\Z[1/p]$.\ssp

\nd {\bf Relative representability.}\ssp

As indicated in 
\cite[pp.200-201]{blr}, representability issues  are best treated
in the $fppf$-topology. The following simple lemma will be used
often. 
\begin{lem}\label{repco}

{\rm{(i)}} Let $F$ be a representable contravariant functor from the category
of schemes over $S$ to sets. Then $F$ is a sheaf with respect to the
${fpqc}$-topology  and, hence, with respect to the $fppf$, \'etale, 
and Zariski topologies.

{\rm{(ii)}} Let $0 \rightarrow F \rightarrow G \rightarrow H \rightarrow
0$ be an exact sequence of sheaves of abelian groups on
$S_{fppf}$. Suppose $F,H$ are representable, and that $F \rightarrow
S$ is an affine morphism {\rm (}i.e. the scheme representing $F$ is affine
over $S${\rm )}. Then $G$ is representable, necessarily by a commutative
group scheme. 
\end{lem}

\begin{proof} (i) \cite[Prop. 1, p.200]{blr}. 

\nd (ii) The proof of \cite[Prop. 17.4]{oort} for $S_{fpqc}$
also works for $S_{fppf}$. \end{proof}

\nd {\bf Picard functor.}\ssp

 Let $f:X \rightarrow S$ be a smooth proper scheme. 

\begin{prop}\label{odb} The sheaves $f_*\mathcal O$, $f_*\mathcal O^*$,
and $R^1 f_*\mathcal O^*$ on $S_{fppf}$ are
representable. The scheme $T_X:=
Hom(D_X, \gm)$  represents $f_*\mathcal O^*$. 
\end{prop}

\begin{proof} The representability of $f_*\mathcal O$ and
$f_*\mathcal O^*$ is rather elementary 
\cite[Cor. 8, Lem. 10, pp.207-208]{blr}. Each character of $D_X$ provides a non-zero function,
constant (since $X$ is 
proper) on each connected component of $X$, i.e., on each irreducible
component of $X$ (since $X$ is smooth). Thus $T_X$
represents $f_*\mathcal O^*$.  The representability of  $R^1
f_*\mathcal O^*$  is due to Murre-Oort \cite[Thm. 3,
p.211]{blr}.\end{proof}{}  

The scheme representing $R^1
f_*\mathcal O^*$ is the Picard scheme $Pic_X$ of
$X$. It is reduced in characteristic
zero but it may not be so in positive characteristic. Its
reduced neutral (= identity) component $Pic^{0, red}_X$ is the classical
Picard variety 
$Pic(X)$. The
N\'eron-Severi group scheme $NS_X$ is the
\'etale group scheme corresponding to the $\G$-module 
$\pi_0(Pic^{red}_X)$; we often write $NS(X)$ for $NS_X(S)$.

\begin{prop}\label{memu1} {{\rm ($p=0$)}} The $S_{fppf}$-sheaves 
$R^if_*\Omega ^{j}$ and $\mathcal
H^i_{DR}(X) = R^if_*\Omega$ given by Hodge and de Rham cohomology
{\rm (}$\Omega$ is the de Rham complex on $X${\rm )} as well as the sheaf
$R^1f_*\Omega^*$ corresponding to the 
multiplicative de Rham complex \cite[3.1.7,
p.31]{mame}
on $X$  $$\Omega^*:= [\mathcal O^*
\xrightarrow{d~{\rm log}} \Omega^1 
\rightarrow \Omega^2 \cdots ]$$ are all 
representable. The first two are representable by vector group
schemes. \end{prop}

\begin{proof} The sheaves
$R^if_*\Omega^j$, $\mathcal H^i_{DR}(X)$ are coherent, free, and
commute with
arbitrary base change \cite[1.4.1.8]{katzdiff},
\cite[p.309-310]{grocry}. For any $t: S' \to S$, we 
have $R^if_*\Omega^j (S') = t^*
H^i(X,\Omega^j)$; similarly for $\mathcal H^i_{DR}(X)$.  Any locally
free $\mathcal O_T$-module $L$ on a scheme $T$ gives rise to a sheaf
on $T_{fppf}$ which is representable by a vector group scheme
\cite[p.1]{mame}. Thus, $R^if_*\Omega^j$, $\mathcal H^i_{DR}(X)$ are
representable by vector group schemes. Similarly, the sheaves
$R^if_*C$  ---  here $C$ is $[\Omega^1 \to \Omega^2
\cdots]$  ---   are also representable by vector group schemes.  
In the exact sequence\footnote{The natural map $f_* \mathcal O^* \to
f_*C$ is zero: the former is represented by a torus and the latter by
a vector group scheme.} (cf. (\ref{gumma}))$$0 \to f_*C \to R^1f_*\Omega^* \to R^1f_* \mathcal O^*
\to R^1f_*C,$$representability is already known for 
all the sheaves other than $R^1f_*\Omega^*$; 
and $f_*C$ is  representable by an affine
scheme. The representability of $R^1f_*\Omega^*$ follows from
(\ref{repco}). \end{proof} 

The scheme $Pic^{\natural}_X$ representing 
$R^1 f_*\Omega^*$ classifies \cite[2.5]{mes} \cite[7.2.1]{katzdiff} isomorphism classes of 
line bundles (=invertible sheaves) on $X$ endowed with an integrable
connection; cf. (\ref{haccha1}). 

\begin{prop}\label{memu2}  Assume that $k$ is of characteristic zero.

{\rm{(i)}} The neutral  component $E^{\natural}$ of
$Pic^{\natural}_X$ is the universal additive {\rm(}= vectorial{\rm)} extension
  of $Pic^0_X = Pic_X^{0, red}$. 

{\rm{(ii)}} The additive group scheme $Lie~E^{\natural}$ represents the sheaf
$\mathcal H^1_{DR}(X)$;  

{\rm Lie} $E^{\natural} = {\rm Lie}~Pic^{\natural}_X
  \xrightarrow{\sim} H^1_{DR}(X):= 
\hi^1(X,\Omega) \xleftarrow{\sim} 
\hi^1(X,[\mathcal  O \rightarrow \Omega^1])$.
\end{prop}
\begin{proof} A proof of (i) for  $X$ an
abelian variety is in \cite[2.1, 2.7, 2.8]{mes}; it can also be
  obtained by  combining propositions 2.6.7, 3.2.3, and 4.2.1 of
  \cite[Chapter I]{mame}.  

For a general smooth proper $X$ over $S$, let $g: E^* \to  E^{\natural}$
  be the map induced   from the universal additive extension $E^*$ of
  $Pic^0_X$.  We need to show that $g$ is an
  isomorphism. Since  both  $E^*$ and $E^{\natural}$ 
  are compatible with base change, we may assume $X(S) \neq
\varnothing$; this provides an Albanese map $u: X \to Alb(X)$. We now
  argue as in \cite[3.0]{mes} using  the standard
isomorphisms $u^*: H^0(Alb(X), \Omega^1) \xrightarrow{\sim} H^0(X,
\Omega^1)$ and $u^*: Pic^0_{Alb(X)} \xrightarrow{\sim} Pic^0_X$.  
Part (ii) follows
from (i) by  \cite[Lem. 2.6.9]{mes}. \end{proof}

\nd {\bf Divisors on a smooth proper variety.}\ssp

 We
recall the classical properties of $Pic(X)$ and $NS(X) = NS_X({S})$
for $f:X \to S$ smooth proper.  
 
\begin{rem}\label{picmd} One has 

(i) ($\ell \neq p$) an isomorphism $H^1_{et}(\bar{X},\Z_{\ell}(1)) \xrightarrow{\sim} 
T_{\ell}Pic(X)$  \cite[p.125]{mi} of
$\G$-modules  provided by the Kummer sequence (\ref{kummer}). 

(ii) ($p =0$)  Lie $Pic(X) \xrightarrow{\sim} H^1(X,\mathcal O)$
\cite[Thm. 1, p.231]{blr}.

(iii) ($p=0$) an exact sequence (\ref{memu2}) $$0 \to H^0(X,\Omega^1) \to {\rm
Lie}~E^{\natural} \to H^1(X,\mathcal O) \to 0.$$
 
(iv) ($k = \C$)  an isomorphism of pure Hodge structures
\cite[pp.156-158]{cime}: 
\begin{equation}\label{hdgpic}
H^1(X(\C),\Z(1)) \xrightarrow{\sim} H_1(Pic(X),\Z)
\end{equation}
provided by the exponential sequence
\begin{equation}\label{exp}
0 \to \Z(1) \to \mathcal O \xrightarrow{exp} \mathcal O^* \to 1;
\end{equation}a 
commutative diagram \cite[Thm. 1.4, p. 17]{dmos} --- vertical maps are
\cite{gaga}:  
$$
\begin{CD}
0 @>>> H^0(X,\Omega^1) @>>>  H^1_{DR}(X) @>>> H^1(X,\mathcal O) @>>>
0\\
@. @V{\wr}VV @V{\wr}VV @V{\wr}VV @.\\
0 @>>> H^0(X^{an},\Omega^1) @>>> H^1(X(\C),\C) @>>>
H^1(X^{an},\mathcal O) @>>>0.\\
\end{CD}$$

(v)  ($\ell \neq p$) the ``cycle class map'' \cite[VI \S 9]{mi} furnishes a
$\G$-equivariant inclusion \cite[3.2.9 (d), p.216]{mi}:  
\begin{equation}\label{algs}
NS(\bar{X})\otimes_{\Z} \Z_{\ell} \hookrightarrow H^2_{et}(\bar{X},
\Z_{\ell}(1));
\end{equation}
numerical and homological
equivalence coincide for divisors (with $\Q_{\ell}$-coefficients).

(vi) ($k =\C$) Lefschetz's $(1,1)$-theorem = the integral Hodge
conjecture 
for divisors  \cite[p.143]{dmot} \cite[p. 156]{cime}:
\begin{equation}\label{h11} 
NS(X)  \xrightarrow{\sim}  {\mbox{Hom}}_{MHS}(\Z,
H^2(X(\C),\Z(1))). 
\end{equation} 

(vii) ($k$ finitely generated) the Tate conjecture \cite[p.72]{tate}
\cite[5.1]{jamm}  for divisors
asks if
\begin{equation}\label{tatec}
NS(X)\otimes_{\Z}\Q_{\ell} \xrightarrow{\sim}
H^2_{et}({X}\times k^{sep},\Q_{\ell}(1))^{\G_{sep}} ?
\end{equation} 
Here $k^{sep}$ is a separable algebraic closure of $k$, $\G_{sep}$
the associated Galois group, and $M^{\G_{sep}}$ denotes the invariants of a $\G_{sep}$-module
$M$. See \cite[\S 5]{tate} for the known cases of (\ref{tatec}).\ssp 
  
(viii) ($p=0$) the Chern class map \cite[2.2.4]{h2} 
gives an injection
\begin{equation}\label{beagle}
c_X: NS(X)\otimes k \hookrightarrow H^1(X, \Omega^1).
\end{equation}

(The map $d~{\rm log}: \mathcal O^* \to \Omega^1$ \cite[2.2.4]{h2} induces 
a map $R^1f_*\mathcal O^* \to R^1f_*\Omega^1$ of representable 
$S_{fppf}$-sheaves; since the second is affine and the identity
component of the first is an abelian scheme, one gets the map $c_X$.)
\qed \end{rem}     

\begin{rem}\label{dcon} Let us consider (\ref{dc1}) for
$X$. By the purity of $H^n(X, \Z(1))$ ($k = \C$), 
both $t^n(X)$ and $I^n(X)$ are
zero for $n>2$. Parts (i)-(iv) of (\ref{picmd}) identify $[0 \to
Pic(X)]$  as
the one-motive $L^1(X/k)$  
whereas parts (v), (vi), and (viii) identify 
$[NS_X/{\rm torsion} \to 0]$ as the one-motive $L^2(X/k)$.
Conjecture \ref{dc1} is known for smooth proper
varieties.\qed \end{rem} 

From (\ref{h11}), one obtains the 
\begin{thm}\label{1-1i} Let $X$ be a smooth proper scheme over $S =
{\rm Spec}~k$. The dimension  
of the $\Q$-vector space $$H^{1,1}_{\Q}(X_{\iota}):= {\rm Hom}_{MHS}(\Q(-1),
H^2(X_{\iota},\Q))$$ is independent of the map $\iota: k \hookrightarrow
\C$. \qed\end{thm}    
 
\nd {\bf Smooth varieties.}\ssp

Let $U$ be the open complement of a strict
divisor $Y$ (\ref{simpxy}) with normal crossings 
in a smooth projective complex variety
$X$. The normalization $\tilde{Y}$ of $Y$ is a smooth projective
scheme. Consider the map $W_Y \to Pic_X$ which sends a divisor $E$ to the
class of the invertible sheaf $\mathcal O(E)$ on $X$. Let $N$ be the
cokernel of the induced map $\lambda: W_Y
\to NS_X$. 

\begin{prop}\label{podade} The $(0,0)$-part of $H^2(U,
\Q(1))$ is $N\otimes\Q$. So $t^2(U)\otimes\Q 
\xrightarrow{\sim} N\otimes\Q$.  
\end{prop} 

\begin{proof} Let $j: U \hookrightarrow X$ denote the
inclusion.        
The cycle  class map \cite[2.2.4-5]{h2} and the Gysin sequence
\cite[3.3]{hp}  provide the  following
commutative  diagram
$$ \begin{CD}
W_{Y} @>{\lambda}>>  NS({X}) @>>> N @>>> 0\\
@| @VVV @VVV @.\\
H^0(\tilde{Y},\Z) @>>> H^2(X,\Z(1)) @>{j^*}>> H^2(U,\Z(1))
@. {}.\\
\end{CD}
$$
By \cite[3.2.17]{h2}, we know that $W_0 H^2(U,
 \Q(1)) = j^*(H^2(X,\Q(1)))$. This implies that the $(0,0)$-part of
$H^2(U,\Q(1))$ is the image under $j^*$ of the $(0,0)$-part of
 $H^2(X,\Q(1))$.  From (\ref{h11}), $NS(X)$ is the 
 $(0,0)$-part of $H^2(X,\Z(1))$. This proves the first claim. The
 second follows immediately because $t^2(U)\otimes\Q$ is the
 $(0,0)$-part of $H^2(U, \Q(1))$: by \cite[3.2.15 (ii)]{h2},
 $W_{-1}H^2(U, \Q(1)) =0$.\end{proof}

\begin{rem}\label{toto} Since $W_{m-1}H^m(U, \Q) =0$ \cite[3.2.15]{h2} for any $m \ge
0$,  only $t^1(U)$ and $t^2(U)$ can be nontrivial. Conjecture
\ref{dc1} for the case of $t^1(U)$ is classically known (\ref{curie});
$L^1(U/k)$ can be identified as the Picard one-motive $[\Ker(\lambda) \to
Pic(X)]$ of \cite{ra}, the Cartier dual of the generalized Albanese
variety of $U$ 
\cite{jp2}.  The  result (\ref{podade}), an analog of
(\ref{h11}) for smooth varieties, identifies $[N \to 0]\otimes\Q$ as 
 $L^2(U/k)\otimes\Q$. Thus, (\ref{dc1}) is known (up to isogeny) for the smooth
variety $U$. \qed \end{rem}

\nd {\bf Simplicial objects.} \cite{curt}\ssp  

 A simplicial object $B_{\bu}$ in a category $\mathcal C$ is a sequence
 of objects $$B_{\bu} = \{B_0, B_1, B_2, ..., B_n, ...\},$$together with
morphisms (face) $ d_i: B_n \rightarrow B_{n-1}$, (degeneracies) $ s_i: B_n
\rightarrow B_{n+1}$ ($0 \le i \le n$) satisfying the simplicial
 identities of which we need only (the degeneracy maps  will not play
 a role in our discussions) 
\begin{equation}\label{att} 
d_i d_j  =  d_{j-1}d_i  \qquad i <j.
\end{equation} 
A truncated simplicial object $B_{\ge m}$ is
the subsequence of objects $\{B_m, B_{m+1}, B_{m+2}, ...\}$,
together with all the maps $d_i$ and $s_j$ between them.\ssp

  Given a simplicial (resp. cosimplicial) commutative group scheme $A_{\bu}$, then the pair
$(A_{\bu}, \delta)$  ---  here 
$\delta_n = \Sigma_{i=0}^{i=n} (-1)^i d_i$
\cite[3.4]{curt}  ---  becomes a chain complex of commutative group
schemes: it follows from (\ref{att}) that $\delta_n \delta_{n+1}
= 0$ (resp. $\delta_{n+1} \delta_n = 0$).\ssp

\nd {\bf Simplicial pairs.}\ssp 

 Let $a:X \rightarrow S$
be a smooth projective morphism. Let $Y$ be a {\emph{strict}}  divisor
with normal crossings on $X$ as in \cite[2.4]{dj}; in
particular, this means that $Y$ is reduced and its irreducible
components $Y^i (i \in I)$ are regular schemes, and of codimension one in
$X$.  Let $U$ be the open subscheme of 
$X$ corresponding to the complement of $Y$. The normalization 
$\tilde{Y}$ of $Y$ is the  disjoint  sum of the $Y^i$'s. 

\begin{defn}\label{simpxy} A simplicial pair $(X_{\bu}, Y_{\bu})$ consists of the data
of    

(i) a simplicial scheme $X_{\bu}$ smooth and projective over $S$. 

(ii) a strict divisor with normal crossings $Y_{\bu}$ of $X_{\bu}$; in
particular, each $Y_m$ is a strict divisor (as defined above) of
$X_m$. 

(These conditions are ($\alpha$), ($\beta$) of \cite[p. 51]{dj}; it
follows from \cite{dj} that any
variety over  $S$ ($k$ perfect) admits a proper hypercovering (\S \ref{varieties})
corresponding to a simplicial pair.)\end{defn} 

In particular, $(X_{\bu}, Y_{\bu})$ is a simplicial object in the category 
of pairs. The schemes $U_{m}:= X_{m} - Y_{m}$ give 
a smooth simplicial subscheme $U_{\bu}$ of $X_{\bu}$ \cite[6.2.6]{h}; let 
$j: U_{\bu} \to X_{\bu}$  be the natural map.   

 A simplicial pair $(X_{\bu}, Y_{\bu})$ gives, for each $l\ge 0$, a  
truncated simplicial pair  $(X_{\ge l},Y_{\ge l}) =
(X_{\bu},Y_{\bu})_{\ge l}$ which
consists of the schemes $(X_m,Y_m)$ (for ${m \ge l}$) and 
the  maps $d, s$. 

One has an evident notion \cite[6.2.8]{h} \cite[p. 75]{sga4} of a 
morphism $\theta: (X_{\bu},Y_{\bu}) \to
 (Z_{\bu},J_{\bu})$  of simplicial pairs; $\theta$ satisfies  
 \cite[1.10]{jannbson} $\theta (Y_m)
 \subset J_m$ and $\theta (X_m - Y_m) \subset (Z_m - J_m)$.   

\begin{rem}\label{eggses} Fix a simplicial pair 
$(X_{\bu}, Y_{\bu})$. The simplicial abelian group scheme
$D_{\bu}$ --- here $D_i = D_{X_i}$ --- gives upon normalization a chain
complex still denoted $D_{\bu}$. Similarly, one obtains  
the chain complexes   of group schemes: $T_{\bu}$ 
 $Pic_{\bu}$, $NS_{\bu}$, $Pic^0_{\bu}$, $W_{\bu}$. The map
$\lambda_i: W_i = W_{Y_i} \to NS_{X_i} = NS_i$ which sends a divisor
$E$ of $X_i$ supported on $Y_i$ to the class $[\mathcal
O(E)]$ of the invertible sheaf $\mathcal O(E)$ gives a map $\lambda:
W_{\bu} \to Pic_{\bu} \to NS_{\bu}$.\end{rem}

\nd {\bf Exact and spectral sequences.}\ssp

 We summarize some results from
\cite[5.1-5.3]{h} about cohomology
of sheaves on a simplicial scheme $Z_{\bu}$. 
  For any sheaf (or complex of sheaves) $F$ on $Z_{\bu}$, there is a
spectral sequence \cite[5.2.3.2, 5.1.12.2]{h} \footnote{Note that the map $d_i$ of
 \cite[5.2.3.2, 5.3.3.2]{h} is denoted here by $\delta_i$ and vice-versa.}
 \begin{equation}\label{zsfbu} 
 E_1^{p,q} = H^q(Z_p,F) \Rightarrow H^{p+q}(Z_{\bu},F)
 \end{equation} with associated low-degree exact sequence
 \begin{equation*}
  0 \rightarrow E^{1,0}_2(F) \rightarrow H^1(Z_{\bu},F) 
 \rightarrow E^{0,1}_2(F) \rightarrow E^{2,0}_2(F). 
 \end{equation*}
The ``filtration b\^ete'' $\sigma$ of Deligne \cite[1.4.7]{h2} is: 
\begin{equation}\label{sigcon}
\sigma_{\ge m}H^*(Z_{\bu}, F):=  \im (H^*(Z_{\bu},
\sigma_{\ge m} F ) \to  H^*(Z_{\bu},  F)).
\end{equation}
 The methods used to deduce (\ref{zsfbu}) in \cite[5.2.3, 5.2.7]{h}
also work for truncated simplicial schemes.
 Given a complex of sheaves $C^{\bu}$ on $Z_{\ge m}$, there is a
spectral sequence  
 \begin{equation}\label{ssfbu}
 E_1^{p+m,q}(C^{\bu}) = H^q(Z_{p+m},C^{\bu}) \Rightarrow
H^{p+q}(Z_{\ge m},C^{\bu}) 
 \end{equation} with associated low-degree exact sequence
 \begin{equation}\label{ssylow}
  0 \rightarrow E^{1+m,0}_2(C^{\bu}) \rightarrow \hi^1(Z_{\ge m},C^{\bu}) 
 \rightarrow E^{m,1}_2(C^{\bu}) \rightarrow E^{2+m,0}_2(C^{\bu}). 
 \end{equation}Let $g_r: Z_r \to S$ and 
$f:Z_{\ge m} \rightarrow S$ be 
 the structure maps. Analogs of (\ref{ssfbu}) and (\ref{ssylow}) 
for the associated $S_{fppf}$-sheaves hold as well:
\begin{equation}\label{fppfbu} 
 E^{p+m,q}_1 (C^{\bu})= R^qg_{p+m ~*}(C^{\bu}) \Rightarrow
 R^{p+q}f_*(C^{\bu})
 \end{equation}
 \begin{equation}\label{fppflow}
  0 \rightarrow E^{1+m,0}_2(C^{\bu}) \rightarrow R^1f_*(C^{\bu}) 
 \rightarrow E^{m,1}_2(C^{\bu}) \rightarrow E^{2+m,0}_2(C^{\bu}). 
 \end{equation}

\section{Construction of $L^n$ for a simplicial pair}\label{constr}
Throughout this section, $k$ will denote a field of characteristic
zero. 

For each simplicial pair $(X_{\bu},Y_{\bu})$ and each non-negative 
integer $n$, we
shall construct one-motives $$L^n =L^n(X_{\bu},Y_{\bu}) = [\mathcal B_n
\xrightarrow{\phi_n} \tilde{\mathcal P}_n]$$ and
$J^n(X_{\bu},Y_{\bu})$;  these  are
contravariant functorial. Their realizations will be analyzed in
subsequent sections. Our construction was
obtained by a careful analysis of \cite[8.1.19.1]{h}.  

Fix a simplicial pair $(X_{\bu},Y_{\bu})$;
we have a diagram $U_{\bu} \overset{j}{\hookrightarrow} X_{\bu}
\overset{i}{\hookleftarrow} Y_{\bu}$. Let $a: X_{\bu}  \rightarrow S$, 
$a_m:X_m \rightarrow S$,  and  $f: X_{\ge
  n-1} \rightarrow S$ be  the structure morphisms.\footnote{We adopt
the notation: $X_{\ge -1} = X_{\bu}$.} We often write
$\mathcal O^*_m$ for  $\mathcal O^*_{X_m}$, $Pic_{m}$ for $Pic_{X_m}$,
etc.\medskip
 
\nd {\bf Construction of $\mathcal P_n$.}

We begin with the construction of a semi-abelian variety $\mathcal
P_n$ which is isogenous to the required $\tilde{\mathcal P}_n$; our
method is similar to \cite[\S3.1]{ra}.  

\begin{prop}\label{reps}
{\rm{(i)}} The sheaf $R^1f_* \mathcal O^* = R^na_*\sigma_{\ge
  n-1} \mathcal O^*_{\bu}$ is representable by a locally algebraic group
scheme $\mathcal G'_n$ with neutral
component $\mathcal G_n$.

{\rm{(ii)}} ${\rm Lie}~ \mathcal G_n \zim H^1(X_{\ge n-1}, \mathcal O)$. 

{\rm{(iii)}} {\rm{(}}$k = \C${\rm{)}} $H^i(X_{\ge n-1},\mathcal F)
\xrightarrow{\sim} H^i(X^{an}_{\ge n-1}, \mathcal F)$, $\mathcal F
= \mathcal O, \mathcal O^*$ and $i =0,1$. 
\end{prop}
 
\begin{proof} (i) By (\ref{fppflow}), the sheaf $R^1f_*\mathcal O^*$
 sits in an exact sequence \begin{equation*}
 0 \rightarrow E^{n,0}_2 \rightarrow R^1f_*\mathcal O^*
 \xrightarrow{\pi} E^{n-1,1}_2 \xrightarrow{\psi} E^{n+1,0}_2.
 \end{equation*}
By (\ref{odb}), the sheaf $a_{m~*} \mathcal O^*_m$ is
representable by the torus $T_{X_m}$. The group (of multiplicative
 type) $Hom(H^m(D_{\bu}),
\gm)$ dual to the 
homology $H^m(D_{\bu})$ (\ref{eggses}) of 
$$D_{X_{m+1}} \xrightarrow{\delta_{m}} D_{X_{m}}
\xrightarrow{\delta_{m-1}} D_{X_{m-1}}$$ 
represents $E^{m,0}_2$. 
Since $Pic_{X_m}$ represents $R^1a_{m~*}\mathcal O^*_m$, the scheme 
$ \mathcal R':= \Ker({\delta^*_{n-1}}:Pic_{X_{n-1}} \rightarrow
Pic_{X_n})$ represents 
$E^{n-1,1}_2$. Since $\Ker(\psi)$ is representable and the scheme
representing  $E^{n,0}_2$ is affine, we can apply (\ref{repco}). 

(ii) It follows from (i) by 
\cite[Lem. 2.6.9]{mes}.

(iii)  Use the GAGA isomorphisms 
$H^i(X_m, \mathcal F) \xrightarrow{\sim} H^i(X_m^{an}, \mathcal F)
\quad (i = 0,1) $ \cite[Prop. 17, 18]{gaga} for
 $\mathcal F = \mathcal O, \mathcal O^*$ and (\ref{ssfbu}),
(\ref{ssylow}) for $X_{\ge n-1}$ and $X_{\ge n-1}^{an}$. \end{proof}  

We call $\mathcal T'$ (resp. $\mathcal Q'$)  the
  group scheme representing $E^{n,0}_2$ (resp. $E^{n+1,0}_2$); we
  denote its 
  neutral component  by $\mathcal T$ (resp. $\mathcal
  Q$). 

The neutral component $\mathcal R$ of $\mathcal R'$ is an
abelian scheme (it is a subscheme of $Pic^0_{X_{n-1}}$)  
and $\mathcal Q$ is affine; so $\psi(\mathcal R) =
0$. Hence  $\mathcal R$ is the neutral component of $\Ker(\mathcal R'
\xrightarrow{\psi} \mathcal Q')$.   Thus we have an exact sequence 
\begin{equation}\label{lcw} 0
\rightarrow \mathcal T' 
\rightarrow \mathcal G_n \xrightarrow{\pi} \mathcal R \rightarrow
0.\end{equation}

\begin{rem}\label{avva} (i) By definition \cite[1.1]{gi}, an
invertible sheaf $L$ on a
  simplicial scheme $Z_{\bu}$ is an invertible sheaf $L_m$ on each
  $Z_m$ such that: for each morphism $\tau: Z_m \to Z_r$ which is a
 composition of $s_i$'s and $d_j$'s, the map
  $\tau^* L_r \to L_m$ is an isomorphism. In (loc. cit), it is 
  shown that $L$ is determined entirely by the data of $L_0$ and the
  isomorphism $\alpha: d_0^* L_0 \zim d_1^*L_0$ satisfying
  $d_2^*(\alpha) \circ d_0^*(\alpha) = d_1^*(\alpha)$; (\ref{reps}),
  for $n=1$, proves the representability of the Picard functor of
  $Z_{\bu}$ which is smooth and proper over $S$. 
The scheme $\mathcal G'_n$ classifies isomorphism classes of pairs
  $(\mathcal L, \alpha)$ where 

(a) $\mathcal L$ is an invertible sheaf on $X_{n-1}$ such that 
$\delta^*_{n-1} \mathcal L$ is isomorphic to $\mathcal O_{X_n}$ (by
(\ref{att}), $\delta^*_{n}\delta^*_{n-1} \mathcal L = \mathcal
O_{X_{n+1}}$);   

(b) $\alpha$ is a trivialization  $\delta^*_{n-1} \mathcal L
\xrightarrow{\sim} \mathcal O_{X_n}$ on $X_{n}$ satisfying a cocycle
condition: $\delta^*_{n} \alpha = 1$; namely,
$\delta_n^* \alpha$ is the identity isomorphism of $\mathcal
O_{X_{n+1}}$. 

(ii) For $n =1$, we can think of $\alpha$ as an isomorphism $d_0^*
\mathcal L \xrightarrow{\sim} d_1^*\mathcal L$. The cocycle condition
$d_2^*(\alpha)\circ d_0^*(\alpha) = d_1^*(\alpha)$ ($\Longleftrightarrow
\delta_1^*\alpha =0$) becomes the 
commutativity of 
 $$
\begin{CD}
d_0^* d_0^* \mathcal L @>{d_0^* (\alpha)}>> d_0^* d_1^*
\mathcal L @ =
d_2^* d_0^* \mathcal L \\
@|    @.   @V{d_2^* (\alpha)}VV \\
d_1^* d_0^* \mathcal L @>{d_1^* (\alpha)}>> d_1^*
d_1^* \mathcal L @=
d_2^* d_1^* \mathcal L\\ 
\end{CD} {} 
$$ 
 (using the identities $d_0d_1=d_0d_0, d_0d_2 
= d_1d_0, d_1d_2 = d_1d_1$).\qed
\end{rem} 

\begin{rem}\label{kaku}
Suppose given a simplicial scheme $Z_{\bu}$, an
invertible sheaf $\mathcal F$ on $Z_m$ and a nowhere vanishing section
$s$ of $\mathcal F$. The identity (\ref{att}) implies
that (i) the sheaf $\delta^*_{m+1}\delta^*_m\mathcal F$ 
is  naturally isomorphic to $\mathcal O_{Z_{m+2}}$; (ii) the section
$\delta^*_{m+1}\delta^*_ms$ of 
$\delta^*_{m+1}\delta^*_m\mathcal F = \mathcal O_{Z_{m+2}}$
corresponds to the identity section of $\mathcal O_{Z_{m+2}}$.\qed  
\end{rem}

\nd {\bf Definition of $\mathcal B_n$ and $\phi_n$.}\ssp
 
We now turn to the construction of the map $\phi_n$. This will take
several steps. 

Let $\lambda'_m: W_m \to Pic_{X_m}$ be the map $E \mapsto \mathcal
O(E)$. 
In general, there is no lifting of $\lambda'_{n-1}: W_{Y_{n-1}}
\to  Pic_{X_{n-1}}$ to a map $W_{Y_{n-1}}
\to \mathcal G'_n \xrightarrow{\pi} R' \rightarrow Pic_{X_{n-1}}$. In
other words, the invertible sheaf $\mathcal O(E)$ on $X_{n-1}$
corresponding to $E \in W_{n-1}$ may not always satisfy conditions
(a), (b) of (\ref{avva}(i)).  
But  a natural lifting 
does exist on the subgroup $K:= \Ker(\delta^*_{n-1}: W_{Y_{n-1}} \to
W_{Y_n})$ of $W_{Y_{n-1}}$. 

\begin{lem}\label{liftss}
There exists a canonical and functorial map $\vartheta': K \rightarrow
\mathcal G'_n$ which fits into a  commutative diagram:
$$
\begin{CD}
K @>{\vartheta'}>> \mathcal G'_n\\
@VVV @VVV\\
W_{Y_{n-1}} @>{{\lambda}'_{n-1}}>> Pic_{X_{n-1}}.\\
\end{CD}
$$ 
\end{lem}

\begin{proof} The scheme $\mathcal G'_n$
classifies pairs $(\mathcal L, \alpha)$ where $\mathcal L$ is an
invertible sheaf on $X_{n-1}$ and $\alpha$ is a trivialization 
$\delta^*_{n-1}\mathcal L \xrightarrow{\sim}\mathcal O_{X_n}$
 satisfying the cocycle condition $\delta^*_{n}\alpha =1$. For any $E 
\in K$, pick a rational section $s_E$ of $\mathcal O(E)$ such that the
divisor $div(s_E)= E$; the set of such sections forms a torsor over
$H^0(X_{n-1},\mathcal O^*)$. Since  $Y_{\bu}$
is a simplicial divisor, the pull-back $\delta^*_{n-1}s_E$ is a
rational section of $\delta^*\mathcal O(E)$ with divisor
$\delta^*_{n-1}E$. As $E$ lies in $K$, we have $\delta^*_{n-1}E
=0$. So $\delta^*_{n-1}s_E$ provides a 
trivialization $\alpha_E: \delta^*_{n-1}\mathcal O(E) \xrightarrow{\sim} 
\mathcal O$ on $X_n$. The trivialization $\delta^*_n \alpha_E$ of
$\delta^*_n  
\delta^*_{n-1}\mathcal O(E)$ on $X_{n+1}$ is provided by the nowhere
vanishing regular section $t_E:= \delta^*_n
\delta^*_{n-1} s_E$. Now the rational section $s_E$ is a nowhere
vanishing section of $\mathcal O(E)$ on the open subscheme
$U_{n-1}$. Applying (\ref{kaku}) to $U_{\bu}$. we see that $t_E$
is the identity section of the sheaf 
$\mathcal O_{U_{n+1}}$. Since $t_E$ is regular on $X_{n+1}$, we
deduce that $t_E$ is the identity section of $\mathcal
O_{X_{n+1}}$.  The map $\vartheta'$ is
defined by $E \mapsto (\mathcal O(E), \alpha_E)$. It is clear that 
modifying $s_E$ by an element of $H^0(X_{n-1},\mathcal O^*)$ 
does not affect the isomorphism class of the
pair $(\mathcal O(E), \alpha_E)$.
\end{proof}

The maps $d_i:X_{\ge n-1}
\ra X_{n-2}$ collectively  provide the morphism 
\begin{equation}\label{lift}
\delta^*: Pic_{X_{n-2}} \to \mathcal G'_n.
\end{equation}
More explicitly,  given an
invertible sheaf $\mathcal L$ on $X_{n-2}$, the pull-back $\mathcal
L':= \delta^*_{n-2} \mathcal L$ is an invertible sheaf on
$X_{n-1}$. As in (\ref{kaku}), (\ref{att}) implies that there
is a canonical trivialization, call it $\beta_{\mathcal L'}$, of
$\delta^*_{n-1} 
\mathcal L'$ on $X_{n}$ satisfying $\delta^*_n \beta_{\mathcal
L'} = 1$. The map $\delta^*$ sends $\mathcal L$
to the (isomorphism class of the) pair $(\mathcal L', \beta_{\mathcal
L'})$. 

\begin{defn}\label{tildel} $\mathcal P'_n$ is the quotient of
 $\mathcal G'_n$ by ${\delta}^*(Pic^0_{X_{n-2}})$; its neutral
component is $\mathcal P_n$.
Note $\pi_0(\mathcal G'_n) \xrightarrow{\sim} \pi_0(\mathcal
P_n')$. 
\end{defn}

\begin{lem}\label{koil} {\rm{(}}$k = \C${\rm{)}} 
${\rm Lie}~\mathcal P_n \xrightarrow{\sim} \frac{H^1(X_{\ge n-1},
\mathcal O)}{\delta^*H^1(X_{n-2}, \mathcal O)}  \xrightarrow{\sim}
\frac{H^1(X_{\ge n-1}^{an},
\mathcal O)}{\delta^*H^1(X^{an}_{n-2}, \mathcal O)}.$
\end{lem}

\begin{proof} Combine (\ref{reps}) (iii) with (\ref{picmd}
(ii)). \end{proof}{} 
 
 The map (\ref{lift}), in turn, induces a map  $NS_{n-2}
\xrightarrow{\mu} \mathcal P'_n$. Combining this with $\vartheta: 
K \xrightarrow{\vartheta'}
\mathcal G'_n \rightarrow \mathcal P'_n$, we get the map 
\begin{equation*}
\rho: K \oplus NS_{X_{n-2}}\rightarrow \mathcal P'_n \qquad
\rho(a,b) = \vartheta(a) + \mu(b).
\end{equation*}

We now turn to the introduction of several group schemes relevant
for the definition  of $\mathcal B_n$. The mapping cone complex of $\lambda: W_{\bu} \to NS_{\bu}$ (\ref{eggses}) is:
\begin{eqnarray}\label{iwas}
W_{Y_{n-2}}\oplus NS_{X_{n-3}} \xrightarrow{\gamma_{n-3}} 
W_{Y_{n-1}}\oplus NS_{X_{n-2}} \xrightarrow{\gamma_{n-2}} 
W_{Y_{n}}\oplus NS_{X_{n-1}}\\
\gamma_{m}:(E,\beta) \mapsto  (\delta^*_{m+1}E,
\delta^*_{m}\beta - \lambda_{m+1}(E)).  
\end{eqnarray} 
Since $W_{\bu}$ is a complex, the image of $\gamma_{n-3}$ is 
contained in $K\oplus NS_{X_{n-2}}$.
One checks that the composite map $\rho \circ \gamma_{n-3}$ is zero. 

\begin{defn}\label{auro} (i) $\mathcal C'_n$ denotes $\frac{K \oplus
NS_{X_{n-2}}}{\gamma_{n-3}(W_{n-2})}$; and $\mathcal C_n$ is the
kernel of the composite map $ \rho': \mathcal C'_n \xrightarrow{\rho}
\mathcal P'_n \rightarrow 
\pi_0(\mathcal P'_n)$.

(ii) $\mathcal A_n$ is the kernel  of the
composite map $\mathcal C'_n \xrightarrow{\rho'}\pi_0(\mathcal P'_n)
\xrightarrow{\sim}  \pi_0(\mathcal G'_n) \to NS_{X_{n-1}}$; this map  ---  
the middle 
isomorphism is the inverse of $\pi_0(\mathcal P'_n)
\xleftarrow{\sim}  \pi_0(\mathcal G'_n)$ --- is the map $\mathcal C'_n
\to NS_{X_{n-1}}$    
induced by $\gamma_{n-2}$. Note $\mathcal C_n$ is a subgroup of
$\mathcal A_n$. 

(iii) $\mathcal B'_n$ is the kernel of the map 
\begin{equation}\label{greece}
\frac{K\oplus NS_{X_{n-2}}}{\gamma_{n-3}(W_{n-2}\oplus NS_{X_{n-3}})}  \xrightarrow{\rho}
\mathcal P'_n \rightarrow \pi_0(\mathcal P'_n), 
\end{equation} 
which is analogous to the cycle class map;  $\pi_0(\mathcal
 P'_n)$ plays the role of a relative N\'eron-Severi
 group. \end{defn} 
We have\begin{equation}\label{bareesh}
\frac{\mathcal C_n}{NS_{X_{n-3}}} \zim \mathcal B'_n.\end{equation}
 Restricting $\rho$ to $\mathcal B'_n$, we obtain a morphism of 
group schemes
$\mathcal B'_n \xrightarrow{\rho} \mathcal P_n.$
Let $\tau_n$ be the torsion subgroup of $\mathcal B'_n$. Replacing
$\mathcal B'_n$ by $\mathcal B_n:= \mathcal B'_n/{\tau_n}$ and 
$\mathcal P_n$ by  
$\tilde{\mathcal P}_n:= \mathcal P_n/{\rho(\tau_n)}$, the map $\rho$
induces a map  
$\phi_n: \mathcal B_n \ra \tilde{\mathcal P}_n$.

\begin{defn} $L^n(X_{\bu},Y_{\bu})$ is the one-motive
\begin{equation}\label{phin}
L^n= L^n(X_{\bu},Y_{\bu}):=  [\mathcal B_n \xrightarrow{\phi_n}
\tilde{\mathcal P}_n]. 
\end{equation}
\end{defn}
 
\begin{rem}\label{zeusp} (i) The one-motive $L^n$ depends only on the
schemes $X_m$, $Y_m$ for $m=n-3, \cdots, n+1$ (and the maps between them).

(ii) Let $k'$ be an extension of $k$ and write
$S'=$~Spec~$k'$. It is clear from the definitions
that $L^n (X_{\bu} \times_S S', Y'_{\bu} \times_S S') = L^n \times_S S'$. 

(iii) By (i) and (ii), there is no loss of generality in assuming that
$k$ is finitely generated over $\Q$. Such a field always admits an
embedding into $\C$. 

(iv) $L^n$ is a contravariant functor for morphisms of simplicial pairs. 
\qed \end{rem}   

\nd {\bf The one-motive $J^n$.}\ssp

\begin{defn}\label{auro2} $\mathcal
V_m$ (resp. $N_m$) is the kernel (resp. cokernel) of $W_m
\xrightarrow{\lambda_m} NS_m$; and $K^0$ (resp. $M$) is the kernel (resp. cokernel) of the
composite map $K \xrightarrow{\vartheta} \mathcal P'_n \rightarrow
\pi_0(\mathcal P'_n)$. 

The one-motive $J^n$ is     
\begin{equation}\label{hosadu}
J^n = J^n(X_{\bu}, Y_{\bu}):= [(\frac{{K^0}/{\mathcal V_{n-2}}}{\rm
torsion}  \xrightarrow{\phi_n} \tilde{\mathcal P}_n].
\end{equation}
\end{defn}
The natural morphism $J^n \to L^n$ is an isomorphism for $n \le 1$;
 $J^n$ is contravariant 
 functorial for morphisms of simplicial pairs. 

The remainder of this section is devoted to
 results which will be used in \S \ref{Tmt}.

\begin{lem}\label{codered} One has 
an exact sequence$$0 \to
\frac{\mathcal V_{n-1} \cap K}{\mathcal V_{n-2}} \to
\mathcal A_n  \to 
\Ker(\delta^*_{n-2}: N_{n-2} \to N_{n-1})\to 0.$$
\end{lem}
\begin{proof}
Taking the quotient of $0 \to K \to K\oplus NS_{n-2} \to NS_{n-2} \to
0$ by the exact sequence $0 \to \mathcal V_{n-2} \to W_{Y_{n-2}} \to
\im(\lambda_{n-2}) \to 0$, we get the exact sequence in the first row of
the diagram
$$
\begin{CD}
 0 @>>> \frac{K}{\mathcal V_{n-2}} @>>> \frac{K\oplus
NS_{n-2}}{\gamma_{n-3}(W_{n-2})} @>>> 
N_{n-2} @>>> 0\\
@. @VVV @VVV @VVV @.\\
0 @>>> \im(\lambda_{n-1}) @>>> NS_{n-1} @>>>
N_{n-1} @>>> 0;\\
\end{CD}$$the required sequence is an easy consequence. 
\end{proof}

\nd {\bf The de Rham cohomology sheaves.}\ssp

We shall use the following complexes on a (simplicial) smooth
projective scheme; the scheme could be either $X_{\bu}$, $X_m$ or
$X_{\ge m}$. We often omit the subscript from the complex if it
is clear from the context which scheme it is on.  
\begin{eqnarray*}
\Omega & := &[ \mathcal O \xrightarrow{d} \Omega^1 \xrightarrow{d}
\Omega^2 \dots]\\
\Omega(log~Y) & :=  &[\mathcal O \xrightarrow{d} \Omega^1(log~Y)
\xrightarrow{d}...]\\
\Gamma_0 :=  [\mathcal O \xrightarrow{d} \Omega^1], & & \Gamma:=
[\mathcal O \xrightarrow{d} \Omega^1(log~Y)]\\ 
\Gamma^*_0  :=  [\mathcal O^* \xrightarrow{d~{\rm log}}
\Omega^1], & & \Gamma^*:= [\mathcal O^* \xrightarrow{d~{\rm log}}
\Omega^1(log~Y)]
\end{eqnarray*}
Let $\Omega(log~Y)$ be the logarithmic de Rham complex
\cite[3.1.3, 6.2.7]{h} 
of $U_{\bu}$ on $X_{\bu}$. 

\begin{lem}\label{dlog} The map $dl:\mathcal G'_n = R^1f_* \mathcal O^* \to
R^1f_*\Omega^1$ induced by $d~{\rm log}: \mathcal O^* \to \Omega^1$
\cite[2.2.4]{h2} yields a map $dlog: \pi_0(
\mathcal G'_n) \to H^1(X_{\ge n-1}, \Omega^1)$. \end{lem}
\begin{proof} Since $\mathcal G_n$ is a semi-abelian variety and
$R^1f_*\Omega^1$ is a vector scheme, the map $dl$ restricted to
$\mathcal G_n$ is zero.\end{proof}  

\begin{prop}\label{gumma} {\rm{(i)}} The sheaf   $R^1f_*\Gamma^*$ 
is representable by a group scheme $\mathcal G^{\diamondsuit}_n$ with
 $\hi^1(X_{\ge n-1},\Gamma^*)$ as its group of $k$-rational points. 

{\rm{(ii)}} Lie~$\mathcal G^{\diamondsuit}_n \zim \hi^1(X_{\ge n-1},
\Gamma)\zim H^1_{DR}(U_{\ge n-1})$.

{\rm{(iii)}} The group $\mathcal  G^{\diamondsuit}_n$ is an extension of a
subgroup {\rm (}containing $\mathcal G_n${\rm)} of $\mathcal G'_n$  by
${}f_*\Omega^1(log ~Y)$.  
\end{prop}

\begin{proof}
(i) In the exact sequence, given by the map $\Gamma^* \to \mathcal O^*$, 
\begin{equation*}
f_* \Gamma^* \xrightarrow{\beta} f_* \mathcal O^* \to 
f_* \Omega^1(log~Y) \xrightarrow{h}
 R^1f_*\Gamma^* \to R^1f_*\mathcal O^* \to R^1f_* \Omega^1(log~Y),
\end{equation*}
$\beta$ is clearly an isomorphism: $\hi^0(X_{\ge n-1}, \Gamma^*) =
H^0(X_{\ge n-1}, \mathcal O^*)$. So $h$ is injective. 
The sheaves $R^if_* \Omega^j(log~Y)$ are locally free and commute with
base change \cite[1.4.1.8]{katzdiff}; so they are representable by
vector group schemes (cf. the proof of (\ref{memu1})). 
Lemma \ref{repco} now provides the
representability of $R^1f_*\Gamma^*$. 

(ii)  the first isomorphism follows from (i) 
\cite[Lemma 2.6.9]{mes} and the second isomorphism from
\cite[1.0.3.7]{katzdiff}, \cite[3.1.8]{h2}; cf. also \cite[3.4]{jamm}.   

(iii) It suffices to show that the map Lie$~\mathcal G^{\diamondsuit}_n \ra$~
Lie$~\mathcal G_n$ is onto. Under the identifications, Lie~$\mathcal
G^{\diamondsuit}_n \zim \hi^1(X_{\ge n-1}, \Gamma)$ and Lie$~\mathcal G_n 
\zim H^1(X_{\ge n-1}, \mathcal O)$, we need the surjectivity of the 
forgetful map $  \hi^1(X_{\ge n-1}, \Gamma) \ra H^1(X_{\ge n-1}, \mathcal
O)$.   We can use (\ref{ssylow}) for $\Gamma$ and $\mathcal
O$ on $X_{\ge n-1}$; since $\hi^0(X_m,\Gamma) = H^0(X_m,\mathcal
O_{X_m})$, we have   $E^{n,0}_2(\Gamma) \xrightarrow{\sim} E^{n,0}_2(\mathcal
O)$.  Thus we need the surjectivity of $E^{n-1,1}_2(\Gamma) \ra
E^{n-1,1}_2(\mathcal O)$.  
Let us show that even the composite map 
$\pi_0: E^{n-1,1}_2(\Gamma_0) \ra E^{n-1,1}_2(\Gamma) \ra
E^{n-1,1}_2(\mathcal O)$ is onto.

Write $E_m^{\natural}$ for the 
universal additive extension of $Pic^0_{X_m}$. We have   
\begin{eqnarray*}
E^{n-1,1}_2(\Gamma_0) \xleftarrow{\sim} \Ker(\delta^*_{n-1}: {\rm Lie}~E^{\natural}_{n-1} \ra
{\rm Lie}~E^{\natural}_{n})\\ 
E^{n-1,1}_2(\mathcal O) \xleftarrow{\sim} \Ker(\delta^*_{n-1}: {\rm Lie}~Pic^0_{X_{n-1}}
\ra {\rm Lie}~Pic^0_{X_{n}}) = {\rm Lie}~\mathcal R.
\end{eqnarray*} 
For any abelian scheme $A$, the Lie algebra
Lie~$A^{\natural}$ of the universal additive (= vectorial) extension
is canonically isomorphic 
to $H^1_{DR}(A^*)$ of the dual abelian scheme $A^*$ \cite[4.1.7,
p. 48]{mame}  (\ref{memu2}).  So the functor ${\natural}: A
\mapsto {\rm Lie}~A^{\natural}$ is exact.  
Thus $E^{n-1,1}_2(\Gamma_0) \xleftarrow{\sim} {\rm 
Lie}~\mathcal R^{\natural}$. 
Hence $\pi_0$ is the natural  (surjective) projection Lie~$\mathcal
R^{\natural} \ra {\rm 
Lie}~\mathcal R$. \end{proof}    

\begin{rem}\label{diksha} 
(i)  For any semi-abelian variety $G$ (with associated abelian variety
$A$ and torus $T$), the universal  additive (= vectorial) extension
$G^{\natural}$ of $G$ is the pullback 
$A^{\natural} \times_A G$ \cite[10.1.7c]{h} of the universal additive
 extension $A^{\natural}$ of $A$. 
The proof of (\ref{gumma})  shows that 
(a) $R^1f_*\Gamma^*_0$ is a vectorial extension of a subgroup (containing 
$\mathcal G_n$) of $\mathcal G'_n$. (b) this vectorial extension
restricted to $\mathcal G_n$ is the universal extension $\mathcal
 G^{\natural}_n =  
\mathcal R^{\natural} \times_{\mathcal
R} \mathcal G_n$ of $\mathcal G_n$. The universal additive extension
$\mathcal P^{\natural}_n$ of $\mathcal P_n$ is the quotient of
$\mathcal G^{\natural}_n$ by $\delta^* E^{\natural}_{n-2}$, defined as
 in the proof of (\ref{gumma} (iii)).

(ii) The analog of (\ref{gumma}) is also true for a
smooth proper variety $Z$ and a strict divisor $V$ of $Z$. For the
choice of the strict divisor  $Y_{m}$ of $X_{m}$, we obtain the
existence of a group scheme  
$Pic^{\diamondsuit}_{X_{m}}$ with $\hi^1(X_{m},\Gamma^*)$ as its
group of  $k$-rational points; one has Lie~$Pic^{\diamondsuit}_{X_{m}}
\zim \hi^1(X_{m},\Gamma) \zim H^1_{DR}(U_m)$, the last isomorphism due
to \cite[3.1.8]{h2}, \cite[1.0.3.7]{katzdiff}. The cokernel $\mathcal
P^{\diamondsuit}_n$ of ${\delta}^*:
Pic^{\diamondsuit}_{X_{n-2}}  \ra 
\mathcal G^{\diamondsuit}_n$ is a vectorial extension of a subgroup
(containing 
$\mathcal P_n$) of $\mathcal P'_n$. 
 \qed \end{rem} 

\begin{lem}\label{ohoh!}  {\rm{(a)}}  {\rm (}$k
=\C${\rm )}  ${\rm Lie}~\mathcal
P^{\diamondsuit}_n \xrightarrow{\sim} 
\sigma_{\ge n-1} H^n(U_{\bu}, \C)$. 

{\rm{(b)}}  {\rm (}$k =\C${\rm )} ${\rm Lie}~\mathcal P^{\natural}_n
\xrightarrow{\sim} 
\sigma_{\ge n-1} H^n(X_{\bu}, \C) \xrightarrow{\sim} W_{-1}H^n(U_{\bu}, \C) $.

{\rm{(c)}} ${\rm Lie}~\mathcal
P^{\diamondsuit}_n \xrightarrow{\sim} 
\sigma_{\ge n-1} H^n_{DR}(U_{\bu})$. 
\end{lem}

\begin{proof} 
Combine $H^1(U_*,
\C) \xleftarrow{\sim} H^1(X^{an}_*, j_{*} \C_{U}) \xrightarrow{\sim}
\hi^1(X^{an}_*, \Omega(log~Y))$ 
(logarithmic Poincar\'e lemma \cite[3.2.2]{h2}),
$\hi^1(X^{an}_*,\Gamma) \xrightarrow{\sim} 
\hi^1(X_*^{an},\Omega(log~Y))$ ($X_* = X_m$  or $X_{\ge
n-1}$) with (\ref{gumma}), (\ref{diksha}).
Part (a) follows from the next isomorphism (the injectivity is a
consequence of the 
degeneration of \cite[8.1.19.1]{h} at $E_2$)  
\begin{equation}\label{patagon}
\frac{H^1(U_{\ge n-1}, \C)}{{\delta}^*_{n-2}H^1(U_{n-2}, \C)}\xrightarrow{\sim} \sigma_{\ge n-1}
H^n(U_{\bu},\C);
\end{equation} the surjectivity is by definition of $\sigma$.  The
first (resp. second) isomorphism in (b) 
is a special case of (a), i.e., for $Y = \varnothing$ (resp. from
(\ref{malli})). Part (c) is proved similarly using (\ref{gumma}),
(\ref{diksha}).\end{proof}

\section{Hodge and de Rham realizations of $L^n$}\label{Tmt} 

We retain the notations of the previous section but we now take $k
=\C$.  
We have a diagram $U_{\bu} \xrightarrow{j} X_{\bu} \xleftarrow{i}
Y_{\bu}$ 
corresponding to our simplicial pair $(X_{\bu}, Y_{\bu})$. 
As before, $f: X_{\ge n-1} \to S$ is the structure map.

The main
results (\ref{main}) (\ref{derham})  of this section  prove the Hodge
and de Rham  part of the conjecture (\ref{dc1}) for $U_{\bu}$; the
\'etale realization will be treated in \S 
\ref{oakland}.\ssp

\nd {\bf Statement of the theorem.}

\begin{lem} The mixed
  Hodge structures $H^*(U_{\bu},\Z)$ are polarizable. \end{lem}  

\begin{proof} Since  $Gr^W_rH^m(U_{\bu},\Q)$  is
a direct sum of subquotients of the cohomology of smooth projective
varieties \cite[8.1.19.1]{h} and the cohomology $H^*(V,\Q)$ of  
a smooth projective complex variety $V$ is polarizable
\cite[2.2.6]{h2}, this is clear. \end{proof} 
 
Denote the largest sub-mixed Hodge structure
of type $(*)$ of 
$H^n(U_{\bu},\Z(1))/{\rm torsion}$ by $t^n(U_{\bu})$; by the previous lemma,  
$Gr^W_{-1}t^n(U_{\bu})\otimes\Q = \Q(1)\otimes Gr^W_1H^n(U_{\bu}, \Q)$ is
 polarizable. By  \cite[10.1.3]{h}, the mixed Hodge structure 
$t^n(U_{\bu})$ corresponds to a one-motive
$I^n(U_{\bu})$ over $\C$. 

\begin{thm}\label{main} There is a 
  canonical and functorial isogeny of one-motives 
  $$ I^n(U_{\bu}) \rightarrow L^n (X_{\bu}, Y_{\bu}) $$ over $\C$, 
i.e., there is a canonical functorial
 isomorphism of $\Q$-mixed Hodge structures:
$$t^n(U_{\bu})\otimes\Q \xrightarrow{\sim} T_{\Z}(L^n)\otimes\Q.$$ 
\end{thm}

\begin{cor} $L^n(X_{\bu}, Y_{\bu})\otimes\Q$ depends only upon 
$U_{\bu}$.
\end{cor}
\begin{proof} Clear. \end{proof}{}

We follow, for the most part, Deligne's arguments in
\cite[10.3]{h}.\ssp

\nd {\bf Proof of the $W_{-1}$-part of Theorem \ref{main}.}\medskip

We begin
with the isogeny one-motive $W_{-1}L^n\otimes\Q = [ 0 \rightarrow
\tilde{\mathcal 
P}_n]\otimes\Q = [ 0 \rightarrow {\mathcal P}_n]\otimes\Q$. 

\begin{prop}\label{kona} $W_{-1}t^n(U_{\bu})\otimes\Q 
\xrightarrow{\sim} T_{\Z} (W_{-1}L^n)\otimes\Q$.\end{prop}
Since $$\mathcal P_n = \frac{\mathcal
G_n}{Pic(X_{n-2})}, \quad H^1(X_{n-2}, \Z(1))
\underset{(\ref{hdgpic})}{\zim} H_1(Pic(X_{n-2}), \Z)$$and $W_{-1} 
t^n(U_{\bu})\otimes\Q = W_{-1}H^n(U_{\bu}, \Q(1))$, this
follows from the next 
\begin{prop}\label{malli} 
{\rm{(i)}} $H':= H^1(X_{\ge n-1},\Z(1)) 
\xrightarrow{\sim} T_{\Z}([0 \rightarrow \mathcal G_n]) = H_1(\mathcal
G_n, \Z)$.

Note $t^1(X_{\ge n-1}) = H^1(X_{\ge n-1},\Z(1))$.\ssp 

\nd${\rm{(ii)}}~\frac{H^1(X_{\ge n-1},\Q(1))}{H^1(X_{n-2}, \Q(1))} \zim 
\sigma_{\ge n-1} H^n(X_{\bu},\Q(1)) = W_{-1}H^n(X_{\bu},\Q(1))
\xrightarrow{\sim} W_{-1}H^n(U_{\bu},\Q(1)).$ 
\end{prop}

(The last isomorphism is an 
analog of a theorem of Grothendieck-Deligne \cite[3.2.16-17]{h2},
\cite[9.1-9.4]{gr}.)

\begin{proof}  (ii) It follows from
  the definition of the weight filtration \cite[p.55]{ez} that the
  image of $H^n(X_{\bu}, \sigma_{\ge n-1}, \Q(1))$ in $H^n(X_{\bu},
  \Q(1))$ is $W_{-1}H^n(X_{\bu},\Q(1))$. This proves the equality.
  The rest of (ii) follows from an  inspection
  of the spectral 
  sequence  \cite[8.1.19.1]{h} and the fact (ibid. 8.1.20 (ii)) that
  it degenerates at $E_2$.  The relevant $E_1$-terms of
  \cite[8.1.19.1]{h} are those with $b= 1$ 
  and $-a = n-1$ (and $r =0$, $p=1$, $q=n-1$) since $W_{-1}H^n(U_{\bu},
  \Q(1))$ is the Tate twist of 
  $W_1H^n(U_{\bu}, \Q)$.  

(i) From (\ref{exp}), we get  the exact sequence 
\begin{equation*}
H^0(X^{an}_{\ge n-1}, \mathcal O)
  \xrightarrow{exp} H^0(X^{an}_{\ge n-1}, \mathcal O^*) \ra
 H'  \xrightarrow{\xi}  H^1(X^{an}_{\ge n-1},
 \mathcal O)  \xrightarrow{exp} H^1(X^{an}_{\ge n-1},\mathcal
 O^*).
\end{equation*}
Since the first map is surjective ($exp: \C \rightarrow \C^*$ is
surjective),  $\xi$ is injective. 
 Proposition \ref{reps} (iii) now gives the required isomorphism
$$ \boxplus: H'= H^1(X_{\ge
n-1}, \Z(1)) \xrightarrow{\sim} T_{\Z}([0 \rightarrow \mathcal G_n]):=
H_1(\mathcal G_n, \Z).$$   
Let us
show that $\boxplus$ is 
compatible with the weight and Hodge filtrations. 
In the sequences (\ref{ssfbu}), (\ref{ssylow}) on $X_{\ge 
n-1}$, we have  

\noindent $\bu$ isomorphisms of Hodge structures
$$\eta_1: E^{n,0}_2(\Q(1)) \xrightarrow{\sim}
H_1(\mathcal T',\Q),  \quad \eta_2: E^{n-1,1}_2(\Q(1))
\xrightarrow{\sim} H_1(\mathcal R, \Q).$$  

For $\eta_1$, use $H^0(X_m^{an}, \Z(1))
\xrightarrow{\sim} H_1(T_{X_m}, \Z)$, a consequence of (\ref{exp}). 
And $\eta_2$ follows from (\ref{picmd} (ii),(iv)).  

\noindent $\bu$ surjectivity of  the map $\pi$ in the commutative diagram
$$
\begin{CD} 
0 @>>> E^{n,0}_2(\Q(1)) @>>{\tau}> H'\otimes \Q @>>{\pi}>
E^{n-1,1}_2(\Q(1))
@>>> 0\\
@. @V{\eta_1}VV  @VV{\boxplus}V @V{\eta_2}VV @.\\ 
0 @>>> H_1(\mathcal T,\Q) @>>> H_1(\mathcal G_n,\Q) @>>> H_1(\mathcal
R,\Q) @>>> 0.
\end{CD}
$$
This follows from the degeneration \cite[8.1.20]{h} at $E_2$ of
(\ref{ssfbu}) for $\Q(1)$. 

The bottom row of the diagram is the the Hodge realization  of the exact
sequence (\ref{lcw}) of isogeny one-motives.  
The image of $\tau$ is $W_{-2}H'\otimes\Q$ \cite[p.55]{ez}. Since
$W_{-2}[ 0 \to \mathcal 
G_n]\otimes\Q = [0 \rightarrow 
\mathcal T]\otimes\Q$ \cite[10.1.4]{h}, we find that $\boxplus$ is compatible
with the weight filtration.  

Since $H'$ is of 
type $(*)$, there is only one nontrivial step in the Hodge filtration,
viz., $F^0$; and,  $F^0(H'\otimes\C) \xrightarrow{\sim} 
F^0((H'/{W_{-2}H'})\otimes \C)$ since $F^0\cap W_{-2} H'\otimes \C =0$. 
Thus, to show that $\boxplus$ is compatible with $F$, it 
suffices to show that $\eta_2$  is a map of Hodge structures. 
This we have done.\end{proof}{}

This finishes the proof of the $W_{-1}$-part of Theorem
\ref{main}.\ssp 

\nd {\bf Interpretation of $H^1$ and its applications.}

\begin{obser}\label{haccha1} (Interpretation of $H^1$) 
Let $d:F \ra G$
be a morphism of abelian sheaves on a space $Z$. 
In \cite[10.3.10]{h}, Deligne notes that $\hi^1(Z, [F
\xrightarrow{d} G])$ can be identified with the set of isomorphism
classes of pairs $(L, \alpha)$ where $L$ is a $F$-torsor and $\alpha$
is a trivialization of the $G$-torsor $dL$. This identification is
based upon the sequence 
\begin{equation*}
H^0(Z, F) \to H^0(Z, G) \ra \hi^1(Z,[F\xrightarrow{d} G]) \ra H^1(Z,F)
\xrightarrow{d} H^1(Z,G)\qed\end{equation*} 
\end{obser} 

Let 
$j^m_* \mathcal O^*_{U}$ \cite[10.3.9(i)]{h} be the subsheaf
of meromorphic sections of $j_*\mathcal O^*_{U}$ on
$X^{an}_{\bu}$. 

\begin{rem}\label{curie} One can prove (\ref{malli}(i)) using 
(\ref{haccha1}); namely, the group $H^1(X_{\ge n-1}, \Z(1)) \zim \hi^1(X_{\ge
n-1}, [\mathcal O \xrightarrow{exp} \mathcal O^*])$ is the set of pairs $(L,
\alpha)$ where $L$ is an $\mathcal O_{X_{\ge n-1}}$-torsor and
$\alpha$ a trivialization of $exp(L)$. Since ${\rm Aut}(L) \to {\rm
Aut}(exp(L))$ 
is onto, $H^1(X_{\ge n-1}, \Z(1))$ is the set of isomorphism
classes $L$ of $\mathcal O_{X_{\ge n-1}}$-torsors (= elements of
Lie~$\mathcal G_n$) with $exp(L) = 0$ in $\mathcal G_n$; thus,
$H^1(X_{\ge n-1}, \Z(1)) \zim H_1(\mathcal G_n, \Z)$.\ssp 

Similar results, based on (\ref{haccha1}), are as follows.\medskip 

\nd (a) the mixed Hodge structure 
$H^1(U_{m},\Z(1))$ is isomorphic to $T_{\Z}[\mathcal
V_{m} \to Pic(X_{m})]$ (\ref{toto}).\medskip

\nd (b)   the mixed Hodge structure 
$H^1(U_{\ge n-1},\Z(1))$ is isomorphic to 
$T_{\Z}[K^0 \to \mathcal G_n]$; note that $[K^0
\to \mathcal G_1]$ ($n=1$) is the Picard one-motive  \cite{ra} of
$U_{\bu}$.\ssp

We refer to \cite{ra} for the details; a sketch of the 
proof of (a) is as follows:    
$H^1(U_{m},\Z(1)) = \hi^1(X_{m}, [\mathcal O_{X}
\xrightarrow{exp} j^m_*\mathcal O^*_{U}])$ can be identified --- as in
\cite[10.3.10c]{h} --- with the set of
isomorphism classes of pairs $(L, \alpha)$ where $L$ is an $\mathcal
O_{X_{m}}$-torsor and $\alpha$ is an isomorphism of the invertible
sheaf $exp(L)$  with $\mathcal O_{X_{m}}(E)$
($E$ is a divisor supported on $Y_{m}$).  Since ${\rm Aut}(L)
\ra {\rm Aut}(exp(L))$ is surjective, 
$H^1(U_{m},\Z(1))$ is the set of pairs $(p,d)$ where $p$ is
an isomorphism class of an  $\mathcal
O_{X_{m}}$-torsor, i.e., an element of Lie~$Pic(X_{m})$, and
$d\in \mathcal V_{m}$ with $exp(p)$ as image in
$Pic(X_{m})$. This gives an isomorphism (as abelian groups) 
$$ H^1(U_{m},\Z(1)) \xrightarrow{\sim} T_{\Z}[\mathcal V_{m} \to
Pic(X_{m})].$$(Combining (a) and (b) yields --- see (\ref{hosadu}) for the definition of
$J^n$ ---  \begin{equation}\label{k3}T_{\Z}J^n\otimes\C \zim
\frac{H^1(U_{\ge n-1}, \C)}{H^1(U_{n-2},
\C)}\underset{(\ref{patagon})}{\zim} \sigma_{\ge 
n-1}H^n(U_{\bu}, \C)\underset{(\ref{ohoh!})}{\xleftarrow{\sim}}
{\rm Lie} \mathcal P^{\diamondsuit}_n.)\end{equation} 

\nd (c) $ \mathcal
G^{\diamondsuit}_n(\C) \zim \hi^1(X_{\ge n-1}, \Gamma^*)$ (\ref{gumma}) 
is the group of isomorphism classes of pairs $(\mathcal
L, \omega)$ with $\mathcal L$ an invertible sheaf on $X_{\ge n-1}$ and 
 $\omega \in H^0(X_{\ge n-1}, \Omega^1(log~Y))$. To relate to
\cite[10.3.10a]{h}, one uses the  ``connections $\theta$ on
invertible sheaves = one-forms $\omega$'' dictionary \cite[2.5]{mes}
\cite[1.5, p.47]{ev}. Thus, $\hi^1(X_{\ge n-1}, \Gamma^*)$
(cf. \cite[7.2.1]{katzdiff}) is the set of
isomorphism classes of pairs $(\mathcal L, \theta)$ where $\mathcal L$
is as before and  $\theta$ a
 connection on $\mathcal L$, holomorphic on $U_{\ge
n-1}$ and allowed 
to have simple poles along $Y_{\ge n-1}$.\ssp 

\nd (d) ${\rm Lie}~ \mathcal
G^{\diamondsuit}_n(\C) \zim \hi^1(X_{\ge n-1},
\Gamma) \zim H^1_{DR}(U_{\ge n-1})$ (\ref{gumma}) is the set of
isomorphism classes of pairs 
$(L, \theta)$ where $L$ is an $\mathcal O_{X_{\ge n-1}}$-torsor and
$\theta$ a connection on $L$ as in (c).\ssp 

\nd (e) Similar results hold for ${\rm
Lie}~Pic^{\diamondsuit}_{X_{m}}(\C)  \zim \hi^1(X_{m},\Gamma) \zim 
H^1_{DR}(U_{m})$ (\ref{diksha})  and  
$Pic^{\diamondsuit}_{X_{m}}(\C) \zim \hi^1(X_{m}, \Gamma^*)$
(\ref{diksha}).\ssp   
  
\nd (f) Parts (d) and (e) yield an interpretation of  
${\rm Lie}~\mathcal P^{\diamondsuit}_n {\xrightarrow{\sim}} 
\sigma_{\ge n-1} \hi^n(U_{\bu},\Gamma)$ (\ref{ohoh!}). \qed  
\end{rem}

While $W_{-1}t^n(U_{\bu})\otimes\Q$ is a quotient of $H^1(X_{\ge
n-1}, \Q(1)))$, $t^n(U_{\bu})\otimes\Q$ is a subquotient of $H^2(U_{\ge
n-2}, \Q(1))$. Thus, more work is necessary to complete the proof of
(\ref{main}).\ssp    
       
\nd {\bf Proof of Theorem \ref{main}.}\ssp

This will be accomplished in the following three steps.\ssp 

$\bu$ {\bf Step 1.} Construction of a certain mixed Hodge structure $h^2_X$.

$\bu$ {\bf Step 2.} Relating $h^2_X$ and $t^n(U_{\bu})$.

$\bu$ {\bf Step 3.} Interpretation of $h^2_X$ using
(\ref{haccha1}).\ssp  

\nd {\bf Step 1.} \emph{Construction of $h^2_X$.}\ssp
    
This uses the truncated complex $\tau_{\le 1} Rj_*\Z(1)_U$
\cite[1.4.6]{h2}. As before, let
$q:\tilde{Y}_m \to Y_m \to X_m$ denote the natural map from the
normalization $\tilde{Y}_m$ of $Y_m$. Since, on
each $X_m$, 
$R^1j_{m,*}\Z(1) = q_*\Z_{\tilde{Y}_m}$ \cite[3.1.9]{h2}, we 
get the triangle in the derived category of sheaves on $X^{an}_{m}$ 
\begin{equation*}\Z(1)_X \to  \tau_{\le 1} Rj_{m,*}\Z(1) \to q_* \Z[-1] \to \Z_X(1)[1]
\to \cdots;\end{equation*}the truncated complexes $\tau_{\le 1} Rj_{m,*}\Z(1)$ for
each $m$ combine to give a complex on $X_{\ge n-2}$ which we denote by
$\tau_{\le 1} Rj_*\Z(1)$. Thus, in the exact sequence on $X_{\ge n-2}$  
\begin{equation*} 0 \to \sigma_{\ge n-1}\Z(1)_X \to  \tau_{\le 1} Rj_*\Z(1) \to
\mathcal F \to 0,\end{equation*} 
$\mathcal F$ on $X_m$ is (quasi-isomorphic to) the complex $q_*
\Z[-1]$ for $m > n-2$ and $q_*\Z[-1] \to \Z_{X_{n-2}}(1)[1]$ for $m
=n-2$. 
The associated cohomology sequence (here we use that $H^i(X_{\ge n-2},
\sigma_{\ge n-1}\Z(1)) =  H^{i-1}(X_{\ge n-1}, \Z(1))$)    
\begin{equation}\label{c2}
\frac{H^1(X_{\ge n-1}, \Z(1))}{H^1(X_{n-2}, \Z(1))} \hra \hi^2(X_{\ge
n-2}, \tau_{\le 1} Rj_*\Z(1))  \to \hi^2(X_{\ge n-2}, \mathcal F) \xrightarrow{\delta^*}
   H^2(X_{\ge n-1}, \Z(1)) 
\end{equation}
fits into the commutative
diagram (defining $h^2_X$ by pullback via $\nu$)
$$
\begin{CD}
\frac{t^1(X_{\ge n-1})}{t^1(X_{n-2})} @>>> h^2_X @>>>
\frac{K \oplus NS_{n-2}}{W_{n-2}}\\
@V{\wr}V{(\ref{malli})}V @VVV @V{\nu}VV\\
\frac{H^1(X_{\ge n-1}, \Z(1))}{H^1(X_{n-2}, \Z(1))} @>>> \hi^2(X_{\ge
n-2}, \tau_{\le 1} Rj_*\Z(1)) @>>>  
\frac{H^0(\tilde{Y}_{\ge n-1}, \Z) \oplus H^2(X_{n-2}, \Z(1))}{H^0(\tilde{Y}_{n-2},
\Z)}\\
@VVV @VVV @VVV \\
\frac{\hi^2(X_{\ge n-2}, [\sigma_{\ge n-1} \mathcal O \xrightarrow{exp}
j_*^m\mathcal O^*_U])}{H^1(X_{n-2}, \mathcal O)} @>>>  \hi^2(X_{\ge
n-2}, [\mathcal O \xrightarrow{exp}
j_*^m\mathcal O^*_U]) @>>> H^2(X_{n-2}, \mathcal O)\\
@VVV @V{a}VV @|\\
\frac{\hi^2(X_{\ge n-2}, \mathcal K)}{H^1(X_{n-2}, \mathcal O)} @>>>   
 \hi^2(X_{\ge n-2}, \Gamma) @>>> H^2(X_{n-2}, \mathcal O)\\
@V{b}VV @. @.\\
\frac{H^1(X_{\ge n-1}, \mathcal O)}{H^1(X_{n-2}, \mathcal O)}
@<{\sim}<{(\ref{koil})}< {\rm Lie}~\mathcal P_n @. (**)\\
\end{CD} 
$$ 
\begin{rem}\label{penang}
\nd (i) The map $\nu$, induced by (\ref{h11}) and the isomorphism
$H^0(\tilde{Y}_i, \Z) \cong W_{Y_i}$, is \emph{injective}.

\nd (ii) $h^2_X$ is a mixed Hodge structure of type $(*)$; note  
$W_{-1}h^2_X = \frac{H^1(X_{\ge n-1}, \Z(1))}{H^1(X_{n-2}, \Z(1))}$. 

(An argument similar to that in the proof of (\ref{kapadu}(i)) shows that
the natural map $h^2_X\otimes\Q \to H^2(U_{\ge n-2}, \Q(1))$ is 
injective with image $t^2(U_{\ge n-2})\otimes\Q$.)     

\nd (iii) The map $a$ is induced by the (first) morphism of complexes 
$$
\begin{CD}
\mathcal O_{X} @>{exp}>> j^m_*\mathcal O^*_{U} @>>> 0\\
 @| @VV{d~{\rm log}}V @.\\
\mathcal O_{X} @>{d}>> \Omega^1(log~ Y) @>>> 0\\
 @| @| @.\\
 \mathcal O_{X} @>{d}>> \Omega^1(log~ Y) @>{d}>>
\Omega^2(log~Y) \rightarrow \cdots; \\
\end{CD}
$$it is compatible with the natural map\begin{equation*}\hi^2(X_{\ge
n-2},  Rj_*\Z(1)) \zim H^2(U_{\ge n-2}, \Z(1)) \hra \hi^2(X_{\ge n-2},
\Omega(log~Y)) \xrightarrow{\sim} H^2(U_{\ge
n-2},\C).\end{equation*}The last isomorphism is the logarithmic
Poincar\'e lemma \cite[3.2.2]{h2}.     

\nd (iv) On $X_{\ge n-1}$, $\mathcal K = \Gamma$ and,
on $X_{n-2}$, $\mathcal K = [ 0 \to \Omega^1(log~Y)]$. 

\nd (v) Since $\mathcal C'_n = \frac{K \oplus NS_{n-2}}{W_{n-2}}$ is
of type $(0,0)$, the map
$h^2_X \to H^2(X_{n-2}, \mathcal O)$ in $(**)$ is
zero. We deduce an injection $h^2_X \hra \frac{\hi^2(X_{\ge n-2}, [\sigma_{\ge n-1}
\mathcal O \xrightarrow{exp} 
j_*^m\mathcal O^*_U])}{H^1(X_{n-2}, \mathcal O)}$  (with finite
cokernel --- see (\ref{devaki}) below) and a
commutative diagram   
$$\begin{CD}
h^2_X @>>> h^2_X\otimes\C @<<< W_{-1}h^2_X\otimes\C\\
@VVV @VVV @V{\wr}V{(ii)}V\\
\frac{\hi^2(X_{\ge n-2}, [\sigma_{\ge n-1}
\mathcal O \xrightarrow{exp} j_*^m\mathcal O^*_U])}{H^1(X_{n-2},
\mathcal O)} @>{a}>> \frac{\hi^2(X_{\ge n-2}, \mathcal
K)}{H^1(X_{n-2},\mathcal O)} @<<< \frac{\hi^1(X_{\ge n-1},
\Gamma_0)}{\hi^1(X_{n-2}, \Gamma_0)}\\
@VVV @VV{b}V @V{(\ref{ohoh!})}V{\wr}V\\
\frac{H^1(X_{\ge
n-1}, \mathcal O)}{H^1(X_{n-2}, \mathcal O)} 
@<{\sim}<{(\ref{koil})}< {\rm Lie}~\mathcal P_n @<<< 
{\rm Lie}~\mathcal P^{\natural}_n.\\ 
\end{CD}
$$
\nd (vi) The Hodge filtration \cite[3.2.2]{h2}, restricted
to $\mathcal K$, is $F^0 \mathcal K = \mathcal K$ and $F^1\mathcal K =
[0 \to \Omega^1(log~Y)]$ on $X_{\ge n-2}$. It induces the map
$b$ --- see $(**)$ --- thereby defining the Hodge filtration
on $h^2_X\otimes\C$.  Because $h^2_X$ is defined using cohomology with
$\Z(1)$-coefficients, $F^1$ on $\mathcal K$ corresponds to $F^0$ on
$h^2_X\otimes\C$.  
\qed\end{rem}

\begin{lem}\label{samos} One has 
 a commutative diagram
$$
\begin{CD}
\mathcal C'_n @>{\rho}>> {\mathcal P}'_n @>>> \pi_0(\mathcal P'_n)\\ 
@V{\nu}VV @. @A{\wr}AA\\ 
\frac{H^0(\tilde{Y}_{\ge n-1}, \Z) \oplus H^2(X_{n-2}, \Z(1))}{H^0(\tilde{Y}_{n-2},
\Z)} @>{\delta^*}>{(\ref{c2})}> H^2(X_{\ge n-1}, \Z(1)) @<<<
\pi_0(\mathcal G'_n).\\
\end{CD}
$$
\end{lem}

\begin{proof} Straightforward. \end{proof}{}

Combining (\ref{samos}), (\ref{kona}) and (\ref{malli}), we obtain the
exact sequence  
\begin{equation}\label{kvk} 0 \to  H_1(\mathcal P_n, \Q) \to h^2_X\otimes\Q \to
\mathcal C_n\otimes\Q  \to 0.\end{equation}

\nd {\bf Step 2.} \emph{Relating $h^2_X$ and $t^n(U_{\bu})$.}\ssp

Using the next lemma, we shall show that $t^n(U_{\bu})$ is a quotient
of $h^2_X$ by $NS_{n-3}$. 

\begin{lem}\label{kapadu} {{\rm (i)}} One has an isomorphism $\nu:\mathcal
B_n\otimes\Q  
\zim Gr^W_0 t^n(U_{\bu})\otimes \Q.$

{{\rm (ii)}} $\mathcal C_n$ has finite index in $\mathcal A_n$;
cf. {\rm (\ref{auro})}.   

{{\rm (iii)}} $K^0:= \Ker(K \xrightarrow{\vartheta}
\pi_0(\mathcal G'_n))$ has finite index in $K \cap \mathcal V_{n-1}:=
\Ker(K \xrightarrow{\lambda_{n-1}} NS_{n-1})$. 
\end{lem}

\begin{proof} (i) By \cite[8.1.19.1]{h}, $Gr^W_{b}H^n(U_{\bu},\Q)$ is
a subquotient of $E_1^{n-b, b}$. Since $Gr^W_0 t^n(U_{\bu})\otimes\Q$ is a
Tate twist of the ($1,1$)-part of $Gr^W_2 H^n(U_{\bu},\Q)$, the
relevant $E_1$-terms correspond to $b=2$ and $-a= n-2$ (with
$p,q,r$ satisfying $p+2r = 2$ and $q-r = n-2$). If written
explicitly  (as is
done immediately after 8.1.19.1 in \cite{h})\footnote{We note a
harmless typo in (loc. cit): Gysin maps go from $H^{b-2r}$ to 
$H^{b-2(r-1)}$ and not $H^{b-2(r-2)}$.}, we
obtain --- using the degeneration \cite[8.1.20 (ii)]{h} at $E_2$ of
\cite[8.1.19.1]{h} --- that  $Gr^W_2 H^n(U_{\bu},\Q)$ is (the
coefficients in (\ref{soya}) are $\Q$)    
\begin{equation}\label{soya}
\frac{\Ker(H^0(\tilde{Y}_{n-1})(-1) \oplus H^2(X_{n-2})
\xrightarrow{t_{n-1}}  H^0(\tilde{Y}_{n})(-1) \oplus
H^2(X_{n-1}))}{\im(H^0(\tilde{Y}_{n-2})(-1) \oplus H^2(X_{n-3}) 
\xrightarrow{t_{n-2}}  H^0(\tilde{Y}_{n-1})(-1) \oplus H^2(X_{n-2}))};\end{equation} 
here $t_{m}(a,b) = (\delta^*_{m} a, \delta^*_{m-1}b
-\lambda_m(a))$. This is easily compared with (\ref{iwas}).   

By \cite[8.1.20 (ii)]{h}, the group in (\ref{soya}) is unchanged if
one replaces $H^2(X_{n-1})$ by $H^2(X_{\ge n-1})$ and uses the map
$\delta^*$ of (\ref{c2}).   
The $(1,1)$-part of  $Gr^W_2H^n(U_{\bu},\Q)$ can be identified using
(\ref{h11}). The lemma now follows from $\pi_0(\mathcal P'_n)
\cong \pi_0(\mathcal G'_n)$,  $\pi_0(\mathcal G'_n) \hra
H^2(X_{\ge n-1}, \Z(1))$  ---  a consequence of
(\ref{exp}), and (\ref{samos}).

(ii) and (iii) also follow from the degeneration \cite[8.1.20 (ii)]{h}
at $E_2$ of \cite[8.1.19.1]{h}. \end{proof}

Using (\ref{kapadu}(i)) and (\ref{kona}), we can rewrite 
\begin{equation}\label{palani}
0 \to W_{-1}t^n(U_{\bu})  \to t^n(U_{\bu}) \to Gr^W_0
t^n(U_{\bu}) \to 0\end{equation} 
as the 
sequence (exact modulo finite groups)  of mixed Hodge structures
\begin{equation}\label{horanadu}
 0 \to H_1(\mathcal P_n, \Z) \to
t^n(U_{\bu}) \to \mathcal B_n \to 0. \end{equation}

\begin{prop}We have an isomorphism of mixed Hodge structures
\begin{equation}\label{nyaya}
\frac{h^2_X\otimes\Q}{\delta^*NS(X_{n-3})\otimes\Q} \zim
t^n(U_{\bu})\otimes\Q.\end{equation}\end{prop}  

\begin{proof} Since $h^2_X$ is of type $(*)$, the natural map $$H^2(X_{\ge
n-2}, \tau_{\le 1} Rj_*\Z(1)) \to H^2(U_{\ge
n-2},\Z(1)) \to H^n(U_{\bu}, \Z(1))$$gives a map $h^2_X \to
t^n(U_{\bu})$. The proposition follows from 
 (\ref{horanadu}), (\ref{bareesh}), (\ref{kapadu}), and
(\ref{kvk}). \end{proof}   

\nd {\bf Step 3.} \emph{Interpretation of $h^2_X$ using}
(\ref{haccha1}).\ssp

\begin{thm} One has a natural isomorphism $$\frac{\hi^2(X_{\ge n-2}, [\sigma_{\ge n-1}
\mathcal O \xrightarrow{exp} 
j_*^m\mathcal O^*_U])}{H^1(X_{n-2}, \mathcal O)}\zim
T_{\Z}[\mathcal C_n \xrightarrow{\rho} \mathcal P_n].$$\end{thm}  
  
\begin{rem}\label{merext}     
Every $j_*^m\mathcal
O^*_{U_{r}}$-torsor  $P$ gives an invertible sheaf on $U_{r}$ and this
extends to an invertible sheaf $P'$ on $X_{r}$; given any two such
extensions $P'$ and $P''$, one has an isomorphism of invertible
sheaves $P'' \zim P' \otimes \mathcal O(v)$ ($v\in W_{Y_{r}}$) on
$X_{r}$.\qed \end{rem}

\begin{proof}
The group $\hi^2(X_{\ge n-2}, [\sigma_{\ge n-1}
\mathcal O \xrightarrow{exp} 
j_*^m\mathcal O^*_U])$ is actually a $\hi^1$ in disguise 
(this shift occurs for reasons of degree for $j_*^m\mathcal O^*_U$ --- 
this is in degree one --- and of truncation for $\mathcal O$): 
it sits in an exact sequence 
\begin{multline*}
\to \hi^2(X_{\ge n-2}, [\sigma_{\ge n-1}
\mathcal O \xrightarrow{exp} 
j_*^m\mathcal O^*_U]) \to H^1(X_{\ge n-1}, \mathcal O_X)\oplus
H^1(X_{n-2}, j^m_*\mathcal O^*_U) \to\\ \to H^1(X_{\ge n-1},
j^m_*\mathcal O^*_U) \to \end{multline*}
Via (\ref{haccha1}),
we can interpret $\hi^2(X_{\ge n-2}, [\sigma_{\ge n-1}
\mathcal O \xrightarrow{exp} 
j_*^m\mathcal O^*_U])$ as the set of isomorphism classes of triples 
$(L, P, \alpha)$ where $L$ is an $\mathcal O_{X_{\ge n-1}}$-torsor and $P$ is
a $j_*^m\mathcal O^*_{U_{n-2}}$-torsor and $\alpha$ is an
isomorphism\footnote{We write $exp_m(L)$ for the $j_*^m\mathcal O^*_{U_{\ge
n-1}}$-torsor and $exp(L)$ for the $\mathcal O^*_X$-torsor associated
with $L$.}  
$exp_m(L) \zim \delta^*P$ of $j_*^m\mathcal O^*_{U_{\ge
n-1}}$-torsors. Via (\ref{merext}),  $\alpha$ can also be thought of
as  an isomorphism 
(as in \cite[10.3.10c]{h})    
$exp(L) \zim \delta^*P'\otimes \mathcal O(D')$ of invertible sheaves 
on $X_{\ge n-1}$ with $D' \in K$. If $P''$ is another invertible
sheaf on $X_{n-2}$ corresponding to $P$, then we have $\alpha'':exp(L)
\zim \delta^* P''\otimes \mathcal O(D'')$ and $P'' \zim P' \otimes
\mathcal O(v)$ ($v\in W_{n-2}$); so one can rewrite $\alpha$ as  
$exp(L) \zim \delta^*(\mathcal O(-v)\otimes P'') \otimes
\mathcal O(D''+ \delta^*v)$. As $(P',D')- (P'',
D'') = \gamma_{n-3}(v) \in NS_{n-2}\oplus K$,  the element $w=
(P', D') \in \mathcal C'_n$ depends only on the isomorphism 
class of the triple $(L, P, \alpha)$. Thus, we may associate the
pair $(w, [L])\in \mathcal C'_n\oplus {\rm Lie}~\mathcal G_n$ with the
isomorphism class of $(L, P, \alpha)$.  

The map $\delta^*:H^1(X_{n-2}, \mathcal O) \to \hi^2(X_{\ge n-2},
[\sigma_{\ge n-1} 
\mathcal O \xrightarrow{exp} 
j_*^m\mathcal O^*_U])$ can be described as follows. 
The class $[I]\in H^1(X_{n-2}, \mathcal O)$ of an $\mathcal
O_{X_{n-2}}$-torsor $I$ is mapped to the class of the triple 
 $(\delta^* I, P= exp_m(I),
\alpha_I)$ with $\alpha_I: exp_m(\delta^*I) \zim \delta^*(exp_m(I))$
the tautological 
isomorphism; associated with this triple is the invertible sheaf $P' =exp(I)$
(\ref{merext}), the element $D =0$ of $K$ and the isomorphism  
$\alpha_I: exp(\delta^*I) \zim
\delta^*(exp (I))$. If $P''$ is another invertible sheaf 
corresponding to the $j_*^m\mathcal O^*_U$-torsor $exp_m(I)$, then as
before $P'' \zim P'\otimes \mathcal O(v)$ ($v\in W_{n-2}$); thus,
to the triple $(\delta^* I, exp_m(I),
\alpha_I)$, we may attach the pair $(0,[\delta^*I]) \in \mathcal
C'_n\oplus {\rm Lie}~\mathcal G_n$.    

Taking into account  the
surjectivity of $H^0(X_{\ge n-1}, \mathcal O) \xrightarrow{exp} H^0(X_{\ge n-1},
\mathcal O^*)$, we find that elements of 
$$\frac{\hi^2(X_{\ge n-2}, [\sigma_{\ge n-1}
\mathcal O \xrightarrow{exp} 
j_*^m\mathcal O^*_U])}{H^1(X_{n-2}, \mathcal O)}$$ can be identified with 
 pairs $(w, L)$ where $L$ is an element of Lie~$\mathcal P_n$ 
 and (as before) $w \in \mathcal C'_n$ with 
$exp(L) = \rho(w)$ in $\mathcal
P_n$. This last equality forces, by (\ref{samos}),  $w$ to be  an
element of $\mathcal C_n$. Recalling 
the definition of $T_{\Z}$ \cite[10.1.3.1]{h}, we deduce the required 
isomorphism.\end{proof} 
Consider the composite
inclusion \begin{equation}\label{devaki}
h^2_X \hra \frac{\hi^2(X_{\ge n-2}, [\sigma_{\ge n-1}
\mathcal O \xrightarrow{exp} 
j_*^m\mathcal O^*_U])}{H^1(X_{n-2}, \mathcal O)}~ {\zim}~
T_{\Z}[\mathcal C_n \xrightarrow{\rho} \mathcal
P_n];\end{equation}$h^2_X\otimes\Q$ (resp. $T_{\Z}[\mathcal C_n
\xrightarrow{\rho} \mathcal P_n]\otimes\Q$)   is an extension 
 (\ref{kvk}) (resp. \cite[10.1.3.1]{h}) of $\mathcal
C_n\otimes\Q$ by $H_1(\mathcal P_n, \Q)$. Thus, we obtain  
 an isomorphism of abelian
groups $$h^2_X \xrightarrow{\sim}~
T_{\Z}[\mathcal C_n \xrightarrow{\rho} \mathcal P_n] \quad
\textrm{(modulo finite groups)}.$$From (\ref{bareesh}) and (\ref{nyaya}), we deduce  
an isomorphism  of abelian groups 
\begin{equation*}
\Lambda: t^n(U_{\bu})
\xrightarrow{\sim}~ T_{\Z}L^n \quad \textrm{(modulo finite groups)}.
\end{equation*}

\nd {\bf Compatibilities of $\Lambda$.}\medskip

$\bu$ \emph{Weight filtration}: Since $NS(X_{n-3})\otimes\Q = 
t^2(X_{n-3})\otimes\Q$ is of type $(0,0)$,
from (\ref{nyaya}) we  have $W_{-1}h^2_X\otimes\Q \xrightarrow{\sim}
W_{-1}t^n(U_{\bu})\otimes\Q$.
The elements of $W_{-1}h^2_X$ are characterized, by
(\ref{palani}), (\ref{horanadu}) and (\ref{curie}), as
those corresponding to pairs  $(L, \beta)$ with $exp(L) = 0
=\beta$. This proves the compatibility of $\Lambda$ with the weight
filtration.\medskip

$\bu$ \emph{Hodge filtration}:
From \cite[10.1.3.1]{h}, the map $\alpha: T_{\Z}L^n \ra {\rm
Lie}~\mathcal P_n$ (used to construct $T_{\Z}L^n$) gives the
Hodge filtration. Namely, $$F^0(T_{\Z}L^n \otimes \C)  = \Ker(\alpha_{\C}:
(T_{\Z}L^n)\otimes\C   \ra {\rm
Lie}~\mathcal P_n).$$ Therefore, we obtain $$\alpha_{\C}: \frac{(T_{\Z}L^n)\otimes\C}{F^0} \xrightarrow{\sim} {\rm Lie}~\mathcal P_n.$$    
Since $t^n(U_{\bu})$ is of type $(*)$, we get \cite[10.1.3.3]{h}
the isomorphism  
$$\frac{t^n(U_{\bu})\otimes\C}{F^0} \xleftarrow{\sim}
\frac{W_{-1}t^n(U_{\bu})\otimes\C}{F^0 \cap W_{-1}}.$$ 
Now (\ref{ohoh!}) and (\ref{kona})  together imply that
$$W_{-1} t^n(U_{\bu})\otimes\C \xrightarrow{\sim} \frac{\hi^1(X_{\ge n-1}, \Gamma_0)}{\hi^1(X_{n-2}, \Gamma_0)} \xleftarrow{\sim} {\rm Lie}~\mathcal
P^{\natural}_n.$$   
Noting that the Hodge filtration on $H^*(U_{\bu}, \C)$ is induced by
the filtration (\cite[3.2.2]{h2}, \cite[8.1.8, 8.1.12]{h})
$$F^i \Omega(log~Y):=  0 \to 0 \to \cdots \to \Omega^i(log~Y)
\xrightarrow{d} \Omega^{i+1}(log~Y) \to \qquad ,
$$we obtain --- cf. (\ref{penang} (v), (vi)) ---  
the commutativity of
the following diagram  
$$ 
\begin{CD}
{} @. W_{-1}t^n(U_{\bu})\otimes\C @<{\sim}<< {\rm Lie}~\mathcal
 P^{\natural}_n @>{\sim}>{(\ref{ohoh!})}> \frac{\hi^1(X_{\ge n-1}, \Gamma_0)}{\hi^1(X_{n-2}, \Gamma_0)}\\
@. @VVV @VVV @VVV\\
\frac{t^n(U_{\bu})\otimes\C}{F^0} @<{\sim}<< \frac{W_{-1}t^n(U_{\bu})\otimes\C}{F^0 \cap W_{-1}} @<{\sim}<< {\rm
 Lie}~\mathcal P_n @>{\sim}>{(\ref{koil})}> \frac{H^1(X_{\ge n-1}, \mathcal
 O)}{H^1(X_{n-2}, \mathcal O)},\\
\end{CD}
$$
which fits into a larger commutative 
diagram 
$$
\begin{CD}
H^n(U_{\bu},\Z(1))/{\rm torsion} @<<< t^n(U_{\bu}) @>{\Lambda}>>
T_{\Z}L^n @>{\alpha}>> {\rm Lie}~\mathcal P_n\\
@VVV @VVV @. @| \\
H^n(U_{\bu},\C) @<<< t^n(U_{\bu})\otimes\C @>>>
\frac{t^n(U_{\bu})\otimes\C}{F^0} @>>> {\rm Lie}~\mathcal P_n.\\
\end{CD}
$$
This diagram shows that $\Lambda$ is compatible with the Hodge
filtration thereby finishing the proof of Theorem
\ref{main}.\qed\medskip  

Theorem \ref{main} partly proves Conjecture \ref{dc1} (up to isogeny) for the
simplicial scheme $U_{\bu}$; it remains to prove the statements
concerning the  
de Rham and \'etale realizations. We first treat the de Rham
realization. The \'etale realization is dealt with in \S
\ref{oakland}.\ssp

\nd {\bf The de Rham realization of $L^n$.}\ssp

Using $T_{DR}L^n \zim T_{\C}L^n$ \cite[10.1.8]{h},  
$T_{\C}L^n \zim t^n(U_{\bu})\otimes\C$ (\ref{main}), and
$t^n(U_{\bu})\otimes\C \hra H^n(U_{\bu}, \C)  \zim H^n_{DR}(U_{\bu})$
\cite[3.1.8, 3.2.2]{h2} gives us a 
map $$T_{DR}L^n \to H^n_{DR}(U_{\bu});$$ our next task is to show that
this map can be constructed purely algebraically (\ref{derham}) and
this, over
any field of characteristic zero.  

In the remainder of this section, $k$ denotes a field of characteristic zero. Our next main
result 
(\ref{derham}) of this section requires us to construct  
a group scheme $\tilde{\mathcal U}$,  a map $\tilde{\psi}: \mathcal B_n \to
\tilde{\mathcal U}$ such that, when $k =\C$, $\tilde{\psi}$ lifts to a
compatible homomorphism $\psi: T_{\Z}L^n \to {\rm Lie}~\tilde{\mathcal
U}$.  Once
these are acquired, the criterion \cite[10.1.9]{h} may be applied to
deduce  that ${\rm Lie}~\tilde{\mathcal U}$ is the de
Rham  realization $T_{DR}L^n$ of $L^n$.\ssp

\nd {\bf Step 1.} \emph{Construction of a map} $\psi_1: K \to \mathcal
G_n^{\diamondsuit}$.\ssp

Recall the map $q:\tilde{Y}_m \to Y_m \to X_m$ from the 
normalization $\tilde{Y}_m$ of $Y_m$.  
The  Poincar\'e residue sequence \cite[3.1.5.2]{h2} on $X_m$ and $X_{\ge n-1}$  
\begin{equation*} 0 \ra \Omega^1 \rightarrow
\Omega^1(log~Y) \xrightarrow{{\rm Res}} q_*\mathcal 
O_{\tilde{Y}} \rightarrow 0
\end{equation*} gives the exact sequences
\begin{equation}\label{eq:5} 
H^0(\tilde{Y}_{m}, \mathcal O) \to H^1(X_{m}, \Omega^1)
\to H^1(X_{m}, \Omega^1(log ~Y)),\end{equation}
\begin{equation}\label{eq:2} 
H^0(\tilde{Y}_{\ge n-1}, \mathcal O) \to H^1(X_{\ge n-1}, \Omega^1)
\to H^1(X_{\ge n-1}, \Omega^1(log ~Y)).\end{equation} 
Since the first map of (\ref{eq:5}) is the composition of $$
W_{Y_m}(S)\otimes k \zim H^0(\tilde{Y}_{m}, \mathcal O) \xrightarrow{\lambda_m} NS_{X_m}(S)\otimes k
\underset{(\ref{beagle})}{\xrightarrow{c_{X_m}}} 
H^1(X_m, \Omega^1),$$ there is a natural injection --- recall $N_m:=
\Coker(\lambda_m)$ (\ref{auro2}) --- 
\begin{equation}\label{eq:7}
\kappa_m: {N_m}(S) \otimes k \hra H^1(X_{\ge n-1},
\Omega^1(log ~Y)).\end{equation}   
The boundary map of (\ref{eq:2}) is induced by the  composite map $$K(S) \xrightarrow{\vartheta'}
H^1(X_{\ge n-1}, \mathcal O^*) = \mathcal G'_n (S) \to \pi_0(\mathcal G'_n)
\underset{(\ref{dlog})}{\xrightarrow{dlog}} H^1(X_{\ge
n-1}, \Omega^1)$$via $K(S)\otimes k \zim H^0(\tilde{Y}_{\ge n-1},
\mathcal O)$.
 Thus, in the exact sequence 
$$\hi^1(X_{\ge n-1}, [\mathcal O^* \xrightarrow{d~{\rm log}}
\Omega^1(log ~Y)]) \to H^1(X_{\ge n-1}, \mathcal O^*) \to
H^1(X_{\ge n-1}, \Omega^1(log ~Y)),$$ one finds, by the exactness of
(\ref{eq:2}), 
that the map $\vartheta'$ admits a  lifting  to $\hi^1(X_{\ge n-1},
\Gamma^*)$.  

In fact, one can easily construct a natural lifting 
$\psi_1: K(S) \to \hi^1(X_{\ge n-1}, \Gamma^*) = \mathcal
G_n^{\diamondsuit}(S)$ of $\vartheta'$ using \v{C}ech cohomology; let
$C^*(G)$ denote 
\v{C}ech cochains of $X_{\ge n-1}$ with coefficients in $G$ (relative to a
suitable open cover $\{U_i\}$) and $\partial$ the \v{C}ech differential.
By (\ref{curie}c), for each element $E \in K$, we have to construct an
invertible sheaf $L$ on $X_{\ge n-1}$ and a 
connection $\nabla$ with at most logarithmic poles along $Y_{\ge
n-1}$. In terms of cocycles \cite[2.5]{mes} \cite[7.2]{katzdiff}, using $C^1(\Gamma^*) =
C^1(\mathcal O^*) \oplus C^0(\Omega^1(log~Y))$, 
if $\{s_{ij}\}\in C^1(\mathcal O^*)$ represents
$L$ and $\{\omega_i\} \in
C^0(\Omega^1(log~Y))$ represents $\nabla$, then these
satisfy $ds_{ij}/{s_{ij}} = \omega_j - \omega_i$; $\nabla$ is integrable if and
only if $\omega_i$ is closed.    

Let $E\in K$ be an effective divisor. If $\{f_i\}$ are local equations
for $E$, consider the cochain $(s_{ij}, df_i/{f_i})\in 
C^1(\mathcal O^*)\oplus
C^0(\Omega^1(log~Y))$ where $f_j =
f_i s_{ij}$; clearly, $df_j/{f_j} - df_i/{f_i} =
ds_{ij}/{s_{ij}}$. Thus, the pair $(s_{ij}, df_i/{f_i})$ is a cocycle;
and, it represents an
element of $\hi^1(X_{\ge n-1}, \Gamma^*)$. 
 If $\{g_i\}$ are different local equations for $E$, then
one gets an element $(t_{ij}, dg_i/{g_i})\in C^1(\mathcal O^*)\oplus
C^0(\Omega^1(log~Y))$ with $g_j = g_i t_{ij}$. 
There exist $u_i \in C^0(\mathcal O^*)$ with $f_i = u_i
g_i$. Since $s_{ij}/t_{ij} = u_j/u_i$ and $ds_{ij}/{s_{ij}} -
dt_{ij}/{t_{ij}} =  du_j/{u_j}  - du_i/{u_i}$, i.e., the 
cocycle $(t_{ij}, dg_i/{g_i}) -(s_{ij}, df_i/{f_i})$ is 
a coboundary ($\partial u$), the element $\psi_1(E) = (s_{ij},
df_i/{f_i})$ of  $\hi^1(X_{\ge n-1}, \Gamma^*)$ depends only on $E$ but not
on the choice of the local defining equations.   
The association $E \mapsto \psi_1(E)$ defined for effective $E$ easily
extends to a homomorphism $\psi_1: K \to \mathcal
G_n^{\diamondsuit}$.\ssp

\nd {\bf Step 2.} \emph{Construction of the scheme} $\tilde{\mathcal
U}$.\ssp 

Let $g$ be the 
structure map $X_{\ge n-2} \to S$; recall $a_m$ is the structure map
 $X_m \to S$. 
\begin{lem}\label{sec:de-rham-realization-1} The $S_{fppf}$-sheaf $R^2g_*\mathcal K^*$ 
associated with the complex $\mathcal K^*:= [\sigma_{\ge n-1} \mathcal O^*
\xrightarrow{d~{\rm log}} \Omega^1(log~Y)]$ on $X_{\ge n-2}$  is
representable. One has $\xi_*: \pi_0(R^2g_*\mathcal K^*) \hra
\pi_0(\mathcal G'_n)$; here $R^2g_*\mathcal K^*$  denotes the
associated representing group scheme.\end{lem}
\begin{proof}  One has an exact sequence  
\begin{multline*}
H^0(X_{\ge n-1},
\mathcal O^*) \to H^1(X_{\ge n-2}, \Omega^1(log~Y)) \to \hi^2(X_{\ge
n-2}, \mathcal K^*) \xrightarrow{\xi} \\
  \to H^1(X_{\ge n-1}, \mathcal O^*) \to 
H^2(X_{\ge n-2}, \Omega^1(log~Y));\end{multline*} 
the sheaf $R^1f_*\mathcal O^*$ is representable by 
$\mathcal G'_n$, $f_*\mathcal O^*$ by a torus, and 
 $R^ig_* \Omega^1(log~Y)$ are 
representable by vector group schemes. Since $\Hom(f_*\mathcal O^*,
R^1g_* \Omega^1(log~Y)) = 0$, one can now invoke (\ref{reps}) to
obtain the representability of $R^2g_*\mathcal K^*$. 

Since $\Ker(\xi)$ is connected, the map $\xi_*$ is injective.\end{proof}
The inclusion $\sigma_{\ge n-1}\Gamma^* \hra \mathcal
K^*$ induces  an exact sequence 
\begin{equation}\label{eq:6}
0 \to  \frac{\mathcal
G_n^{\diamondsuit}}{a_{n-2~*} \Omega^1(log~Y)} \to 
R^2g_*\mathcal K^* \xrightarrow{\upsilon} R^1a_{n-2~*}
\Omega^1(log~Y). 
\end{equation} 

Consider the commutative diagram 
$$ 
\begin{CD}
H^1(X_{n-3}, \Omega^1(log~Y)) @>{\delta^*}>> H^1(X_{\ge n-2},
\Omega^1(log~Y))\\
@VV{\delta^*_{n-3}}V @VVV\\
H^1(X_{n-2}, \Omega^1(log~Y)) @<{\upsilon}<< \hi^2(X_{\ge n-2},
\mathcal K^*),\\
\end{CD}
$$and the map $$Pic^0_{X_{n-2}}(S) \hra H^1(X_{n-2}, \mathcal O^*) 
\to \hi^2(X_{\ge
n-2}, \mathcal K^*)$$induced by the inclusion $\mathcal K^* \hra \Gamma^*$ on
$X_{\ge n-2}$. Taking the quotient of $\hi^2(X_{\ge
n-2}, \mathcal K^*)$ by the images of $H^1(X_{n-3}, \Omega^1(log~Y))$
and $Pic^0_{X_{n-2}}(S)$ under these maps, and pulling back via
(cf. (\ref{eq:6}) (\ref{eq:7}))  
\begin{equation}\label{eq:8}\frac{\Ker(\delta^*_{n-2}:N_{n-2}\to N_{n-1})(S)\otimes
k}{NS_{n-3}(S)\otimes k} \to 
\frac{N_{n-2}(S)\otimes
k}{NS_{n-3}(S)\otimes k} \to 
\frac{H^1(X_{n-2}, \Omega^1(log~Y))}{H^1(X_{n-3}, \Omega^1(log~Y))}\end{equation}
(actually we do this at the level of the group schemes which represent
the associated sheaves) gives us a group scheme $\mathcal U'$. Its
identity component is denoted $\mathcal U$.  

 \begin{lem}\label{sec:de-rham-realization} The map $\xi_*:\pi_0(\mathcal U') 
\hra \pi_0(\mathcal P'_n)$ is injective.\end{lem} 
\begin{proof} It is
clear that the previous operations of quotient and pullback do not affect
$\pi_0$ and thus $\pi_0(R^2g_*\mathcal K^*) \cong \pi_0(\mathcal
U')$. The natural map $\xi: R^2 g_*\mathcal K^* \to \mathcal G'_n$
provides a map $\xi: \mathcal U' \to \mathcal P'_n$; by
(\ref{sec:de-rham-realization-1}), one has an injection $\xi_*:
\pi_0(\mathcal U') \hra \pi_0(\mathcal P'_n)$.\end{proof} 

\begin{lem} One has an exact sequence
\begin{equation}\label{eq:9} 0 \to {\rm Lie}~\mathcal
P^{\diamondsuit}_n \to   {\rm Lie}~\mathcal U \xrightarrow{\upsilon}
\frac{\Ker(\delta^*_{n-2}:N_{n-2}\to N_{n-1})(S)\otimes 
k}{NS_{n-3}(S)\otimes k}.\end{equation}\end{lem} 
\begin{proof} 
By (\ref{diksha}), ${\rm Lie}~\mathcal P^{\diamondsuit}_n$  
is a quotient of ${\rm Lie}~\mathcal G^{\diamondsuit}_n$ by
${\rm Lie}~Pic^{\diamondsuit}_{X_{n-2}}$.  Since ${\rm
Lie}~Pic^{\diamondsuit}_{X_{n-2}}$ is an extension of ${\rm Lie}~Pic^0_{X_{n-2}}$ by
$a_{n-2~ *}\Omega^1(log ~Y)$, we obtain (\ref{eq:9}) noting
(\ref{eq:6}). One just has to observe that the image of $Pic^0_{X_{n-2}}$ in
$R^2g_*\mathcal K^*$ is equal to that of the neutral
component of $Pic^{\diamondsuit}_{X_{n-2}}$ under the natural map
$ Pic^{\diamondsuit}_{X_{n-2}} \xrightarrow{\delta^*} \mathcal
G^{\diamondsuit}_n \to R^2g_*\mathcal K^*$.\end{proof} 

Since the composite map $$ H^1(X_{n-2}, \mathcal O^*) \to \hi^2(X_{\ge
n-2}, \mathcal K^*) \xrightarrow{\upsilon}H^1(X_{n-2},
\Omega^1(log~Y))$$is the map induced by $d~{\rm log}: \mathcal O^* \to
\Omega^1 \to \Omega^1(log~Y)$, the natural map $NS_{n-2} \to
\frac{R^2g_*\mathcal K^*}{Pic^0_{n-2}}$induces a map $$\psi_2:NS_{n-2} \to
\mathcal U'.$$By
definition (\ref{eq:8}) of $\mathcal U'$, 
$\im(\delta^*_{n-3}: NS_{n-3} \to NS_{n-2})$ is contained in 
$\Ker({\psi_2})$.

\begin{lem} The map $\psi': K \oplus NS_{n-2}
\to \mathcal U'$ defined by $(u,v) \mapsto \psi_1(u) + \psi_2(v)$
provides a map  $\psi':\mathcal B'_n
\to \mathcal U$.\end{lem}   
\begin{proof}  Since $\psi' = \psi_2$ on  
the subgroup $\gamma_{n-3}(NS_{n-3}) = \delta^*_{n-3}(NS_{n-3})$, 
$NS_{n-3}$ is contained in $\Ker({\psi'})$. 
It is straightforward to check that $\gamma_{n-3}(W_{Y_{n-2}})$
is in  $\Ker({\psi'})$. This gives a map $\psi': \mathcal B'_n \to
\mathcal U'$. Note  that
the composite map $$ K \oplus 
NS_{n-2} \xrightarrow{\psi'} \mathcal U' \xrightarrow{\xi} \mathcal
P'_n$$ is  $\rho$. By 
(\ref{sec:de-rham-realization}), the image under $\psi'$
of $\mathcal B'_n$ is contained in $\mathcal U$.\end{proof} 

Put $\tilde{\mathcal U} = \mathcal U/{\psi'(\tau_n)}$; here, 
$\tau_n$ is the torsion subgroup of $\mathcal B'_n$. The map $\psi'$
induces  a map $\tilde{\psi}: \mathcal B_n \to \tilde{\mathcal U}$.   

\begin{thm}\label{derham} {\rm (i)  ($k=\C$)} There is a natural commutative diagram
$$\begin{CD}
\mathcal B_n @>{\tilde{\psi}}>> \tilde{\mathcal U} @>{\xi}>> \tilde{\mathcal P}_n\\
@AAA @AAA @AAA\\
T_{\Z}L^n @>{\psi}>> {\rm Lie}~\tilde{\mathcal U} @>{\xi}>> {\rm Lie}~\tilde{\mathcal
P}_n\\
\end{CD}
$$
whose exterior square is \cite[10.1.3.1]{h} for $L^n$; further,
${\psi}$ induces an isomorphism
$\psi_{\C}:T_{\C}L^n \zim {\rm Lie}~\tilde{\mathcal U}$.

{\rm (ii)} $T_{DR} L^n \zim  {\rm Lie}~\tilde{\mathcal U}$.

{\rm (iii)} There is a natural map ${\rm Lie}~\tilde{\mathcal U} \to
H^n_{DR}(U_{\bu})$.\end{thm}

\begin{proof}(i) ($k = \C$)   
From (\ref{codered}), (\ref{kapadu}) and (\ref{hosadu}), we obtain an
exact sequence of isogeny one-motives over $\C$
$$0 \to  \frac{[K^0 \to \mathcal G_n]}{[\mathcal V_{n-2} \to
Pic^0_{n-2}]}\otimes\Q \to [\mathcal C_n
\xrightarrow{\rho} \mathcal P_n]\otimes\Q  \to [\Ker(\delta^*_{n-2}: N_{n-2} \to N_{n-1}) \to
0]\otimes\Q \to 0.$$ The proof (in Step 3.) of Theorem \ref{main} also shows
that (\ref{devaki}) is an isomorphism (modulo finite groups) of mixed
Hodge structures. Combining this with (\ref{k3}), we obtain
that 
$h^2_X\otimes\C \zim T_{\C}[\mathcal C_n \xrightarrow{\rho} \mathcal
P_n]$ sits in an exact sequence
$$ 0 \to {\rm Lie}~\mathcal P^{\diamondsuit}_n \to h^2_X\otimes\C \to
\Ker(\delta^*_{n-2}: N_{n-2} \to
N_{n-1})\otimes\C \to 0.$$  
 Recall that there is a map $h^2_X\otimes\C 
 \to \hi^2(X_{\ge n-2}, \mathcal K) = {\rm Lie}~R^2g_*\mathcal K^*$ 
 and that the latter sits in an exact sequence
$$ 0 \to \frac{\hi^1(X_{\ge n-1}, \Gamma)}{H^0(X_{n-2},
\Omega^1(log~Y))} \to \hi^2(X_{\ge n-2}, \mathcal K) \to H^1(X_{n-2},
\Omega^1(log~Y)) \to \hi^2(X_{\ge n-1}, \Gamma).$$ 
Comparing with (\ref{eq:9}), the map $h^2_X\otimes\C \to {\rm
Lie}~R^2g_*\mathcal K^*$ gives a map 
$$\psi_{\C}:T_{\C}L^n\underset{(\ref{nyaya}), (\ref{main})}{\xleftarrow{\sim}} \frac{h^2_X\otimes\C}{NS(X_{n-3})\otimes\C} \to {\rm Lie}~\tilde{\mathcal
U},$$whose composition  with $T_{\Z}L^n \hra T_{\C}L^n$ gives us the
map $\psi$. Since the composite map $ K \oplus
NS_{n-2} \xrightarrow{\psi'} \mathcal U' \xrightarrow{\xi} \mathcal
P'_n$ is $\rho$ --- see (\ref{greece}), it follows easily that the composite map 
$\xi\circ \tilde{\psi}: \mathcal B_n \to \tilde{\mathcal P}_n$ is
$\phi_n$. Thus, the exterior square is the one that
intervenes in the definition \cite[10.1.3.1]{h} of $T_{\Z}L^n$.

We leave it to the reader to check that $\psi_{\C}$ is an isomorphism and
that it is compatible with  the weight
and Hodge filtrations. 

(ii) This follows from (i) by \cite[10.1.9]{h}. 

(iii) Compose the natural maps $${\rm Lie}~\tilde{\mathcal
U} \to \hi^2(X_{\ge n-2}, \mathcal K) \to \hi^n(X_{\bu},
\Omega(log~Y)) \zim H^n_{DR}(U_{\bu}),$$the last being a consequence
of \cite[3.1.8]{h2}, \cite[1.0.3.7]{katzdiff}, \cite[6.11.4]{jamm}.  
\end{proof}

\section{Etale realization of $L^n$}\label{oakland} 
Fix a simplicial pair $(X_{\bu}, Y_{\bu})$ over a field $k$ of
characteristic zero. When $k = \C$,
we have a map  (see below for notations) $$TL^n\otimes\Q
\zim  (T_{\Z}L^n)\otimes_{\Z}\A
\underset{(\ref{main})}{\zim}T_{\Z}I^n(U_{\bu})\otimes_{\Z}\A \hra
H^n(U_{\bu}, \Z(1))\otimes_{\Z} \A
\underset{(\ref{artgro})}{\zim}H^n_{et}(U_{\bu}, \A(1));$$ the first
isomorphism is from \cite[10.1.6]{h}. We shall show that this map can
be constructed purely algebraically for all $k$; combined with
(\ref{main}), (\ref{derham}), this will prove (\ref{dc1}) up to isogeny
for $U_{\bu}$.   By
(\ref{zeusp}), we may and do assume that $k$ is a finitely generated
extension of $\Q$.  

We adopt the following notations: $r$ is a positive integer,
$\ell$ is a positive prime integer,  $H^i_{et}(\bar{V},\hat{\Z}(1))$
is ${\varprojlim}_r
H^i_{et}(\bar{V},\1_r)$ and  $H^i_{et}{}(\bar{V},\mathbb A(1))$ is 
$H^i_{et}{}(\bar{V},\hat{\Z}(1))\otimes\Q$. 
For any one-motive $M$, recall the finite $\Z/{r\Z}$-module
$T_{\Z/{r\Z}}M:= H^0(M\stackrel{L}{\otimes}{\Z/{r\Z}})$
\cite[10.1.5]{h}; note $T_{\ell}M:= {\varprojlim}_d
T_{\Z/{\ell^d\Z}}M$, $M\otimes\hat{\Z} = TM =  {\varprojlim}_r
T_{\Z/{r\Z}}M$, and $M\otimes\A:= TM\otimes\Q$. 
For any commutative group (scheme) $\mathcal A$, $_r\mathcal A$ is the 
$\Ker(\mathcal A \xrightarrow{r}\mathcal A)(\bar{S})$ and 
$\mathcal A_{tor}$ the torsion subgroup (scheme). Let $\Z/{r\Z}$
denote the constant group. All maps in this
section are $\G$-equivariant.

Recall the Kummer sequence of \'etale sheaves on a (simplicial) scheme
$V$ \begin{equation}\label{kummer}
0 \rightarrow \1_{r} \rightarrow \gm \xrightarrow{r}
\gm \rightarrow 0. 
\end{equation}
For proper $V$, this gives, by the divisibility of
$H^0_{et}{}(\bar{V},\gm)$,    
\begin{equation}\label{kap}
H^1_{et}{}(\bar{V},\1_r) \zim
~_r H^1_{et}{}(\bar{V},\gm);\end{equation} for $V$ smooth and proper, we also
obtain  
\begin{equation}\label{kumns}
\kappa_V: NS(\bar{V})\otimes{\Z/{r\Z}} \hra
H^2_{et}{}(\bar{V},\1_r).\end{equation}  
Similarly, using the divisibility of $\mathcal G_n(\bar{S})$, one
obtains an injection
\begin{equation}\label{gandhi}
\pi_0(\mathcal G'_n)(\bar{S})\otimes{\Z/{r\Z}} \hra  H^2_{et}{}(\bar{X}_{\ge
n-1},\1_r).\end{equation} 

\begin{rem}\label{specet}
The spectral sequence \cite[8.1.19.1]{h} has an  \'etale
analogue $E_1^{et}$ \cite[Introduction]{dj} \cite[\S 6]{h1}  \cite[\S
14]{hp}; it calculates the \'etale cohomology of $\bar{U}_{\bu}$ with 
$\Q_{\ell}$-coefficients.

The degeneration of the \'etale analogue at $E_2$ (i.e., that $E_2^{et} =
E_{\infty}^{et}$) and the definition of the weight filtration $W$ 
on $H^*_{et}(\bar{U}, \Q_{\ell}(1))$  is a 
consequence of \cite{weil}; see \cite[\S 6,7]{h1}, \cite[\S 14]{hp}
for these and the fact that the weight filtrations in mixed
Hodge theory are compatible with those in \'etale cohomology, compared
via (\ref{artgro}). 

Recall the Artin-Grothendieck
comparison theorem \cite[\S 4]{art} \cite[p.23]{ka}: for any variety $V$ over
$\C$, one has canonical isomorphisms \begin{equation}\label{artgro}
H^i_{et}(V, \Z/{r\Z}) \zim H^i(V(\C), \Z/{r\Z})\end{equation}
 between the \'etale and classical (singular)
cohomology of $V$.
Any imbedding $\iota: k \hra \C$ gives an isomorphism of $\bar{k}$ with the
algebraic closure $\overline{\iota(k)}$ of $\iota(k)$ in $\C$; this, in
turn,  provides an isomorphism
$H^*_{et}(U_{m}\times_k \bar{k}, \Z/{r\Z}) \zim
H^*_{et}(U_{m}\times_k \overline{\iota(k)}, \Z/{r\Z}) \zim H^*_{et}(U_{m}\times_{\iota} \C,\Z/{r\Z})$; the last
isomorphism is from \cite[VI 2.6]{mi}.     
Now  (\ref{artgro}) 
 and the degeneration \cite[8.1.20 (ii)]{h} of \cite[8.1.19.1]{h} at $E_2$
together provide another proof that   $E_2^{et} =
E_{\infty}^{et}$.\qed \end{rem}  

\begin{lem}\label{bombay} One has {\rm{(i)}}   
$$\frac{H^1_{et}{}(\bar{X}_{\ge n-1},\A(1))}
{H^1_{et}{}(\bar{X}_{n-2}, \A(1))} \zim \sigma_{\ge n-1}
H^n_{et}{}(\bar{X}_{\bu},\A(1)) =
W_{-1}H^n_{et}{}(\bar{X}_{\bu},\A(1)) \zim  W_{-1}
H^n_{et}{}(\bar{U}_{\bu},\A(1)).$$  
{\rm{(ii)}} $H^1_{et}{}(\bar{X}_{\ge
n-1},\hat{\Z}(1))\xrightarrow{\sim} T[0 \to {\mathcal
G}_n]$ and  $T[0 \to {\mathcal P}_n]\otimes\Q \xleftarrow{\sim} W_{-1}
H^n_{et}{}(\bar{U}_{\bu},\A(1))$. 
\end{lem}

\begin{proof} (i) As in (\ref{malli}), this  follows 
from the degeneration at $E_2$ of the \'etale analogue of
\cite[8.1.19.1]{h} for $\Q_{\ell}(1)$ --- see (\ref{specet}).  

(ii) The first isomorphism follows from (\ref{kap}) for $V = X_{\ge
n-1}$ by taking the
inverse limit over $r$ (note $\pi_0(\mathcal G')_{tor}$ is
finite). Similarly, one has (\ref{picmd} (i))
 $H^1_{et}{}(\bar{X}_{n-2},\Z_{\ell}(1))
\xrightarrow{\sim}T_{\ell}Pic(X_{n-2})$.
 Combining this with (i) gives the second isomorphism. 
\end{proof}

\begin{prop}\label{finel} {\rm{(i)}} $H^0([W_{Y_m} \xrightarrow{\lambda_{m}}
Pic_{X_{m}}]\stackrel{L}{\otimes} \Z/{r\Z}))
\xleftarrow{\sim} H^1_{et}{}(\bar{U}_{m}, \1_{r})$.

{\rm{(ii)}} $H^0([K \xrightarrow{\vartheta'}
  \mathcal G'_n]\stackrel{L}{\otimes
} \Z/{r\Z}) \xleftarrow{\sim} H^1_{et}{}(\bar{U}_{\ge n-1},\1_{r})$.
\end{prop}

\begin{proof} It is exactly identical to \cite[10.3.6]{h}; we repeat the proof of (i) here for the
convenience of the reader. 

By (\ref{haccha1}), $H^1_{et}{}(\bar{U}_{m}, \1_{r}) \zim
\hi^1_{et}(\bar{U}_{m}, [\gm \xrightarrow{r} \gm])$ is 
the set of isomorphism classes of pairs $(\mathcal L, \alpha)$
where $\mathcal L$ is an invertible sheaf 
on $\bar{U}_m$ together with an isomorphism $\alpha: \mathcal
L^{\otimes r} \zim
\mathcal O$. Let $(\mathcal L, \alpha)$ be such a pair. The invertible
sheaf $\mathcal L$ extends to an invertible sheaf $\tilde{\mathcal L}$ on
$\bar{X}_m$, and there exists a divisor $E$ of $\bar{X}_m$ with support in $\bar{Y}_m$ such that
$\alpha$ extends to an isomorphism $\alpha': \tilde{\mathcal
L}^{\otimes r} \zim \mathcal
O(E)$. If there exists an isomorphism $\beta$ of $\tilde{\mathcal
L}^{\otimes r}$ with $\mathcal O(E)$, this isomorphism is uniquely
determined up to multiplication by an element of the divisible group
$H^0_{et}(\bar{X}_m, \gm)$. One deduces that the pair $(\tilde{\mathcal L},
E)$ determines $(\mathcal L, \alpha)$ up to isomorphism. For a pair
$(\tilde{\mathcal L}, E)$  to come from a suitable  $(\mathcal L,
\alpha)$, it is necessary and sufficient that $r[\tilde{\mathcal L}] =
[\mathcal O(E)]$ in $Pic_{\bar{X}_m}$.  It comes from $(\mathcal
O_{\bar{X}_m}, 0)$ if and only if it is of the form $(\mathcal O_{\bar{X}_m}(E), rE)$. This
identifies $H^1_{et}{}(\bar{U}_{m}, \1_{r})$ with $H^0$ of the complex
which is the tensor product of $[W_{Y_m}(\bar{S}) \xrightarrow{\lambda_{m}}
Pic(\bar{X}_{m})]$ (degrees $0$ and $1$) and $ [\Z \xrightarrow{r}
\Z]$ (degrees $-1$ and  $0$):
$$
\begin{CD}
W_{Y_m}(\bar{S}) @>{\lambda_{m}}>>   Pic(\bar{X}_{m})\\
@AA{r}A   @AA{-r}A\\
W_{Y_m}(\bar{S}) @>{\lambda_{m}}>>  Pic(\bar{X}_{m})\\
\end{CD}
$$
This proves the first isomorphism.\end{proof}

\begin{lem}\label{adivasi} 
{\rm{(i)}} The maps $[\mathcal B_n
\xrightarrow{\phi_{n}}\tilde{{\mathcal P}}_n] \xleftarrow{b} 
 [\mathcal B'_n \xrightarrow{\rho} {\mathcal P}_n] \xrightarrow{a}  [\frac{\mathcal
C'_n}{NS_{n-3}}
\xrightarrow{\rho} \mathcal P'_n]$ of complexes {\rm (}concentrated in
 degrees $0$ and $1${\rm)} 
induce isomorphisms {\rm(}modulo
finite groups{\rm)} 
\begin{equation}\label{maharaja}
TL^n \xleftarrow{\sim}~{\varprojlim}_rH^0([\mathcal B'_n \xrightarrow{\rho} {\mathcal P}_n]\stackrel{L}
{\otimes} {\Z/{r\Z}}) \zim~{\varprojlim}_rH^0([\frac{\mathcal
C'_n}{NS_{n-3}}
\xrightarrow{\rho} \mathcal P'_n]\stackrel{L} {\otimes} {\Z/{r\Z}}).
\end{equation}

{\rm{(ii)}} The map $[\mathcal C_n \xrightarrow{\rho} \mathcal P_n]$ 
to $[\mathcal C'_n \xrightarrow{\rho} \mathcal P'_n]$ induces an
isomorphism  {\rm(}modulo
finite groups{\rm)}  
\begin{equation}\label{maharani}
{\varprojlim}_rH^0([\mathcal C_n \xrightarrow{\rho} {\mathcal P}_n]\stackrel{L}
{\otimes} {\Z/{r\Z}}) \zim~{\varprojlim}_rH^0([\mathcal C'_n
\xrightarrow{\rho} \mathcal P'_n]\stackrel{L} {\otimes} {\Z/{r\Z}}).
\end{equation}
\end{lem}

\begin{proof} (cf. \cite[10.3.5 ii]{h}) (i) The first isomorphism is
clear because $b$ is an isogeny;  
the cokernel $[\im(\rho') \to \pi_0(\mathcal P'_n)]$ of the map
$a$ is quasi-isomorphic to $[ 0 \to \Coker(\rho')]$ and
$\Coker(\rho')_{tor}$ is finite. This yields the second
isomorphism. 

(ii) The argument is similar to that of (i).\end{proof}

The map $[\sigma_{\ge n-1}\gm \xrightarrow{r} \gm] \to [\gm
\xrightarrow{r} \gm]$ of \'etale complexes on $U_{\ge n-2}$  
gives the exact sequence 
\begin{equation}\label{hippo} 
\frac{\hi^2_{et}{}(\bar{U}_{\ge n-2}, [\sigma_{\ge n-1} \gm \xrightarrow{r}
\gm])}{H^1_{et}(\bar{U}_{n-2}, \gm)} \to \hi^2_{et}{}(\bar{U}_{\ge n-2}, [\gm \xrightarrow{r}
\gm]) \to H^2_{et}(\bar{U}_{n-2}, \gm)
\end{equation}

\begin{thm}\label{biget}{\rm{(i)}}  One has a canonical isomorphism
$$\Pi_r: \frac{\hi^2_{et}{}(\bar{U}_{\ge n-2},
[\sigma_{\ge n-1} \gm \xrightarrow{r} 
\gm])}{H^1_{et}(\bar{U}_{n-2}, \gm)} \zim H^0([\mathcal C'_n
\xrightarrow{\rho} \mathcal 
P'_n]\stackrel{L} {\otimes} {\Z}/{r\Z}).$$This gives an inclusion
$\Lambda': T[\mathcal C_n\xrightarrow{\rho} \mathcal P_n]\hra H^2_{et}(\bar{U}_{\ge
n-2}, \hat{\Z}(1))$. 

{\rm{(ii)}} The map $\Lambda'$ of {\rm{(i)}} induces an inclusion
$\Lambda_{et}: TL^n\otimes\Q \hra  H^n_{et}{}(\bar{U}_{\bu},\A(1)).$ 
  
{\rm{(iii)}} One has $\Lambda_{et}(T(W_{j}L^n \otimes\Q)) \subseteq W_j
H^n_{et}{}(\bar{U}_{\bu},\A (1))$ with  equality for $j = -2, -1$.
\end{thm}

\begin{proof} (i) The proof follows that of
(\ref{finel}).  Via (\ref{haccha1}), $\hi^2_{et}{}(\bar{U}_{\ge n-2},
[\sigma_{\ge n-1} \gm \xrightarrow{r} 
\gm])$ is the group of isomorphism classes of triples $(I, \mathcal L,
\alpha)$ where $I$ is a $\gm$-torsor on $\bar{U}_{n-2}$ and $\mathcal
L$ is a $\gm$-torsor on $\bar{U}_{\ge n-1}$ and $\alpha: \delta^*I
\zim \mathcal L^{\otimes r}$. The class $[L] \in H^1_{et}(\bar{U}_{n-2},
\gm)$ of a $\gm$-torsor $L$ gets mapped to the class $\delta^*[L]$ 
of the triple $(L^{\otimes r}, \delta^*L, \alpha_L)$ where $\alpha_L$ 
is the canonical isomorphism $\delta^*(L^{\otimes r}) \zim
(\delta^*L)^{\otimes r}$. 

Fix the class $f$ of a triple $(I, \mathcal L,
\alpha)$. The torsors ${I}$ and $\mathcal L$ extend to torsors $\tilde{I}$ and
$\tilde{\mathcal L}$ on $\bar{X}_{n-2}$ and $\bar{X}_{\ge n-1}$
respectively and  there exists a divisor $E$ of $\bar{X}_{\ge n-1}$
supported on $\bar{Y}_{\ge n-1}$ such that $\alpha$ extends to an
isomorphism $\tilde{\alpha}: \delta^*\tilde{I}\otimes \mathcal O(E)
\zim \tilde{\mathcal L}^{\otimes r}$. This says that the elements $w=
(E, [\tilde{I}]) \in K \oplus NS_{n-2}$ and $[\tilde{\mathcal L}] \in
\mathcal P'_n$ satisfy $\rho(w) = r [\tilde{\mathcal L}]$; in other
words, the pair $(w, [\tilde{\mathcal L}])$ 
 defines an element $\Pi_r(f)$ of  $H^0([\mathcal C'_n
\xrightarrow{\rho} \mathcal 
P'_n]\stackrel{L} {\otimes} {\Z}/{r\Z}).$ 

This element $\Pi_r(f)$ is clearly zero 
if our triple is isomorphic to $\delta^*[L]$: if $L$ extends to a
$\gm$-torsor $\tilde{L}$ on $\bar{X}_{n-2}$, then $\Pi_r(f)$ is given
by $w = (0, r[\tilde{L}])$ and the class of $\rho([\tilde{L}]) \in
\mathcal P'_n$.  If $\tilde{I}'$ is another extension of $I$, then there is an
 isomorphism  $\tilde{I} \zim \tilde{I'} \otimes\mathcal
O(d)$ ($d \in W_{n-2}$). Using $\tilde{I}'$ instead of
 $\tilde{I}$ modifies $w$ by $\gamma_{n-3}(d)$. If $\tilde{\mathcal
L'}$ is another extension of $\mathcal L$, then there is an
 isomorphism $\tilde{\mathcal L} \zim  \tilde{\mathcal
L'}\otimes \mathcal O(F)$ ($F \in K$). Using $ \tilde{\mathcal
L'}$ instead of $\tilde{\mathcal L}$ modifies $w$ by $\pm rF$.  
This shows that 
 the map $\Pi_r$ is well-defined; on the image of the map
$\frac{H^1_{et}(\bar{U}_{\ge n-1},
\1_r)}{H^1_{et}(\bar{U}_{n-2}, \1_r)} \to \hi^2_{et}{}(\bar{U}_{\ge n-2},
[\sigma_{\ge n-1} \gm \xrightarrow{r} 
\gm]),$ the map $\Pi_r$ reduces to the maps (isomorphisms) of
(\ref{finel}). 

Given an element $g$ of $H^0([\mathcal C'_n
\xrightarrow{\rho} \mathcal 
P'_n]\stackrel{L} {\otimes} {\Z}/{r\Z})$, pick representatives 
 $u\in K$, $\mathcal J_1 \in Pic_{n-2}$, and $\mathcal J \in \mathcal G'_n$ for
$g$. By definition, $\mathcal O(u)\otimes\delta^* \mathcal J_1 \zim
\mathcal J^{\otimes r}\otimes \mathcal J_2$ where $\mathcal J_2$ is an element
of $Pic^0_{n-2}$; by the divisibility of $Pic^0_{n-2}(\bar{S})$, the
relation can be rewritten as $\mathcal O(u)\otimes\delta^* \mathcal J_1 \zim
\mathcal J^{\otimes r}$ for an appropriate element of $\mathcal
G'_n$. Restricting to $\bar{U}_{\ge n-1}$, we obtain an isomorphism
$\beta: \delta^* \mathcal J_1 \zim
\mathcal J^{\otimes r}$; the triple $(\mathcal J_1, \mathcal J,
\beta)$ defines an element of $\hi^2_{et}{}(\bar{U}_{\ge n-2},
[\sigma_{\ge n-1} \gm \xrightarrow{r} 
\gm])$.  One easily verifies  that this is an inverse to the map
$\Pi_r$ and that $\Pi_r$ is an isomorphism. The map $\Lambda'$ is
obtained by composing the inverse of ${\varprojlim}_r\Pi_r$  with the
maps in (\ref{hippo}) and  (\ref{maharani}).   

\nd (ii) Since $\delta^*H^2_{et}(\bar{U}_{n-3}, \A(1))$ is contained in $\Ker(H^2_{et}(\bar{U}_{\ge
n-2}, \A(1)) \to  H^n_{et}(\bar{U}_{\bu},
\A(1)))$, the  
following 
diagram gives us $\Lambda_{et}$:
$$\begin{CD}
H^2_{et}(\bar{U}_{n-3}, \A(1)) @>{\delta^*}>> H^2_{et}(\bar{U}_{\ge
n-2}, \A(1))\!  \!  \! \!  \!  \!
@>>> H^n_{et}(\bar{U}_{\bu}, \A(1))\\
@A{(\ref{kumns})}AA @AA{{\rm(i)}}A @.\\
NS(\bar{X}_{n-3})\otimes \A  @>{\delta^*}>> T[\mathcal C_n
\xrightarrow{\rho} {\mathcal P}_n]\otimes\Q;\\ 
\end{CD}
$$note $NS(\bar{X}_{n-3}) \otimes\A \zim T[
NS(\bar{X}_{n-3}) \to 0]\otimes\Q$. 
Since $W_{-1}H^2_{et}(\bar{U}_{n-3}, \A(1)) =0$ \cite[7.2]{hp}
\cite[3.2.15]{h2}, the
intersection $$\delta^*H^2_{et}(\bar{U}_{n-3}, \A(1)) ~\cap~ W_{-1} H^2_{et}(\bar{U}_{\ge
n-2}, \A(1))$$is zero. Now (\ref{bombay} (ii)) shows that $\Lambda_{et}$ 
is injective on $W_{-1}TL^n\otimes\Q$. It remains to show  the
injectivity of either
$$\Lambda_{et}: \mathcal B_n(S)\otimes\A \to 
Gr^W_{0}H^n_{et}(\bar{U}_{\bu}, \A(1))$$ or its $\Q_{\ell}$-analogue. 
The base change of this map via any
$\iota: k \hra \C$ can be identified, via (\ref{artgro}) and the
compatibility of the \'etale and classical cycle class maps \cite[5.3]{jamm}, with the
\emph{injective} map
$\nu$ of (\ref{kapadu}(i)):$$\mathcal
B_n(\C)\otimes\Q_{\ell} \underset{\nu \otimes{id}}{\xrightarrow{\sim}}
Gr^W_0t^n(U_{\bu, \iota})\otimes\Q_{\ell} \hra Gr^W_0H^n(U_{\bu,
\iota}, \Q(1))\otimes\Q_{\ell}
\underset{(\ref{artgro})}{\xleftarrow{\sim}} Gr^W_0H^n_{et}(U_{\bu,
\iota},  \Q_{\ell}(1)).$$

\nd (iii) Follows from the definition of the filtration $W$
(\cite{h1}, \cite{weil}, \cite[\S 14]{hp}) and
(\ref{bombay}). \end{proof}

\begin{defn} $t^n_{\ell}(U_{\bu}):= \Lambda_{et}(T_{\ell}L^n
\otimes\Q)$ is a $\G$-invariant
$\Q_{\ell}$-subspace of $H^n_{et}{}(\bar{U}_{\bu};
\Q_{\ell}(1))$. Similarly, $t^n_{\A}(U_{\bu}):=
\Lambda_{et}(TL^n\otimes\Q)$ is a $\G$-invariant $\A$-submodule of
$H^n_{et}{}(\bar{U}_{\bu}; \A(1))$.    
\end{defn}   

\begin{rem}\label{smalet} Over $k = \C$, one has a commutative diagram  
$$\begin{CD}
 TL^n\otimes\Q @>{\Lambda_{et}}>> H^n_{et}({U}_{\bu},\A(1))\\    
@V{\wr}VV @V{\wr}V(\ref{artgro})V\\
T_{\Z}L^n\otimes\A @>{\Lambda}>{(\ref{main})}> H^n(U_{\bu}(\C),
\A(1));\\ 
\end{CD}
$$
the left vertical map is from \cite[10.1.6]{h}. This diagram, together
with (\ref{main}), also provides the injectivity of $\Lambda_{et}$
in (\ref{biget}(ii)). It also provides canonical isomorphisms
$t^n_{\ell}(U_{\bu}) \zim 
t^n(U_{\bu})\otimes\Q_{\ell}$ and $t^n_{\A}(U_{\bu}) \zim
t^n(U_{\bu})\otimes\A$.\qed\end{rem}

\section{The one-motives $L^n(V)$}\label{varieties}
In this section, $k$ is any field of characteristic zero. Throughout
the rest of the paper, $V$ will denote a variety (over $S$). 

We shall now
translate our results for simplicial schemes into ones for varieties;
in particular, we show how to 
construct  one-motives $L^n(V)$ associated with any
algebraic variety $V$ over $S$. We
will also show that $L^n(-)$ is a contravariant
functor. 

The main result of this section is the proof of Conjecture
\ref{dc1} for fields of characteristic zero.\ssp  

\nd {\bf Proper hypercoverings.}\ssp 

Let $V$ be a variety over $S$. 
By \cite[6.2.8]{h}, \cite[5.3]{sga4},  there exists a
simplicial pair 
($X_{\bu},Y_{\bu}$)  
with a proper map $\alpha: U_{\bu}\rightarrow V$ which makes $U_{\bu}$
into a proper hypercovering of 
$V$. Namely, we have\ssp 

\nd (i) for any $\iota: k \hookrightarrow \C$,  an isomorphism of mixed
Hodge structures \cite[8.2.2]{h}
\begin{equation}\label{sonatree}
\alpha^*: H^*(V_{\iota},\Z) \zim H^*(U_{\bu, \iota},\Z);
\end{equation}(ii) an isomorphism of $\G$-modules     
\begin{equation}\label{sonaty}
\alpha^*:
H^*_{et}(\bar{V},\hat{\Z}(1)) \zim H^*_{et}(\bar{U}_{\bu},\hat{\Z}(1)); 
\end{equation}(iii)  and an isomorphism of $k$-vector spaces
\begin{equation}\label{sonati}
\alpha^*: H^*_{DR}(V) \xrightarrow{\sim} H^*_{DR}(U_{\bu}). 
\end{equation} 

For (ii), one uses the result \cite[4.3.2]{sga4} that a proper
surjective morphism is a map of universal cohomological descent for
\'etale torsion sheaves and proceeds as in \cite[6.2.5]{h};
cf. \cite[p.306]{be}. 

In (iii), the definition of $H^*_{DR}(V)$ (as indicated in
\cite[p. 89]{pgal}, \cite[6.11.4]{jamm})
is that of \cite{hart} or, by \cite[1.2]{illucie}, as the \emph{crystalline} cohomology
of $V$; see also \cite[Expose
III]{guillen}, especially 1.3, 1.13, 1.14. for a cubical variant.

We denote by $\beta_U$ the inverse of any of these isomorphisms
$\alpha^*$. 
\begin{rem}\label{bertie} Given two such proper hypercoverings of $V$, one can find another
such proper hypercovering which dominates both. More generally, given
any morphism $h:V \rightarrow W$ between two schemes, one can find such
proper hypercoverings of $V$ and $W$ and a morphism between them which
lifts the morphism $h$; we refer to \cite[5.1.4]{sga4},
\cite[6.2.8]{h} for details.    \qed \end{rem} 
Suppose given two proper hypercoverings of $V$ fitting into a diagram  
$$
\begin{CD}
V @<<< U_{\bu} @>>> X_{\bu} @<<< Y_{\bu}\\
@| @V{\theta}VV @V{\theta}VV @.\\
V @<<< 'U_{\bu} @>>> 'X_{\bu} @<<< 'Y_{\bu}\\
\end{CD}
$$
which yields a morphism of one-motives
$$\theta^*_L: L^n('X_{\bu}, 'Y_{\bu}) \rightarrow
L^n(X_{\bu}, Y_{\bu}).$$

\begin{lem}\label{pehla} The morphism $\theta^*_L$ is an isogeny. \end{lem}

\begin{proof} 
The construction of the one-motives $L^n(U_{\bu})$ and $L^n('U_{\bu})$
relies only on a 
finite number 
of schemes and a finite number of associated morphisms between
them. Therefore, we may assume that these one-motives are defined
over a finitely generated subfield $k'$ of $k$.  Such a field $k'$ always
admits an embedding into $\C$. For any such embedding $\iota: k' 
\hookrightarrow \C$, one has a commutative diagram
$$\begin{CD}
T_{\Z}L^n('X_{\bu}, 'Y_{\bu})\otimes\Q @>>> H^n('U_{\bu}\times\C, \Q(1))\\
@V{T\theta^*_L}VV @V{\theta^*}V{\wr}V\\
T_{\Z}L^n(X_{\bu}, Y_{\bu})\otimes\Q @>>> H^n(U_{\bu}\times\C, \Q(1));\\
\end{CD}
$$the isomorphism $\theta^*$ \cite[8.2.2]{h} of mixed Hodge
structures provides an isomorphism $t^n('U_{\bu}) \zim t^n(U_{\bu})$.  
By Theorem \ref{main} applied to $'U_{\bu}\times\C$ and
$U_{\bu}\times\C$,
the map $T\theta^*_L$ is an isomorphism of $\Q$-mixed Hodge
structures. By \cite[10.1.3]{h},  $\theta^*_L$ is an
isogeny.\end{proof}

\begin{defn}
The isogeny one-motive $L^n(V)\otimes\Q$  is the
isogeny class of the one-motive  $L^n(X_{\bu},Y_{\bu})$ of any
simplicial pair $(X_{\bu}, Y_{\bu})$ corresponding to a proper hypercovering  
$U_{\bu}$ of $V$.\end{defn}

\begin{defn}\label{banff}
Let $t^n_{\ell}(V)$ be the image of 
$t^{n}_{\ell}(U_{\bu})$ under the map $\beta_U$, inverse to the
isomorphism   
$\alpha^*: H^n_{et}(\bar{V},\Q_{\ell}(1)) \xrightarrow{\sim}
H^n_{et}(\bar{U}_{\bu},\Q_{\ell}(1))$ for any proper hypercovering $\alpha:
U_{\bu} \rightarrow V$  corresponding to a simplicial pair $(X_{\bu},
Y_{\bu})$. Also define $t^n_{\A}(V) =
\beta_U(t^n_{\A}(U_{\bu}))$, a $\A$-submodule of $H^n_{et}(\bar{V},
\A(1))$. Note that   $t^n_{\A}(V)$ and
$t^n_{\ell}(V)$ are independent
of the choice of the proper hypercovering $U_{\bu}$ of $V$. 
\end{defn}

By (\ref{bertie}), any two proper hypercoverings of $V$ are dominated
by a third one; thus, as in \cite[\S8.2]{h}, the isogeny one-motive
$L^n(V)\otimes\Q$ is well-defined by \cite[5.1.4]{sga4},
\cite[6.2.8]{h}. 

\begin{lem}\label{funter}  
$L^n(-)\otimes\Q$, $t^n_{\ell}(-)$ and
$t^n_{\A}(-)$ are contravariant functorial.\end{lem}  
\begin{proof} Given a morphism $f:V \rightarrow W$, there exists a proper hypercovering
$(X_{\bu},Y_{\bu})$ (resp. $(Z_{\bu},I_{\bu})$) of $V$ (resp. $W$) and
a morphism $\theta$ between them which lifts $f$ \cite[6.2.8]{h}. This induces a  
morphism $\theta^*:L^n(W)\otimes\Q \rightarrow
L^n(V)\otimes\Q$. One has to check that this induced morphism
does not depend upon the choices of the 
auxiliary hypercoverings and for this we refer to
(loc. cit).\end{proof}

The next lemma, inspired by \cite{milram}, will be used to deduce an
integral version (i.e., an actual one-motive $L^n(V)$) of the isogeny
one-motive 
$L^n(V)\otimes\Q$. Let $\mathcal C$ be  
the category whose objects are
 pairs $(R,L)$ where $R$ is an isogeny one-motive over
$k$ and $L$ is a $\G$-stable $\hat{\Z}$-lattice of the $\A$-module 
$TR$, and morphisms from $(R_1, L_1)$ to $(R_2, L_2)$ are those
morphisms  $\phi:
R_1 \to R_2$ of isogeny one-motives such that the induced map $\phi_*: TR_1
\to TR_2$ satisfies $\phi_*(L_1) \subset L_2$.       
\begin{lem}\label{intero} 
The natural functor $M \mapsto (M\otimes\Q, M\otimes\hat{\Z})$ 
from the category $\mathcal M_k$ of one-motives
over $k$ to $\mathcal C$ is an
equivalence of categories. 

\end{lem} 

\begin{proof} This follows from \cite[2.2]{milram}; the main point is
as follows: as $\hat{\Z}$-lattices  of the $\A$-module $TR$ are compact and open
 in $TR$, any two $\G$-stable $\hat{\Z}$-lattices in $TR$  are commensurable.
\end{proof}

 Let $V$ be any variety over $k$. Consider the $\G$-stable
$\hat{\Z}$-lattice  $$ t^n_{\A}(V) \cap 
H^n_{et}(\bar{V}, \hat{\Z}(1))/{\rm torsion} \subset H^n_{et}(\bar{V},
\A(1))$$ in $t^n_{\A}(V)$. Since the map $\Lambda_{et}: 
TL^n(V)\otimes\Q \to H^n_{et}(\bar{V},
\A(1))$ is injective with image $ t^n_{\A}(V)$, 
this defines a $\G$-stable $\hat{\Z}$-lattice of
$TL^n(V)\otimes\Q$ which, by (\ref{intero}), determines a one-motive
$L^n(V) =:L^n(V/k)$ over $k$.

\begin{thm}\label{peddha} Conjecture {\rm \ref{dc1}} is
true for fields of characteristic zero. The one-motives $L^n(V/k)$
defined for an arbitrary variety $V$ over a field $k$ of
characteristic zero are functorial in $V$ and $k$. Furthermore, one has 

{\rm (i)} an injection of $\G$-modules $\Lambda_{et}: TL^n(V/k) \hra
H^n_{et}(\bar{V}, \hat{\Z}(1))/{\rm torsion}$; 

{\rm (ii)} an injection of $k$-vector spaces $T_{DR}L^n(V/k) \hra
H^n_{DR}(V)$; 

{\rm (iii)} a canonical and functorial isomorphism $L^n(V/{\C}) \zim
I^n(V)$. 

If $k$ is finitely generated over $\Q$, then $\Lambda_{et}$ is
compatible with the weight   
filtration: $\Lambda_{et}(W_iL^n(V/k)\otimes\A) \subset W_iH^n_{et}(\bar{V},
\A(1))$ for $i=-2,-1,0$ {\rm(}equality for $i=-2,-1${\rm)}.    
\end{thm}
\begin{proof} The contravariant functoriality of $L^n(V)$ is proved as in
(\ref{funter}) and compatibility with base change  is clear
(\ref{zeusp}).  

 (i) follows from the definition of $L^n(V/k)$,
(\ref{pehla}), (\ref{biget}), (\ref{sonaty}); and (ii) follows from
(\ref{derham}), (\ref{sonati}).  
(iii) follows from the definition of $L^n(V/k)$, (\ref{main}),
(\ref{sonatree}), and (\ref{artgro}).

The last statement follows from (\ref{biget}).
\end{proof}

\begin{rem} Let $N$ be the dimension of $V$. The one-motives $L^n(V)$
vanish for $n > N + 1$; this follows 
from weight considerations \cite[7.3]{hp}. If $V$ is smooth, then
$L^n(V)$ is zero for $n > 2$ (\ref{toto}).  \qed \end{rem} 

\section{Positive characteristic}\label{pos+}
In this section, $k$ is a field of characteristic
$p >0$, $k^{perf}$
its perfection, $k^{sep}$ a separable closure ($\G_{sep}$ is the
associated Galois group), 
$\bar{k}$ (an algebraic closure of $k$) the
compositum of $k^{perf}$ and $k^{sep}$; $\ell$ is  a prime
distinct from $p$.   

Even though we do not assume $k$ to be perfect, all our results
involve passage to $k^{perf}$.  For any variety $V$ over a perfect field $k$ 
and any integer $n$, we use ``neat  hypercoverings'' to  construct
$L^n(V)\otimes\Z[1/p]$ ($0 \le n \le 2$), 
$J^n(V)\otimes\Z[1/p]$ ($ n \ge 0$) --- one-motives up to $p$-isogeny, i.e.,
objects of $\mathcal M_k
\otimes\Z[1/p]$ --- which  are
contravariant functorial; these provide a partial proof of
(\ref{dc1}) for $k$.  
Finally, we reduce Conjecture \ref{dc1} (up to $p$-isogeny), if $k$ is
perfect, to the 
validity of (\ref{tatec}) for surfaces; an analogous result 
(\ref{imper}) holds, even if $k$ is not perfect, under the additional
assumption of ``resolution of singularities''.\medskip  

\nd {\bf Generalized one-motives.}\ssp

 Let  $(X_{\bu}, Y_{\bu})$ be
 a simplicial pair
over $S$ and $f$ be the structure map of $X_{\ge n-1}$. The sheaf 
$R^1f_* \mathcal O^*$ on $S_{fppf}$ is representable (\ref{reps}).  The 
(reduced) neutral component of the  corresponding group scheme 
is not guaranteed to be a semi-abelian variety unless $k$ 
is perfect. 
\begin{thm} {\rm \cite[Corollary 2.3, p.288]{dg}} Let $G$ be a locally algebraic group scheme over a field
$F$. If $F$ is perfect, then the 
reduced scheme $G^{red}$ is a smooth group scheme.\end{thm}

This motivates the following definition. 
\begin{defn}\label{rtp} A generalized one-motive over $k$ 
is a two-term complex  $M =[B \xrightarrow{u} G]$ of group schemes
over $S$ such that, after base change to {\rm
Spec}~$k^{perf}$, 
$[B \xrightarrow{u} G^{red}]$ is a one-motive over {\rm
Spec}~$k^{perf}$.  
\end{defn}

The category $\tilde{\mathcal M}_k$ of generalized one-motives over
$k$ is an additive category. 
The functors $M \mapsto T_{\ell}M$ ($\ell \neq p$), $M\mapsto T^pM$ on
$\mathcal M_k$ extend to the category $\tilde{\mathcal M}_k$.

 Our constructions in Section \ref{constr} provide
generalized one-motives 
$L^n$, $J^n$ over $k$ (and one-motives $L^n$, $J^n$ over $k^{perf}$) 
associated with the simplicial pair $(X_{\bu}, Y_{\bu})$
over $S$; these are contravariant functorial and they are 
compatible with base change.\medskip    

\nd {\bf Relating $L^n$ to \'etale cohomology.}\ssp

For any (simplicial) variety $V$,  set
$H^*_{et}(\bar{V}, \Z^p(1)): =  
{\varprojlim}_r H^*_{et}(\bar{V}, \1_r)$ with $(r,p) = 1$,
$H^*_{et}(\bar{V}, \A^p(1)):= H^*_{et}(\bar{V}, \Z^p(1))\otimes\Q$.  
 
Using the  generalized one-motives
$L^n$, one checks that the proofs of (\ref{finel}), (\ref{adivasi}),
(\ref{biget}(i)) are valid with $\A^p$ instead of $\A$. Also, (\ref{bombay}) is valid: the  \'etale
analogue $E_{*}^{et}$ of \cite[8.1.19.1]{h} 
calculating the \'etale cohomology of $\bar{U}_{\bu}$ with 
$\Q_{\ell}$-coefficients degenerates at $E_2$ (\ref{specet}).
This is a consequence of the weight filtration $W$ \cite[5.3.6]{weil2} 
on  $H^i(\bar{V}, \Q_{\ell})$ for any scheme $V$ over $S$. 
Compatibility with $W$ forces the vanishing of the differentials
$d_r$ of  $E_{*}^{et}$ for $r >1$; cf. \cite[5.3.7]{weil2} for another
application. The point is that given any finite set of (smooth) schemes $X_i$,
$\tilde{Y}_j$, there exist  an integral scheme $\mathcal V$ of finite
type over Spec~$\F_p$, a dominant 
morphism $\eta: S \to \mathcal V$, smooth schemes $\mathcal X_i$, 
 $\tilde{\mathcal
Y}_j$ over Spec~$\F_p$ with smooth proper morphisms $\mathcal X_i \to
\mathcal V$, $\tilde{\mathcal
Y}_j  \to \mathcal V$ such that $\mathcal X_i \times_{\eta} S = X_i$,
$ \tilde{\mathcal Y}_j \times_{\eta}S = \tilde{Y}_j$. The terms of the
spectral sequence $E_{*}^{et}$ are pull-backs via $\eta$ of pure lisse sheaves on
$\mathcal V$ (purity follows from \cite[6.2.6]{weil2});
cf. \cite[p.89, pp.115-117]{jamm}.\medskip 

\nd {\bf Relations with the Tate conjecture.}\ssp

  The proof of the
injectivity of the map 
$\Lambda_{et}: T_{\ell}L^n\otimes\Q \to H^n_{et}(\bar{U}_{\bu},
\Q_{\ell}(1))$  of (\ref{biget}(ii)) uses  
(\ref{kapadu}(i)) which, in turn, is based on (\ref{h11}); the
injectivity is necessary for a proof of 
(\ref{pehla}) in this context. Thus we 
need a proof (valid in positive characteristic) of the injectivity of 
\begin{equation}\label{limbs} 
\nu_{et}:
\mathcal B_n({S})\otimes\Q_{\ell} \to 
Gr^W_{0}H^n_{et}(\bar{U}_{\bu}, \Q_{\ell}(1))\end{equation}where, by the
degeneration (\ref{specet}) of $E_{*}^{et}$ at $E_2$, we have  ---  
(\ref{kapadu}) ---    
$$Gr^W_{0}H^n_{et}(\bar{U}_{\bu},\Q_{\ell}(1)) =
 \frac{\Ker(H^0(\tilde{Y}_{n-1})
\oplus H^2(\bar{X}_{n-2}) 
\xrightarrow{t_{n-1}}  H^0(\tilde{Y}_{n}) \oplus
H^2(\bar{X}_{n-1}))}{\im(H^0(\tilde{Y}_{n-2}) \oplus H^2(\bar{X}_{n-3}) 
\xrightarrow{t_{n-2}}  H^0(\tilde{Y}_{n-1}) \oplus
H^2(\bar{X}_{n-2}))};$$ here $H^0(\tilde{Y}_m) = H^0_{et}(\tilde{Y}_m\times
\bar{S}, \Q_{\ell})$ and $H^2(\bar{X}_m) = H^2_{et}(\bar{X}_m,
\Q_{\ell}(1))$ and the maps $t_m$ are as in (\ref{soya}). By
\cite[2.6, p. 224]{mi},  $H^2_{et}({X}_{m}\times k^{sep}, \Q_{\ell}(1))) \zim
H^2_{et}(\bar{X}_{m}, \Q_{\ell}(1)))$. 
\begin{lem}\label{maam} Let $k$ be a finitely generated extension of the prime
field $\F_p$. The map $\nu_{et}$ of {\rm (\ref{limbs})} 
 is injective if either 

{\rm (i)} $\delta^*_{n-3}H^2_{et}({X}_{n-3}\times
 k^{sep},\Q_{\ell}(1))^{\G_{sep}}$ is equal to  the image of the map
$$\delta^*_{n-3}NS_{n-3}(S)\otimes\Q_{\ell} \hra
NS_{n-2}(S)\otimes\Q_{\ell} \hra H^2_{et}({X}_{n-2}\times k^{sep},
\Q_{\ell}(1)).$$or  

{\rm (ii)} the Tate conjecture {\rm (\ref{tatec})} is true for 
either $X_{n-3}$ or $X_{n-2}$.

\nd In particular, $\nu_{et}$ is injective when $n \le 2$. 
\end{lem} 
\begin{proof} It is straightforward to check that (i) implies the
desired injectivity. Let us show that (ii) implies (i). The 
exact sequences of $\G_{sep}$-modules$$ 0 \to \Ker(\delta^*_{n-3}) \to  
 H^2_{et}({X}_{n-3}\times k^{sep},
\Q_{\ell}(1)) \to \delta^*_{n-3} H^2_{et}({X}_{n-3}\times k^{sep},
\Q_{\ell}(1)) \to 0,$$ 
$$0 \to  \delta^*_{n-3}H^2_{et}({X}_{n-3}\times k^{sep},
\Q_{\ell}(1)) \to  H^2_{et}({X}_{n-2}\times k^{sep},
\Q_{\ell}(1)) \to \Coker(\delta^*_{n-3}) \to 0$$are split: this
follows from \cite[2.10]{tate} --- see the proof of
\cite[5.2(b)]{tate}. Thus, as in (loc. cit), these exact sequences remain
exact after taking $\G_{sep}$-invariants. This, by (\ref{picmd}(v)), 
suffices for the
implication (ii) $\Rightarrow$ (i).\end{proof} 

 Since the map $\frac{K^0(S)}{W_{n-2}(S)}\otimes\Q_{\ell}  
 \to \mathcal B_n({S})\otimes\Q_{\ell} \xrightarrow{\nu_{et}}  
Gr^W_{0}H^n_{et}(\bar{U}_{\bu}, \Q_{\ell}(1))$ is clearly injective,
 we obtain the injectivity of \begin{equation}\label{lambs}
\Lambda_{et}: T_{\ell}J^n(X_{\bu},
 Y_{\bu})\otimes\Q \hra
 H^n_{et}(\bar{U}_{\bu}, \Q_{\ell}(1));\end{equation} its image is denoted
 $s^n_{\ell}(X_{\bu}, Y_{\bu})$. A similar definition gives
 $s^n_{\A^p}(X_{\bu}, Y_{\bu})$.

\begin{rem} (i) Without a condition such as the injectivity of $\nu_{et}$,
 it is hard to show that the Galois representations
 $\Lambda_{et}(T_{\ell}L^n\otimes\Q)$ are ``independent of $\ell$'';
 this is unclear already for their dimensions \cite[pp.27-29]{ka}. 

(ii) If $n \le 2$, then $\nu_{et}$ is injective; set
 $t^n_{\ell}(X_{\bu}, Y_{\bu}):= \im(\Lambda_{et})$, a $\G$-invariant
 subspace of $H^n_{et}(\bar{U}_{\bu}, \Q_{\ell}(1))$. A similar
 definition gives $t^n_{\A^p}(X_{\bu}, 
 Y_{\bu})$.\qed\end{rem} 

\nd {\bf Varieties over perfect fields.}\ssp

 From now on, we assume that $k$ is
perfect.  Let $V$ be an arbitrary variety over $S$.  
The results of \cite{dj} show \cite[\S 1]{dj} \cite[6.3]{be} 
 the existence of a simplicial pair
$(X_{\bu}, Y_{\bu})$ over $S$ 
with a morphism $\alpha: U_{\bu} \ra
V$ which makes $U_{\bu}$  a proper hypercovering of
$V$. One has (\ref{sonaty})
$\alpha^*: H^*_{et}(\bar{V},{\Z}^p(1)){\zim}
H^*_{et}(\bar{U}_{\bu},{\Z}^p(1))$; as before, $\beta_U$ denotes
the inverse of $\alpha^*$. The methods of \cite[6.2]{h} imply, using
\cite{dj},  the
abundance of smooth 
hypercoverings so that (\ref{bertie}) is valid as well.  

\begin{defn}  
Let\footnote{This follows a suggestion of M. Marcolli.} 
$s^n_{\ell}(V)$   be the direct limit of
$\beta_U(s^n_{\ell}(X_{\bu}, Y_{\bu}))$
over all proper hypercoverings 
$\alpha: U_{\bu} \rightarrow V$. It is a
$\Q_{\ell}$-subspace of 
$H^n_{et}(\bar{V},\Q_{\ell}(1))$, with an action of $\G$. A similar
definition holds for $t^n_{\ell}(V)$ ($n \le 2$). 
\end{defn}

Since $H^n_{et}(\bar{V},\Q_{\ell}(1))$ \cite[p. 24]{ka} is a finite
dimensional vector space 
over $\Q_{\ell}$, there is actually a (by no means unique)
proper hypercovering $U_{\bu}$ of $V$ such that
$\beta_U(s^n_{\ell}(X_{\bu}, Y_{\bu})) = s^n_{\ell}(V)$; let us
call such hypercoverings ``neat''.

\begin{prop} The notion of ``neat'' does not
depend on the auxiliary prime $\ell \neq p$. 
\end{prop}

\begin{proof} Clear: the map $\Lambda_{et}$ of (\ref{lambs}) is
injective which means that the dimension of $s^n_{\ell}(X_{\bu},
Y_{\bu})$ is independent of $\ell$; similarly, the dimension of 
$s^n_{\ell}(V)$ is independent of $\ell$.  \end{proof} 

\begin{lem}\label{special} {\rm (i)} Any proper hypercovering of $V$
which dominates a ``neat'' proper hypercovering of $V$ is ``neat''.

{\rm (ii)} Any two ``neat'' proper hypercoverings are dominated by a
``neat'' proper hypercovering. 

{\rm (iii)} More generally, given a pair of
proper hypercoverings one of which is ``neat'', there is a proper
hypercovering (automatically ``neat'') dominating them. 

{\rm (iv)} Any proper
hypercovering is dominated by a ``neat'' proper hypercovering. \qed
\end{lem}

 The advantage of ``neat'' hypercoverings is the

\begin{lem} For any morphism $\theta$ between two ``neat'' proper
hypercoverings $U_{\bu}$ and $'U_{\bu}$ of $V$, 
the induced map  $\theta^*: J^n('X_{\bu}, 'Y_{\bu})\otimes\Q \to
J^n(X_{\bu}, Y_{\bu})\otimes\Q$ is an
isomorphism. 

Here $J^n(-)\otimes\Q$ denotes the
isogeny one-motive obtained from the generalized one-motive --- see {\rm
(\ref{rtp})}.       
\end{lem}

\begin{proof} It is enough to show that the map $   
\theta^*: T_{\ell}J^n('X_{\bu}, 'Y_{\bu})\otimes \Q \ra
T_{\ell}J^n(X_{\bu}, Y_{\bu})\otimes \Q$
is an isomorphism. This map is always injective (for
arbitrary proper hypercoverings). But if the hypercoverings are
``neat'', then both terms are actually equal to
$s^n_{\ell}(V)$.\end{proof}

Our definition of ``neat'' depends upon the integer $n$. But since
$H^i_{et}(\bar{V}, \A^p(1)) = 0$ for $i > 2$~dim~$V$ \cite[pp. 23-24]{ka},
there exist, by  
(\ref{special}), ``neat'' hypercoverings such that
$\beta_U(s^n_{\ell}(X_{\bu}, Y_{\bu})) = s^n_{\ell}(V)$ holds for all
$n$. From now on, let us consider only such ``neat'' hypercoverings. 
\begin{defn} We define $J^n(V)\otimes\Q$ to be the isogeny one-motive
$J^n(X_{\bu}, Y_{\bu})\otimes\Q$  
of any ``neat'' proper hypercovering $U_{\bu}$
of $V$. \end{defn}

The following theorem is a trivial consequence of the previous
definition; the weight filtration on \'etale cohomology is given by 
\cite[5.3.6]{weil2}.
\begin{thm} One has $\G$-equivariant injections
$$\Lambda_{et}: T_{\ell} W_iJ^n(V)\otimes\Q \hra W_iH^n_{et}(\bar{V},
\Q_{\ell}(1));$$ these are isomorphisms for $i = -2, -1$.
Similar statement holds for $\Lambda_{et}:T^pJ^n(V)\otimes\Q \hra
H^n_{et}(\bar{V},\A^p(1))$. \qed    
\end{thm}  

\begin{rem} Note $s^n_{\ell}(V)$ is a
  subspace (possibly 
proper) of $H^n_{et}(\bar{V},\Q_{\ell}(1))$; the allowed weights on 
 $s^n_{\ell}(V)$ are $-2$, $-1$ and $0$; the allowed weights 
on  $H^n_{et}(\bar{V},\Q_{\ell}(1))$  lie between $-2$ and
$2n-2$.  A similar statement is true for $s^n_{\A^p}(V):=
\Lambda_{et}(T^pJ^n(V)\otimes\Q) \subset H^n_{et}(\bar{V},\A^p(1))$.\qed
\end{rem} 

\begin{rem} (Functoriality) Given any morphism $f: Z \ra V$, there is
  an induced morphism $f^*: 
J^n(V)\otimes\Q \ra J^n(Z)\otimes\Q$:  given $f$, pick a ``neat''
proper hypercovering $U_{\bu}$ of $V$. One can find a proper
hypercovering $E_{\bu}$ of $Z$ and a morphism
$F: E_{\bu} \ra U_{\bu}$ lifting $f$. By (\ref{special}),
``neat'' proper hypercoverings of $Z$ are cofinal among proper
hypercoverings of $Z$. So we may choose $E_{\bu}$ to be
``neat''. 
This yields the functoriality of $J^n(-)\otimes\Q$, 
$s^n_{\ell}(-)$, and $s^n_{\A^p}(-)$.\qed \end{rem}

\begin{rem}\label{mem} (i) A variant of (\ref{intero}) shows that the $\G$-invariant 
$\Z^p$-lattice $$(s^n_{\A^p}(V) \cap H^n_{et}(\bar{V},
\Z^p(1))/{\rm torsion}) \subset H^n_{et}(\bar{V},\A^p(1)) $$ provides
a  $\Z[1/p]$-integral structure on $J^n(V)\otimes\Q$, i.e., determines
$J^n(V)\otimes \Z[1/p]$, a one-motive up to $p$-isogeny, defined over
$k$.  Controlling
$p$ would require an integral $p$-adic cohomology for arbitrary
varieties. 

(ii) It is probable that 
the $p$-adic/crystalline realization of
$J^n(V)\otimes\Q$ is related to the rigid cohomology of $V$
\cite{tsu}.\qed\end{rem}   

\nd {\bf The one-motives up to $p$-isogeny
$L^n(V)\otimes\Z[1/p]$.}\smallskip  

Our methods for the construction of $J^n(V)\otimes\Z[1/p]$ show the
following 

\begin{thm}\label{groth} If {\rm (\ref{tatec})}
is true for any
surface over any finitely generated {\rm (}over $\F_p${\rm)} subfield of a
perfect field $k$, 
then {\rm (\ref{dc1}) (up to $p$-isogeny)} is
true for $k$.\end{thm}

\begin{proof} Using the injectivity of  $\nu_{et}$ of (\ref{limbs})
assured by (\ref{maam}), we can apply the previous methods to define
$L^n(V)\otimes\Q$ (using a variant of ``neat'' hypercoverings) and
refine it, as in (\ref{mem}(i)), to  $L^n(V)\otimes\Z[1/p]$. The fact
that the Tate conjecture (\ref{tatec}) for divisors 
reduces to the case of surfaces is well-known; this is proved along
the ideas of the proof of (\ref{maam}) using the weak Lefschetz
theorem \cite[4.1.6]{weil2}.\end{proof}  

\begin{rem}\label{groth2} (i) By (\ref{maam}),  $\nu_{et}$ of (\ref{limbs}) is injective
for $n \le 2$. In this case, our methods provide $L^n(V)\otimes
\Z[1/p]$ (for $n \le 2$), 
one-motives up to $p$-isogeny, associated with $V$ which are
contravariant functorial. It follows from the definitions
(\ref{phin}), (\ref{hosadu}) that $J^n(V)\otimes \Z[1/p] = L^n(V)\otimes
\Z[1/p]$ for $n \le 1$.

(ii) If the field $k$ in (\ref{groth}) is not taken to be perfect,
 then one obtains that $L^n(V)\otimes\Z[1/p]$ (attached to a
 variety $V$ over $k$) is defined over $k^{perf}$.\qed\end{rem}

\begin{thm}\label{imper} Let  $F$ be a field of
characteristic $p >0$; write $F^{perf}$ for its perfection. Assume
that ``resolution of singularities''
holds over $F$ and that 
{\rm (\ref{tatec})} holds for any surface over any finitely generated
{\rm(}over $\F_p${\rm)} subfield of $F$. Then, {\rm (\ref{dc1}) (up to
$p$-isogeny)} is
true for $F$.
\end{thm} 

\begin{proof} Straightforward variant of (\ref{groth}).\end{proof} 

\section{Applications and related results}\label{applix}

\nd {\bf Related work.}\ssp 

For any curve $C$, Deligne \cite[10.3]{h}
constructed a one-motive $H^1_m(C)(1)$ (isomorphic to our $L^1(C)$)
and used it to prove Theorem \ref{main} for the 
$H^1$ of a curve over $\C$.
The one-motive $L^2(V)$ of a projective complex surface
$V$ was already obtained by J. Carlson \cite{ca}. Carlson mentions in \cite{ca2} that
he has constructed other one-motives for a special class of varieties
(over $\C$)  
but these results remain 
unpublished (email, 3 Dec 2001). The one-motive $L^1(X_{\bu}, Y_{\bu})$
is the Picard 
one-motive  $Pic^{+}$ of \cite{bs3} 
and $M^1$ of \cite{ra}. Finally, \cite{brs} contains independent
proofs of some of our
results.\ssp

\nd {\bf Motivic principles: illustrations.}
\cite[1.7, 2.5]{ja}\ssp 

Let $V$ be a variety over a field $k$ of characteristic zero; by
(\ref{zeusp}), we may assume $k$ to be finitely generated over $\Q$
without loss of generality.
\begin{prop}\label{swaraj} The rank of
$$ H^{1,1}_{\Q}(V_{\iota})^n:= \Hom_{MHS}(\Q(0), Gr^W_0
H^n(V_{\iota},\Q(1))),$$  
{\rm (the so-called $(1,1)$-part)} is independent of the imbedding $\iota: k
\hookrightarrow \C$. \end{prop}   

\begin{proof} This follows from Theorem \ref{main} since the
dimension of the left hand side is the rank of 
$\mathcal B_n(\C)$ of the one-motive $L^n(V)$.\end{proof} 

\begin{rem} Put $k =\C$; (\ref{swaraj}) 
 yields an algebraic characterization of the
$(1,1)$-part of $H^n(V,\Q)$ (\ref{1-1i}). The
analogous statement for the $(m,m)$-part of $H^n(V,\Q)$ is not known
(for $m >1$) in general; it is part of the (homological) Hodge conjecture
\cite[7.2]{jamm}.\qed \end{rem}     

Consider
the Galois representations $M_{\ell}:=H^n_{et}(\bar{V}, \Q_{\ell}(1))$
(these have a weight filtration $W$ by Galois submodules 
\cite[\S 13, \S14]{hp}.  For each prime $\ell$, $h_{\ell}(V):=
W_{-1}M_{\ell}$ defines an 
element $\zeta_{\ell} \in
Ext^1_{\G}(Gr^W_{-1}M_{\ell},Gr^W_{-2}M_{\ell})$. For each $\iota: k
\hra \C$, the substructure $h_{\iota}(V):= W_{-1}M_{\iota}(V)$ of
 the mixed Hodge structure $M_{\iota}:=
H^n(V_{\iota},\Q(1))$ defines an element $\zeta_{\iota}
\in Ext^1_{MHS}(Gr^W_{-1}M_{\iota},Gr^W_{-2}M_{\iota})$.  

The isogeny one-motive $W_{-1}L^n\otimes\Q$ can be viewed as 
extension $\zeta$ of two isogeny one-motives given by the
exact sequence (\ref{lcw}): 
$$\zeta: 0 \rightarrow [0 \rightarrow \mathcal T]\otimes\Q \rightarrow [0
\rightarrow \tilde{\mathcal P}_n]\otimes\Q \rightarrow
 [ 0 \rightarrow \mathcal R]\otimes\Q \rightarrow 0$$

The element $\zeta \in Ext^1(\mathcal R, \mathcal T)\otimes\Q$ is zero
if and only if $\tilde{\mathcal P}_n$ is isogenous to the
direct product $\mathcal R \times \mathcal Q$. 

\begin{cor}\label{jap} One has

{\rm (a)} $\zeta$ is zero $\Leftrightarrow  \zeta_{\iota}$ is zero for all
$\iota  \Leftrightarrow  \zeta_{\iota}$ is zero for one
$\iota:k\hookrightarrow \C$.

{\rm (b)} $\zeta$ is zero 
$\Leftrightarrow \zeta_{\ell}$ is zero  for all primes $\ell
\Leftrightarrow \zeta_{\ell}$ is zero for one $\ell$.
\end{cor}

\begin{proof} The relation between $\zeta$, $\zeta_{\ell}$, and
$\zeta_{\iota}$ is given by (\ref{biget}) and (\ref{kona}). Statements
(a) and (b)  follow easily from 
\cite[Thm. 4.3]{ja}  and  \cite[10.1.3]{h} that the
$\ell$-adic (or Hodge) realizations of an isogeny one-motive is a direct
sum of pure objects if and 
only if the isogeny one-motive is isogenous to a direct
product.\end{proof}\ssp

\nd {\bf Independence of the prime $\ell$ in \'etale cohomology.}\ssp 

 We now take $k$ to 
be a finite field; and let $\Phi \in \G$ be the Frobenius
automorphism.  

The weight filtration $W$ on $H^n_{et}(\bar{V},\Q_{\ell})$
 is defined via the endomorphism  $\Phi^*_{\ell}$ of
 $H^n_{et}(\bar{V},\Q_{\ell})$ induced by $\Phi$.   
Let $b_{i,\ell,n}(V)$ be the dimension of the 
  $\Q_{\ell}$-vector space $Gr^W_{i}H^n_{et}(\bar{V},\Q_{\ell})$; we
recall the well-known results:
\begin{lem}\label{wayout} {\rm (Deligne)} \cite{weil} \cite[3.3.9]{weil2}  If $V$ is smooth and
proper, then   

{\rm (1)} $b_{i,\ell,n}(V) =0$ if $i \neq n$. 

{\rm (2)} $b_{n,\ell,n}(V)$ is independent of $\ell$; thus, we may set
$b_n(V) = b_{n,\ell,n}(V)$. 
 
{\rm (3)}  the characteristic polynomial of $\Phi^*_{\ell}$ on
$H^n_{et}(\bar{V},\Q_{\ell})$  has coefficients in $\Q$ independent of
the prime $\ell$.

{\rm (4)} \cite[4.1.5]{weil2} {\rm (}$V$ projective{\rm)}
$b_{2n+1}(V)$ is even.\footnote{Deligne (loc. cit) remarks that this is not
yet known for $V$ proper smooth.}\qed\end{lem}

\begin{quest} {\rm (N. Katz)} \cite[pp. 27-29]{ka} Which parts of {\rm
(\ref{wayout})} are valid for 
$b_{i,\ell,n}(V)$, $Gr^W_{i}H^n_{et}(\bar{V},\Q_{\ell})$ for general 
$V$, i.e., for $V$ possibly singular and non-proper?
\end{quest} 

This is not answered by a formal application of \cite{weil} and
resolution of singularities;
cf. \cite[\S1]{dj}. Genuinely new ingredients
are needed for an answer in general. 
Katz \cite{ka3} has proved it for
 certain  smooth varieties; no general results are known
 for singular varieties. A consequence of (\ref{zund}(v)) is: 

\begin{thm}
Fix an arbitrary variety $V$ and an integer $n \ge 0$.  The following systems  of
Galois representations  
satisfy {\rm (2)} and {\rm (3)} of {\rm (\ref{wayout})}:  
  
\nd {\rm (i)} $W_{-2}H^n_{et}(\bar{V},\Q_{\ell}(1))$. {\rm (ii)}
$Gr^W_{-1}H^n_{et}(\bar{V},\Q_{\ell}(1))$. {\rm (iii)}
$W_{-1}H^n_{et}(\bar{V},\Q_{\ell}(1))$.\qed    
\end{thm}
Let  $c_n(V)$ denote the rank (over $\Z$) of the Galois module
  $I_n(\bar{S})$; here $I_n$ is the ``lattice'' in the one-motive
  $J^n(V) = [I_n \ra 
  \mathcal P_n]$.

\begin{thm}\label{zund}

{\rm (i)} The integers $b_{-2,\ell ,n}(V)$ and $b_{-1,\ell ,n}(V)$ are
  independent of $\ell\neq p$.

{\rm (ii)}  $b_{-1, \ell ,n}(V)$ is an \ub{even} integer;
cf. \cite[16.1]{hp}.   

{\rm (iii)} $b_{0, \ell , n(V)} \ge c_n(V)$.

{\rm (iv)} Let $f:V \ra V$ be any morphism. The characteristic
  polynomial of the induced endomorphism $f^*_{\ell}: s^n_{\ell}(V)
  \ra s^n_{\ell}(V)$  has $\Q$-coefficients which are independent of $\ell$.

{\rm (v)} The characteristic polynomial of the automorphism $\Phi_{\ell}$ on
the $\G$-submodule $s^n_{\ell}(V)$ of
$H^n_{et}(\bar{V},\Q_{\ell}(1))$ has $\Q$-coefficients which
are independent of $\ell$.

{\rm (vi)} For any
  morphism $g: V \ra W$, we have a commutative diagram  
$$
\begin{CD} 
s^n_{\ell}(W) @>{g^*_{\ell}}>> s^n_{\ell}(V) \\
@VVV @VVV\\
H^n_{et}(\bar{W},\Q_{\ell}(1)) @>>> H^n_{et}(\bar{V},\Q_{\ell}(1));\\
\end{CD}
$$
the characteristic polynomials of
  $\Phi_{\ell}$ on the $\G$-modules $\Ker(g^*_{\ell})$ and $\Coker(g^*_{\ell})$
  have $\Q$-coefficients which are independent of $\ell$.
\end{thm}

\begin{proof} Write $\mathcal P$ as an extension of an abelian
  variety $\mathcal A$ by a torus $\mathcal T$; it
  follows from (\ref{biget}) that $b_{-2,\ell ,n}(V)$ (resp.
  $b_{-1,\ell ,n}(V)$) are the dimensions of $\mathcal T$
  (resp. $\mathcal A$). This proves (i) and (ii). (iii) also follows
  from (\ref{biget}) since we know (\ref{lambs}) that
  $I(\bar{S})\otimes_{\Z}\Q_{\ell}$  injects into $Gr^W_{0}H^n_{et}(\bar{V},
  \Q_{\ell}(1))$. 

The endomorphism
  algebra $\End(J^n(V)\otimes\Q)$ is a finite dimensional $\Q$-algebra. By
  functoriality of $J^n(-)\otimes\Q$, the morphism $f$ induces an
  element $f^* \in \End(J^n(V)\otimes\Q)$. The characteristic
  polynomial of $f^*$ is a polynomial with $\Q$-coefficients. Since
  $s^n_{\ell}(V)$ is the $\ell$-adic realization of $J^n(V)\otimes\Q$,
  by functoriality, the map $f^*_{\ell}$  on $s^n_{\ell}(V)$ has the
  same characteristic 
  polynomial. This proves (iv). The same argument proves (v): the
  Frobenius morphism $F_V$ of $V$ and the geometric Frobenius
  $\Phi^{-1}$ induce the same endomorphism on
  $H^*(\bar{V},\Q_{\ell}(1))$ \cite[1.15]{weil}.         

The map $g: V \to W$ induces a map $g^*: J^n(W)\otimes\Q \to
J^n(V)\otimes\Q$. Since the category of isogeny one-motives is abelian,
we have the isogeny one-motives $\Ker(g^*)$ and $\Coker(g^*)$. Their
$\ell$-adic realizations are the  Galois
modules $\Ker(g^*_{\ell})$ and $\Coker(g^*_{\ell})$. Apply the argument
in the previous paragraph to $\End(\Ker(g^*))$ and
$\End(\Coker(g^*))$. This proves (vi). \end{proof}\ssp

\nd {\bf Rationality of systems of certain $\ell$-adic Galois
  representations.}\ssp

\begin{quest} {\rm (J.-P. Serre)} \cite[I-10]{jp} \cite[12.5?]{se2}
For a fixed variety $V$ over a number field $k$ and integer 
 $n$, is  the system of Galois  
representations $H^n_{et}(\bar{V},\Q_{\ell})$ ``rational''?
\end{quest}  

The same question for the systems
$W_iH^n_{et}(\bar{V},\Q_{\ell}(1))$ clearly refines the above one. For
smooth proper $V$,  Deligne's theorem (Weil conjectures) 
\cite{weil2} provides an affirmative 
answer (see \cite{se2} Exemple after 12.5? on page 393). 
 For any imbedding $\iota: k \hookrightarrow \C$,
the weight filtrations on $H^n(V_{\iota}, \Q(1))$ and 
$H^n_{et}(\bar{V},\Q_{\ell}(1))$ are compatible
\cite[\S14]{hp}. But this does not imply the  
``rationality'' of the system $W_iH^n_{et}(\bar{V},\Q_{\ell}(1))$ in general.
However, the combination of (loc. cit) and the proof of 
(\ref{zund} (v)) yields  
\begin{thm} The system of Galois representations $W_j
  H^n_{et}(\bar{V},\Q_{\ell}(1))$ is ``rational'' for $j =-2,-1$ as is
  the system $t^n_{\ell}(V)$. \qed \end{thm}

\nd \emph{Acknowledgements.} I heartily thank S. Lichtenbaum for his
constant encouragement, guidance, and support. I would like to express 
my gratitude  to S. Bloch, P. Deligne, H. Gillet, J. Gordon,
M. Marcolli, J. Milne, M. Nori, and
B. Totaro. The debt to the work of Deligne and Carlson should be
evident. I would like to thank the
referees for many constructive remarks and suggestions.\medskip

\hfill ``iyaM visRSTiryata AbabhUva yadi vA dadhe yadi vA na

\hfill yo asyAdhyakSaH parame vyoman so aNga veda yadi vA naveda''

\hfill Nasadiya Sukta (Rigveda X 129).

%

\providecommand{\bysame}{\leavevmode\hbox to3em{\hrulefill}\thinspace}
\providecommand{\MR}{\relax\ifhmode\unskip\space\fi MR }
\providecommand{\MRhref}[2]{%
  \href{http://www.ams.org/mathscinet-getitem?mr=#1}{#2}
}
\providecommand{\href}[2]{#2}
\bibliographystyle{amsplain}

\end{document}